\documentclass[preprint,3p]{elsarticle}
\makeatletter
\def\ps@pprintTitle{%
\let\@oddhead\@empty
\let\@evenhead\@empty
\def\@oddfoot{\reset@font\hfil\thepage\hfil}
\let\@evenfoot\@oddfoot
}
\makeatother

\usepackage[utf8]{inputenc}

\usepackage{amsmath, amsfonts, amsthm, algorithm, algorithmicx, algpseudocode, booktabs}
\newtheorem{remark}{Remark}
\usepackage{bm}
\usepackage{caption}
\usepackage{subcaption}
\usepackage{siunitx}
\DeclareSIUnit\inch{in}
\DeclareSIUnit\foot{ft}
\DeclareSIUnit\hour{hr}
\DeclareSIUnit\minute{min}

\DeclareMathOperator*{\argmin}{arg\,min\,}
\let\max\relax
\let\min\relax
\DeclareMathOperator*{\max}{max\,}
\DeclareMathOperator*{\min}{min\,}
\DeclareMathOperator*{\mean}{mean\,}
\providecommand{\abs}[1]{\lvert#1\rvert}

\providecommand{\norm}[1]{\lVert#1\rVert}

\begin{document}

\begin{frontmatter}

\title{Mesh sampling and weighting for the hyperreduction of nonlinear Petrov-Galerkin reduced-order models with local reduced-order bases\protect}

\author[1]{Sebastian Grimberg\corref{cor1}}
\cortext[cor1]{Corresponding author.}
\ead{sjg@stanford.edu}
\author[1,2,3]{Charbel Farhat}
\ead{cfarhat@stanford.edu}
\author[1]{Radek Tezaur}
\ead{rtezaur@stanford.edu}
\author[4]{Charbel Bou-Mosleh}
\ead{cboumosleh@ndu.edu.lb}

%\authormark{GRIMBERG \textsc{et al.}}

\address[1]{Department of Aeronautics and Astronautics, Stanford University, Stanford, CA 94305, U.S.A.}
\address[2]{Department of Mechanical Engineering, Stanford University, Stanford, CA 94305, U.S.A.}
\address[3]{Institute for Computational and Mathematical Engineering, Stanford University, Stanford, CA 94305, U.S.A.}
\address[4]{Department of Mechanical Engineering, Notre Dame University-Louaize, Zouk Mosbeh, Lebanon}

%\address[1]{\orgdiv{Department of Aeronautics and Astronautics}, \orgname{Stanford University},
%\orgaddress{Durand Building, 496 Lomita Mall, Stanford, \state{CA} 94305-4035, \country{USA}}}
%\address[2]{\orgdiv{Department of Mechanical Engineering}, \orgname{Stanford University},
%\orgaddress{Building 530, 440 Escondido Mall, Stanford, \state{CA} 94305-3030, \country{USA}}}
%\address[3]{\orgdiv{Institute for Computational and Mathematical Engineering}, \orgname{Stanford University},
%\orgaddress{Huang Engineering Center, 475 Via Ortega, Suite 060, Stanford, \state{CA} 94305-4042, \country{USA}}}
%\address[4]{\orgdiv{Department of Mechanical Engineering}, \orgname{Notre Dame University-Louaize},
%\orgaddress{P.O. Box 72, Zouk Mosbeh, \country{Lebanon}}}

%\corres{*Sebastian Grimberg, Department of Aeronautics and Astronautics, Stanford University, Durand Building, 496 Lomita Mall, Stanford, %\state{CA}, 94305-4035, \country{USA}\\ \email{sjg@stanford.edu}}

%\presentaddress{This is sample for present address text this is sample for present address text}

\begin{abstract}
\small
The energy-conserving sampling and weighting (ECSW) method is a hyperreduction method originally developed for accelerating the performance of Galerkin
projection-based reduced-order models (PROMs) associated with large-scale finite element models, when the underlying projected operators need to be frequently recomputed as in parametric and/or
nonlinear problems. In this paper, this hyperreduction method is extended to Petrov-Galerkin PROMs
where the underlying high-dimensional models can be associated with arbitrary finite element, finite volume, and finite difference semi-discretization methods. Its scope is also extended to cover
local PROMs based on piecewise-affine approximation subspaces, such as those designed for mitigating the Kolmogorov $n$-width barrier issue associated with convection-dominated
flow problems. The resulting ECSW method is shown in this paper to be robust and accurate. In particular, its offline phase is shown to be fast and parallelizable, and the potential of its online phase
for large-scale applications of industrial relevance is demonstrated for turbulent flow problems with $O(10^7)$ and $O(10^8)$ degrees of freedom. For such problems, the online part of the ECSW method
proposed in this paper for Petrov-Galerkin PROMs is shown to enable wall-clock time and CPU time speedup factors of several orders of magnitude while delivering
exceptional accuracy.
\end{abstract}

%\keywords{hyperreduction, local basis, machine learning, nonlinear model reduction, Petrov-Galerkin, reduced mesh}

\end{frontmatter}

%\maketitle

%\tableofcontents

%

\section{Introduction}
\label{sec:intro}

Projection-based model order reduction (PMOR) is an enabling computational technology for dramatically reducing the solution time and storage requirements associated with the analysis of 
high-fidelity, high-dimensional (large-scale), differential-equation-based computational models in various parametric settings. Time-critical applications requiring real-time, near real-time, or simply 
faster than ever performance such as computational-based design and optimization, statistical analysis, and optimal control stand to benefit from the reduced computational complexity of the outcome of 
PMOR -- that is, from projection-based reduced-order models (PROMs). Unfortunately, the computational cost associated with the construction of a PROM of dimension $n$ -- and more specifically, the 
underlying projected vectors and matrices (a.k.a. the reduced-order vectors and matrices) -- typically scales with both $n$ and the dimension of the underlying high-dimensional computational model (HDM) 
$N \gg n$. Thus, unless this issue is mitigated, a typical PMOR method does not necessarily deliver the expected speedup factor.

A popular approach for addressing the aforementioned scaling issue is to decompose the computation of projected vectors and matrices in two parts: one whose computational complexity scales with the
dimension of the HDM but can be executed {\it offline} to pre-compute once-for-all some numerical quantities -- that is, before the parametric or repetitive PROM computations are performed; and another 
part which exploits the aforementioned pre-computed quantities to perform all {\it online} computations with a complexity that is independent of the large dimension of the HDM. This offline-online 
decomposition has been extensively demonstrated for parametric, linear problems \cite{antoulas2005, benner2015} and for nonlinear problems exhibiting a low-order polynomial dependence on the chosen 
degrees of freedom (DOFs) \cite{rowley2004, barbic2005, nguyen2005}. For some HDMs with non-polynomial nonlinearities, lifting transformations have been advocated to enable the offline-online 
decomposition approach associated with low-order polynomial nonlinearities \cite{kramer2019}. However, for highly nonlinear HDMs, such transformations are in general computationally intractable, in 
which case the aforementioned offline-online decomposition is simply not an option.

Alternatively, a computational approach known as hyperreduction has gained wide acceptance for eliminating the computational bottlenecks associated with the repeated re-evaluations of parametric 
reduced-order vectors and matrices. A comprehensive review of this approach, which is equally applicable to linear and arbitrarily nonlinear parametric PROMs, can be found in \cite{farhat2020}. Generally 
speaking, hyperreduction introduces in the construction of a PROM an additional approximation layer that enables the efficient computation of the projected, low-dimensional vectors and matrices defining
a PROM without sacrificing the desired level of accuracy. Specifically, the computational complexity of the processing of a PROM after it is hyperreduced is independent of the size of the underlying
HDM.

Hyperreduction methods can be classified in two types \cite{farhat2020}: {\it approximate-then-project} methods, which were developed first and have a longer record of successful applications; and 
{\it project-then-approximate} methods, which have been proposed more recently and are more robust. The origins of approximate-then-project hyperreduction methods can be traced back to the gappy proper 
orthogonal decomposition (POD) method \cite{everson1995} developed in the context of image reconstruction. As their label implies, such methods approximate first a nonlinear high-dimensional 
quantity -- for example, by interpolating information at a carefully selected subset of the nodes of a computational mesh using a small number of empirically-derived basis functions -- then 
compute the exact projection of the approximation onto the left reduced-order basis (ROB) underlying the PROM of interest. The first examples of this type of hyperreduction methods include the 
empirical interpolation method (EIM) \cite{barrault2004, grepl2007} and the related best points interpolation method \cite{nguyen2008}, which are derived at the continuous level and enjoy some level 
of theoretical support for elliptic problems. Variant approaches based on function sampling have also been proposed, including the missing point estimation method \cite{astrid2008} and the collocation 
method \cite{ryckelynck2005, ryckelynck2009}. In this first type of hyperreduction methods however, the discrete empirical interpolation method (DEIM) \cite{chaturantabut2010}, which can be viewed 
as a variant of the discrete form of the aforementioned EIM, is certainly the most popular method: it has seen widespread adoption along with subsequent adaptations \cite{tiso2013, peherstorfer2014}.

Hyperreduction methods of the project-then-approximate type differ from their predecessors by directly approximating the reduced-order vectors and matrices associated with the projection of HDM quantities
onto the left ROB associated with the PROM of interest. As such, they seek to achieve a better approximation of the PROM than approximate-then-project hyperreduction methods. They can be interpreted as 
generalized quadrature rules in which the set of quadrature ``points'' and associated weights are learned in a supervised procedure on an empirical set of training data, which illustrates one of the 
many connections between PMOR and machine learning. The earliest hyperreduction methods of this type include the cubature-based approach of \cite{an2008} and the energy-conserving sampling and weighting 
(ECSW) method \cite{farhat2014}. More recent examples include the empirical cubature method (ECM) \cite{hernandez2017} and the linear program-based empirical quadrature method \cite{yano2019}. All such 
methods compute a reduced mesh -- that is, a subset of the elements or other entities of the computational mesh underlying the given HDM -- whose entities define the quadrature points. In this sense, all 
hyperreduction methods of the project-then-approximate type perform {\it mesh sampling}. Among such methods, ECSW distinguishes itself from alternatives in that for second-order hyperbolic problems such 
as structural dynamics and wave propagation problems, it preserves the Lagrangian structure associated with Hamilton's principle \cite{farhat2015}. As such, if a time-integrator applied to a second-order 
hyperbolic PROM is provably unconditionally stable, it is guaranteed that this time-integrator will remain numerically stable on the hyperreduced counterpart PROM (HPROM) produced by ECSW. In particular,
it has been shown that ECSW leads to numerically stable and accurate HPROMs when applied to realistic structural dynamics problems for which state-of-the-art hyperreduction methods of 
the approximate-then-project type fail to do so \cite{farhat2015}. 

All hyperreduction methods of both types outlined above have been developed in the context of Galerkin PMOR methods -- that is, PMOR methods where the left and right ROBs are identical,
or equivalently, the subspaces of the test and approximation functions are the same. Only two of them have been explored for the acceleration of Petrov-Galerkin (PG) PMOR methods -- that is, PMOR
methods where the left and right ROBs differ. Both of these methods are of the approximate-then-project type: a gappy-POD-like method that is similar to DEIM, was developed for the Gauss-Newton with 
approximated tensors (GNAT) PMOR method \cite{carlberg2011, carlberg2013}, and is referred to in the remainder of this paper as the the gappy-POD-based GNAT method; and a least-squares variant of the 
collocation approach \cite{ryckelynck2005} that was first explored in \cite{legresley2006}; then tailored in \cite{washabaugh2016, zahr2016} for PG PROMs in computational fluid dynamics (CFD). 
Furthermore, hyperreduction methods of the project-than-approximate type have been developed so far only for second-order hyperbolic problems semi-discretized by a finite element (FE) method.

Yet, the PMOR of many parametric, nonlinear HDMs including those for which the tangent matrices are nonsymmetric calls for PG-based reduction methods. For example, it was recently
shown in \cite{grimberg2020} that the real culprit behind most if not all numerical instabilities reported in the literature for PROMs of convection-dominated laminar and turbulent flows is 
the Galerkin framework used for constructing such PROMs; and that alternatively, numerically stable and accurate PROMs for such laminar and turbulent flows can be constructed using a PG framework
without resorting to additional closure models or tailoring of the subspace of approximation. However, it was also shown in \cite{zahr2016} that PG HPROMs of convection-dominated 
flows constructed using the two aforementioned hyperreduction methods of the approximate-than-project type at best deliver a problem dependent performance and at worst perform poorly 
when applied outside of their training data. Furthermore, the mesh reduction algorithms underlying all hyperreduction methods of the approximate-then-project type are based on
suboptimal greedy sampling procedures that require as input the size of the reduced mesh -- which is unknown {\it a priori} -- thereby hindering practicality for many applications.

For all reasons summarized above, it remains to develop for nonlinear PG PROMs a reliable and practical hyperreduction method of the project-then-approximate type that is applicable to first-order 
hyperbolic problems semi-discretized by any preferred technique. This paper focuses on filling this gap. Specifically, it presents an ECSW-type method for the hyperreduction of PG PROMs where the 
underlying HDMs arise from the semi-discretization of first-order hyperbolic problems by any preferred scheme. It emphasizes CFD applications associated with convection-dominated viscous flows for 
the following reasons:
\begin{itemize}
	\item They are prime candidates for PG-based PMOR. 
	\item At high Reynolds numbers, they entail very large-scale meshes that are occasionally perceived to challenge the computational tractability of hyperreduction methods or their fidelity
		\cite{kramer2019} and therefore are excellent candidates for assessing their sheer wall-clock time performance and accuracy. 
	\item Due to the large Kolmogorov $n$-width \cite{pinkus1985} of the solution manifold associated with their HDMs, the PMOR of such applications typically requires the construction of multiple 
		local, piecewise-affine subspace approximations instead of a single global approximation \cite{peherstorfer2014, amsallem2012}, which raises the issue of how to perform hyperreduction
		in this case?
\end{itemize}

To this end, the remainder of this paper is organized as follows. In Section \ref{sec:pg}, the context is set to the PG-based PMOR of nonlinear, first-order dynamical systems using local subspace 
approximations. The computational bottlenecks are overviewed in order to motivate the concept of hyperreduction. In Section \ref{sec:ecsw}, the ECSW hyperreduction method previously developed for
Galerkin PROMs of HDMs associated with FE-based approximations of second-order dynamical systems is generalized to PG PROMs of HDMs constructed using any preferred semi-discretizations of first-order 
dynamical systems. In Section \ref{sec:sampling}, first the ECSW mesh sampling procedure is reviewed along with the parallel computational kernels required for minimizing the offline wall-clock time for 
very large-scale problems. Next, a comprehensive approach for constructing the reduced mesh resulting from the aforementioned sampling procedure is described. This approach covers many popular spatial 
discretization methods. In Section \ref{sec:app}, the generalized ECSW method is applied to the hyperreduction of three PG PROMs associated with convection-dominated,
viscous flow problems: an academic, unsteady, laminar flow problem that is easy to reproduce by the interested reader; an unsteady wake flow problem associated with the geometry of the so-called Ahmed 
body \cite{ahmed1984}; and a very large-scale turbulent flow problem with $O(10^8)$ unknowns associated with an F-16C/D Block 40 aircraft configuration with external stores at a high angle of attack.
The latter application is chosen to demonstrate both the computational tractability of the proposed ECSW hyperreduction method and its accuracy for a challenging, very large-scale application.
Performance comparisons between the proposed hyperreduction method, the gappy-POD-based GNAT alternative \cite{carlberg2011, carlberg2013}, and the least-squares collocation 
hyperreduction method \cite{legresley2006} equipped with the mesh sampling algorithm of the GNAT method are also presented for the first application. Finally, conclusions are offered in 
Section \ref{sec:conclusions}.

\section{Nonlinear Petrov-Galerkin projection-based model order reduction based on local reduced-order bases}
\label{sec:pg}

\subsection{Local subspace approximation and Petrov-Galerkin projection}
\label{sec:pg1}

Here, the focus is set on the first-order, $\bm{\mu}$-parametric, $N$-dimensional, highly nonlinear, semi-discrete problem
\begin{equation}
	\begin{split}
		\label{eqn:hdm0}
		\bm{M}(\bm{\mu})\dot{\bm{u}}(t; \bm{\mu}) + \bm{f}(\bm{u}(t; \bm{\mu}); \bm{\mu}) &= \bm{g}(t; \bm{\mu}) \\
		\bm{u}(0; \bm{\mu}) &= \bm{u}^0(\bm{\mu})
	\end{split}
\end{equation}
where $t \in [0, T_f]$ denotes time, the dot denotes the derivative with respect to time, $\bm{\mu} \in \mathcal{P}\subset\mathbb{R}^p$ is a $p$-dimensional vector of parameters, $\mathcal{P}$ denotes 
the bounded parameter space of interest, $\bm{u}(t; \bm{\mu})\in\mathbb{R}^N$ is the time-dependent solution vector, and $\bm{u}^0(\bm{\mu})\in\mathbb{R}^N$ is its initial condition. Throughout this 
paper, $\bm{M}(\bm{\mu}) \in \mathbb{R}^{N\times N}$ is a parametric mass matrix and is reasonably assumed to be symmetric positive definite (SPD), $\bm{f}(\bm{u}(t; \bm{\mu}); \bm{\mu})\in\mathbb{R}^N$
is a nonlinear function resulting from the semi-discretization of the partial differential equation (PDE) of interest and referred to as the nonlinear flux vector, and 
$\bm{g}(t; \bm{\mu})\in\mathbb{R}^N$ is a time-dependent, parametric, source term vector that may or may not be zero depending on the application. The HDM underlying problem (\ref{eqn:hdm0}) can 
represent the semi-discretization by any preferred approximation technique of any PDE -- given that higher-order ordinary differential equations can be rewritten in first-order form by expanding the 
vector of DOFs $\bm{u}(t; \bm{\mu})$. However, all numerical examples discussed in this paper pertain to first-order systems of conservation laws.

Following \cite{amsallem2012b}, the solution manifold associated with the above HDM-based problem -- which is assumed to be highly nonlinear -- is approximated here in the time and parameters domain of 
interest by a collection of $N_c$ piecewise-affine subspaces of $\mathbb{R}^N$. The dimension of each $k$-th affine subspace is denoted by $n_k$, with $n_k \ll N$ $\forall k\in\{1, \ldots, N_c\}$.  
The local approximation based on the $k$-th approximation subspace is written as
\begin{equation}
	\label{eqn:subspace}
	\bm{u}(t; \bm{\mu}) \approx \bm{u}_{o,k} + \bm{V}_k\bm{y}_k(t; \bm{\mu}) \quad \forall (t,\bm{\mu}) \in [0, T_f] \times \mathcal{P}
\end{equation}
where $\bm{u}_{o,k}\in\mathbb{R}^N$ is a fixed vector defining an affine offset of the $k$-th local subspace represented by the right ROB $\bm{V}_k \in \mathbb{R}^{N\times n_k}$
and $\bm{y}_k(t; \bm{\mu})\in\mathbb{R}^n$ is the vector of reduced (or generalized) coordinates of the representation of the solution in the $k$-th local approximation subspace. Throughout the
remainder of this paper, each right ROB $\bm{V}_k$ is assumed to be orthonormal -- that is, $\forall k$, $\bm{V}_k^T\bm{V}_{k}=\bm{I}$, where the superscript $T$ denotes the transpose operation.

\begin{figure}[h!]
	\centering
	%\vglue 0.1 truein
	\includegraphics[width=0.5\textwidth]{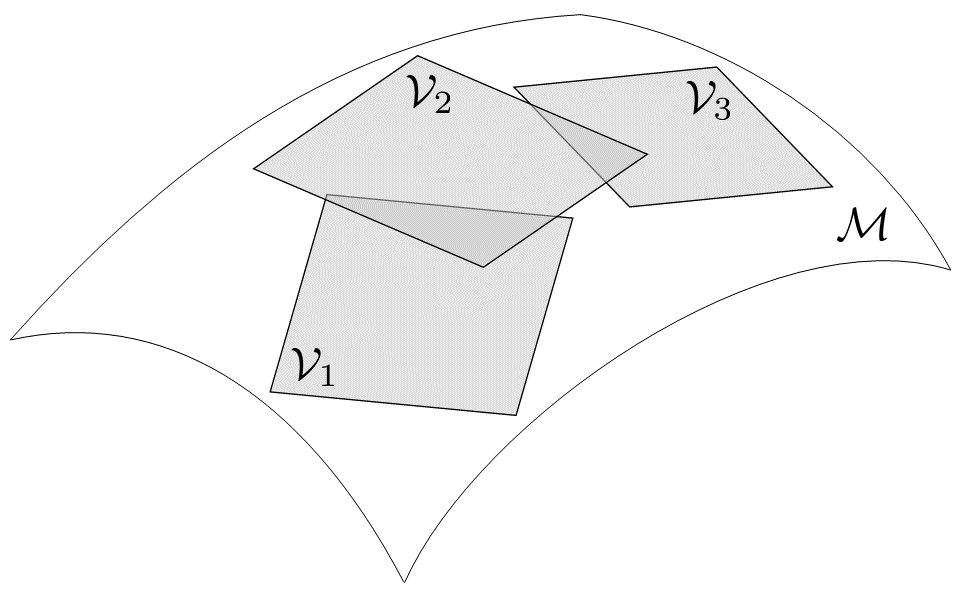}
	%\vglue 0.1 truein
	\caption{Approximation of the HDM-based solution manifold, $\mathcal{M}$, in the time and parameter domain of interest by piecewise-affine subspaces with $N_c = 3$.}
	\label{fig:manifold}
\end{figure}

For the sake of completeness and in order to keep this paper as self-contained as possible, the construction of the $N_c$ local subspaces and the online selection for a queried (time, parameter) 
instance $(t,\bm{\mu})$ of the local approximation subspace are discussed in Sections \ref{sec:pg1.1} and Section \ref{sec:pg1.2}, respectively. The concept of approximating the HDM-based solution 
manifold using a collection of local affine subspaces is illustrated in Figure \ref{fig:manifold}.

The HDM-based problem (\ref{eqn:hdm0}) can be rewritten in residual form as follows
\begin{equation}
	\label{eqn:hdm}
	\bm{r}\left(\bm{u}(t; \bm{\mu}),\dot{\bm{u}}(t; \bm{\mu}),t; \bm{\mu}\right) = \bm{M}(\bm{\mu})\dot{\bm{u}}(t; \bm{\mu}) + \bm{f}(\bm{u}(t; \bm{\mu}); \bm{\mu}) - \bm{g}(t; \bm{\mu}) = 0 \,.
\end{equation}
Substituting (\ref{eqn:subspace}) into (\ref{eqn:hdm}) and projecting the resulting overdetermined nonlinear system of equations onto the left subspace defined by the left ROB 
$\bm{W}_k \in \mathbb{R}^{N\times n_k}$, $k = 1, \ldots, N_c$, yields an {\it implicit} database of $N_c$ local nonlinear PG PROMs associated with the HDM underlying (\ref{eqn:hdm}). Here, the word
{\it implicit} is used to specify that the $N_c$ local nonlinear PG PROMs are not pre-computed and stored in the database; they are constructed online using the underlying local ROBs
which are pre-computed and stored in the database. Each local, nonlinear PG PROM-based problem can be written as
\begin{equation}
	\label{eqn:prom}
	\bm{r}_{k}\left(\bm{y}_k(t; \bm{\mu}),\dot{\bm{y}}_k(t; \bm{\mu}),t; \bm{\mu}\right) = \bm{W}_k^T\bm{r}\left(\bm{u}_{o,k}+\bm{V}_k\bm{y}_k(t; \bm{\mu}),\bm{V}_k\dot{\bm{y}}_k(t; \bm{\mu}),t; \bm{\mu}\right) = 0\vspace{\belowdisplayskip}
\end{equation}
where $\bm{r}_{k}\left(\bm{y}_k(t; \bm{\mu}),\dot{\bm{y}}_k(t; \bm{\mu}),t; \bm{\mu}\right) \in \mathbb{R}^{n_k}$.

For many applications, it has been shown that using a PG projection instead of a Galerkin one, where $\bm{W}_k=\bm{V}_k$ $\forall k$, to construct a PROM of the form given in (\ref{eqn:prom}) endows the 
PROM with superior numerical stability and accuracy properties. For example for steady-state linear problems, PG projection has been advocated to guarantee the asymptotic stability of the resulting PROM 
by bounding the approximation error \cite{buithanh2008}. For time-dependent linear problems, PG projection has been proposed to guarantee that the resulting linear PROM satisfies the Lyapunov stability
criterion \cite{amsallem2012}. For nonlinear problems, it was shown that a for a specific class of PG projections associated with a specific method for constructing a left ROB $\bm{W}_k$, the resulting
PROMs produce solutions that minimize the time-discrete nonlinear HDM-based residual over the affine approximation subspace, even when the HDM is characterized by non-SPD tangent matrices
\cite{grimberg2020} -- a condition which typically destroys the optimal properties of a Galerkin projection.

\subsection{Construction of a piecewise-affine local subspace of approximation}
\label{sec:pg1.1}

Several methods have been proposed for constructing the local approximation subspace (\ref{eqn:subspace}). All of them rely on partitioning the collected set of $N_s$ training snapshots 
$\mathcal{S} = \{\bm{u}^{(s)}\}_{s=1}^{N_s}$, where $\bm{u}^{(s)} = \bm{u}^m(\bm{\mu}^q)$ is a discrete approximation of $\bm{u}(t^m; {\mu}^q)$, $t^m \in [0, T_f]$, and 
$\bm{\mu}^q \in \mathcal{P}$. However, these methods differ by how they specifically partition $\mathcal S$ into subsets of solution snapshots. For example, snapshot partitioning has been performed by 
simply partitioning the time \cite{drohmann2011} or parameter \cite{haasdonk2011, he2020} domain. Alternatively, state space (or HDM-based solution manifold) partitioning has been advocated and 
realized by clustering and compressing the solution snapshots \cite{amsallem2012b}: this enables the construction of local, nonlinear PROMs capable of capturing the different regimes and features 
(e.g., discontinuities and fronts) that may be experienced by the solution of an HDM such as (\ref{eqn:hdm0}) \cite{amsallem2012b}, as well as capturing the effects on this solution of variations in the 
parameters of such an HDM \cite{washabaugh2012}. All of these approaches accommodate the hyperreduction process.

In this work, all local subspace approximations of the form (\ref{eqn:subspace}) are performed by partitioning the state space (or solution manifold) as first proposed in \cite{amsallem2012b}.
For this purpose, the collected solution snapshots are first divided into $N_c$ non-overlapping clusters $\mathcal{S}_k$, $k = 1, \ldots, N_c$, where 
$\bigcup\limits_{k = 1}^{N_c} \mathcal{S}_k = \mathcal{S}$ and $\mathcal{S}_k$ $\bigcap$ $\mathcal{S}_l = \emptyset$ $\forall k \neq l$, using the $k$-means algorithm\cite{hastie2009} equipped with the
standard Euclidean distance. Then, the postprocessing procedure summarized in Algorithm \ref{alg:overlap} \cite{washabaugh2012} is applied to introduce some amount of overlap between the generated 
clusters to promote the continuity of the approximation in the transitional regions. The modified snapshot clusters are denoted by $\{{\mathcal S}_k^+\}_{k=1}^{N_c}$. Next, for each affine 
subspace approximation, the affine offset $\bm{u}_{o,k}$ is specified -- for example, as the problem initial condition $\bm{u}^0(\bm{\mu})$ or the cluster centroid $\bm{u}_{c,k}$ defined here as follows
\begin{equation*}
	\bm{u}_{c,k} = \frac{1}{N_{s, k}} \sum_{s=1}^{N_{s, k}} \bm{u}^{(s, k)},  \qquad k = 1, \ldots, N_c
\end{equation*}
where $\bm{u}^{(s, k)}$ denotes a generic training snapshot $\bm{u}^{(s)}$ assigned to the $k$-th cluster and $N_{s, k}$ denotes the total number of training snapshots assigned to this cluster.
Finally, each right ROB $\bm{V}_k$ appearing in (\ref{eqn:subspace}) is constructed by shifting the snapshots collected in $\mathcal{S}_k$ by $\bm{u}_{o,k}$, then compressing them using the singular
value decomposition (SVD) method.

\begin{remark}
	In Algorithm \ref{alg:overlap}, $\phi \in [0, 1]$ is a user-specified tolerance for controlling the amount of overlap between the generated clusters. Alternatively, a different clustering 
	algorithm such as fuzzy $c$-means \cite{bezdek1984} may be used to generate in one shot a set of overlapping clusters, albeit at an increased computational cost due to an increase in the number of 
	distances to be computed at each iteration of the algorithm.
\end{remark}

\begin{algorithm}
	\caption{Introduction of overlapping into a set of initially disjoint snapshot clusters (based on \cite{washabaugh2012})}
	\textbf{Input:} $\mathcal{S}_k$, $k=1,\ldots,N_c$, $\phi\in[0,1]$\\
	\textbf{Output:} $\mathcal{S}^+_k$, $k=1,\ldots,N_c$
	\begin{algorithmic}[1]
		  \State $\mathcal{S}^+_k \leftarrow \mathcal{S}_k$, $nei_k \leftarrow \emptyset$, for $k=1,\ldots,N_c$
		  \For {$\bm{u}^{(s)} \in \mathcal{S}$} \Comment build inter-cluster connectivity
		    \State $(k,l)\leftarrow$ closest two cluster centers to $\bm{u}^{(s)}$
		    \State $nei_k \leftarrow nei_k \cup l$
		    \State $nei_l \leftarrow nei_l \cup k$
		  \EndFor
		  \For {$k\in\{1,\ldots,N_c\}$} \Comment augment clusters to add overlap
		    \For {$l\in nei_k$}
		      \State $\mathcal{I} \leftarrow$ closest $\phi \times \abs{\mathcal{S}_l}$ snapshots to cluster $k$ in $\mathcal{S}_l$
		      \State $\mathcal{S}^+_k \leftarrow \mathcal{S}^+_k \cup \mathcal{I}$
		    \EndFor
		  \EndFor
	\end{algorithmic}
	\label{alg:overlap}
\end{algorithm}

\subsection{Online selection of the most-appropriate local subspace of approximation}
\label{sec:pg1.2}

At any given time $t$, the most-appropriate local affine subspace in which to approximate the HDM-based solution -- or equivalently, which ROB $\bm{V}_k$ and affine offset $\bm{u}_{o,k}$,
$k = 1, \ldots, N_c$, to select for constructing and solving the most-appropriate PROM-based problem $\bm{r}_k = 0$ (see (\ref{eqn:prom})) in order to approximate the HDM-based solution --
can be determined online as described below, based on the distances between the PROM-based solution at a sufficiently close previous time $t^-$ and the cluster centroids $\bm{u}_{c,k}$.

At time $t = 0$, the closest cluster to the initial condition and the corresponding most-appropriate initial PROM-based solution are given by 
\begin{equation*}
	k^0 = \argmin_{l \in \{1,\ldots,N_c\}} d\left(\bm{u}^0(\bm{\mu}), \bm{u}_{c,l}\right)
\end{equation*}
where $d(\bm{v},\bm{w}) = \norm{\bm{v}-\bm{w}}_2$, and 
\begin{equation*}
	\bm{y}_{k^0}(0; \bm{\mu}) = \bm{V}_{k^0}^T\left(\bm{u}^0(\bm{\mu})-\bm{u}_{o,k^0}\right) \,.
\end{equation*}
Hence, if at any time $t>0$ the PROM-based solution $\bm{y}_{k^-}(t; \bm{\mu})$ is governed by the subspace approximation (\ref{eqn:subspace})
indexed by $k^-$, the closest cluster to the corresponding approximate HDM-based solution is defined as
\begin{equation}
	\label{eqn:distswitch}
	k^+ = \argmin_{l \in \{1, \ldots, N_c\}} d\left(\bm{u}_{o, k^-}+\bm{V}_{k^-}\bm{y}_{k^-}(t; \bm{\mu}), \bm{u}_{c, l}\right) \,.
\end{equation}
When $k^+ \neq k^-$, the approximate solution $\bm{u}_{o,k^-}+\bm{V}_{k^-}\bm{y}_{k^-}(t; \bm{\mu})$ is projected onto the $k^+$-th local affine subspace in order to maintain consistency in the 
approximation, which leads to
\begin{equation}
	\label{eqn:projectswitch}
	\bm{y}_{k^+}(t; \bm{\mu}) = \bm{V}_{k^+}^T\left(\bm{u}_{o,k^-}+\bm{V}_{k^-}\bm{y}_{k^-}(t; \bm{\mu})-\bm{u}_{o,k^+}\right) \,.
\end{equation}

Most importantly, it is noted here that for any pair of indices $(k,l) \in \{1, \ldots, N_c\}\times\{1, \ldots, N\}$, the square of the distance between the approximate solution $\bm{u}_{o,k}+\bm{V}_{k}\bm{y}_{k}(t; \bm{\mu})$ and the centroid of 
the $l$-th cluster of solution snapshots $\bm{u}_{c, l}$ can be written as
\begin{equation*}
	\begin{split}
		d\left(\bm{u}_{o,k}+\bm{V}_{k}\bm{y}_{k}(t; \bm{\mu}), \bm{u}_{c, l}\right)^2 &= \left(\bm{u}_{o, k}+\bm{V}_{k}\bm{y}_{k}(t; \bm{\mu})-\bm{u}_{c, l}\right)^T \left(\bm{u}_{o,k}+\bm{V}_{k}\bm{y}_{k}(t; \bm{\mu})-\bm{u}_{c,l}\right) \\
		&= \norm{\bm{y}_k(t; \bm{\mu})}_2^2+\underbrace{2(\bm{u}_{o,k}-\bm{u}_{c,l})^T\bm{V}_k}_{\textnormal{pre-computable}\,\in\,\mathbb{R}^{1\times n_k}}\bm{y}_k(t; \bm{\mu})+\underbrace{\norm{\bm{u}_{o,k}-\bm{u}_{c,l}}_2^2}_{\textnormal{pre-computable}\,\in\,\mathbb{R}} \,.
	\end{split}
\end{equation*}
Hence, by pre-computing once for all (offline) the $O(N_c^2)$ fixed, cluster-related reduced-order quantities $(\bm{u}_{o,k}-\bm{u}_{c,l})^T\bm{V}_k$ and distances $\norm{\bm{u}_{o,k}-\bm{u}_{c, l}}_2^2$, 
the closest cluster of solution snapshots to the approximate solution $\bm{u}_{o,k}+\bm{V}_{k}\bm{y}_{k}(t; \bm{\mu})$, in the sense defined in (\ref{eqn:distswitch}), can be identified online,
in $O(N_c  n_k)$ operations. Furthermore, for any pair of indices $(k, l)$, the projection (\ref{eqn:projectswitch}) can be expressed as
\begin{equation*}
	\bm{V}_{l}^T\left(\bm{u}_{o,k}+\bm{V}_{k}\bm{y}_{k}(t; \bm{\mu})-\bm{u}_{o,l}\right) = \underbrace{\bm{V}_{l}^T\bm{V}_{k}}_{\textnormal{pre-computable}\,\in\,\mathbb{R}^{n_l\times n_k}}\bm{y}_{k}(t; \bm{\mu})+\underbrace{\bm{V}_{l}^T(\bm{u}_{o,k}-\bm{u}_{o,l})}_{\textnormal{pre-computable}\,\in\,\mathbb{R}^{n_l}}
\end{equation*}
which shows that by pre-computing the additional $O(N_c^2)$ reduced-order matrices $\bm{V}_{l}^T\bm{V}_{k}$ and reduced-order vectors 
$\bm{V}_{l}^T(\bm{u}_{o,k}-\bm{u}_{o,l})$, the projection (\ref{eqn:projectswitch}) can also be computed online, in $O(n_k n_l)$ operations.

In summary, all computations associated with the online selection of a most-appropriate local subspace of approximation can be performed in real-time during the online solution of a PROM-based problem 
of the form given in (\ref{eqn:prom}).

\subsection{Computational bottlenecks}
\label{sec:pg2}

Nevertheless, although the dimension of the PROM-based problem (\ref{eqn:prom}) is $n_k \ll N$, the cost of solving this problem scales with both dimensions $n_k$ and $N \gg n_k$, which typically
prevents the PROM from achieving real-time or near real-time performance. To see this, consider first the case where (\ref{eqn:prom}) is discretized using an {\it explicit} time-integration scheme.
The evaluation of the PROM-based solution at the $m$-th computational time-step for a queried parameter point $\bm{\mu}^{\star} \in \mathcal{P}$ requires in this case the assembly of the following
reduced-order quantities:
\begin{itemize}
	\item $\bm{M}_{k}(\bm{\mu}^{\star}) = \bm{W}_k^T\bm{M}(\bm{\mu}^{\star})\bm{V}_k \in \mathbb{R}^{n_k\times n_k}$, which is time-independent and therefore must be performed only once for a queried
		parameter point $\bm{\mu}^{\star}$, and requires $O(N^2 n_k)$ operations \big(though exploiting sparsity can reduce this computational complexity to $O(N n_k)$\big). 
	\item $\bm{b}_{k}^{m}\left(\bm{y}_k^m(\bm{\mu}^{\star}), t^m; \bm{\mu}^{\star}\right) = \bm{W}_k^T\left(\bm{f}(\bm{u}_{o,k}+\bm{V}_k\bm{y}_k^m(\bm{\mu}^{\star}); \bm{\mu}^{\star})
		-\bm{g}^m(\bm{\mu}^{\star})\right)\in \mathbb{R}^{n_k}$, which is time-dependent, must be performed at each computational time-step, and requires $O(N n_k)$ operations. Note that here and
		throughout the remainder of this paper, the superscript $m$ designates as earlier a fully-discrete quantity evaluated at time $t^m$.

\end{itemize}
For an {\it implicit} time-integration scheme, the solution by Newton's method -- or any similar method -- of the $n_k$-dimensional nonlinear system of algebraic equations arising at each computational 
time-step requires, in addition to the assembly of the reduced-order quantities mentioned above, the construction at each Newton iteration of the following reduced-order quantities:
\begin{itemize}
	\item $\bm{J}_{k}^m\left(\bm{y}_k^m(\bm{\mu}^{\star}); \bm{\mu}^{\star}\right) = \bm{W}_k^T{\bm{J}^m}\left(\bm{u}_{o,k}+\bm{V}_k\bm{y}_k^m(\bm{\mu}^{\star}); 
		\bm{\mu}^{\star}\right)\bm{V}_k \in \mathbb{R}^{n_k\times n_k}$, where $\bm{J}^m(\bm{u}^m(\bm{\mu}); \bm{\mu})\in\mathbb{R}^{N\times N}$ is the Jacobian matrix of the discrete
		counterpart of the HDM-based residual (\ref{eqn:hdm}) with respect to $\bm u$ -- that is,
	\begin{equation}
		\label{eqn:hdmjac}
		\bm{J}^m\left(\bm{u}^m(\bm{\mu}); \bm{\mu}\right) = \bm{M}(\bm{\mu})\frac{\partial \hat{\dot{\bm{u}}}}{\partial \bm{u}}\left(\bm{u}^m(\bm{\mu})\right)+\frac{\partial \bm{f}}{\partial \bm{u}}\left(\bm{u}^m(\bm{\mu}); \bm{\mu}\right)
	\end{equation}
	where $\hat{\dot{\bm{u}}}\left(\bm{u}^m(\bm{\mu})\right)$ is the fully-discrete approximation of the time-derivative $\dot{\bm{u}}^m(\bm{\mu})$ at time-instance $t^m$. For most time-integration 
	schemes -- including linear multistep and diagonally implicit Runge-Kutta (DIRK) schemes -- $\hat{\dot{\bm{u}}}(\bm{u}^m(\bm{\mu}); \bm{\mu})$ can be written as
	\begin{equation*}
		\hat{\dot{\bm{u}}}\left(\bm{u}^m(\bm{\mu})\right) = \bm{p}\left(\bm{u}^m(\bm{\mu}), \bm{u}^{m-1}(\bm{\mu}), \ldots, \bm{u}^{m-s}(\bm{\mu}), \Delta t^{m}, \Delta t^{m-1}, \ldots, 
		\Delta t^{m-s} \right)
	\end{equation*}
	where $\bm{p}$ is the function characterizing the time-discretization scheme, $s$ is its stencil width, and $\Delta t$ is the time-step size. Hence, the computation of the above reduced-order 
	Jacobian matrix requires in general $O(N^2 n_k)$ operations -- or $O(N n_k)$ operations when sparsity patterns can be perfectly exploited.
\end{itemize}

From the above discussion, it follows that even for a single parameter $\bm{\mu}^{\star}$, the processing of the the time-dependent PROM-based problem (\ref{eqn:prom}) suffers a computational bottleneck 
due to a computational complexity that scales not only with the smaller dimension $n_k$ of this PROM, but also with the much larger dimension $N$ of the HDM.

As explained in Section \ref{sec:intro}, it is sometimes possible to mitigate the aforementioned computational bottleneck using a divide and conquer strategy. For many applications however
-- for example, for compressible flow problems formulated in terms of conservation state variables  -- this is simply not possible, particularly when the left ROB $\bm{W}_k$, which may
depend on the vector of generalized coordinates \cite{carlberg2011, carlberg2013}, introduces complex nonlinearities in the PROM. For such applications, achieving computational efficiency in the 
solution of the PROM-based problem (\ref{eqn:prom}) requires in general the hyperreduction of the underlying PROM and more specifically, the hyperreduction at each computational time-step of the 
projected or reduced-order quantities $\bm{W}_k^T\left(\bm{f}(\bm{u}_{o,k}+\bm{V}_k\bm{y}_k^m(\bm{\mu}^{\star}); \bm{\mu}^{\star}) -\bm{g}^m(\bm{\mu}^{\star})\right)\in \mathbb{R}^{n_k}$ and
$\bm{W}_k^T\bm{J}^m\left(\bm{u}_{o,k}+\bm{V}_k\bm{y}_k^m(\bm{\mu}^{\star}); \bm{\mu}^{\star}\right)\bm{V}_k \in \mathbb{R}^{n_k\times n_k}$ discussed above.

\section{Hyperreduction via mesh sampling and weighting}
\label{sec:ecsw}

\subsection{Hyperreduction of a reduced-order nonlinear residual vector}
\label{sec:ecsw1}

In this section, an ECSW-type hyperreduction method is developed for nonlinear, PG PROMs to enable their processing in a computational complexity that scales only with the small dimension 
$n_k \ll N$ of the local PROM,  $k \in \{1, \ldots, N_c\}$. Specifically, the hyperreduction method described herein is independent of the semi-discretization scheme 
underlying the $\bm{\mu}$-parametric, HDM (\ref{eqn:hdm0}) or its residual form (\ref{eqn:hdm}).

Let $\mathcal{E} = \left\{e_1,\ldots,e_{N_e}\right\}$, with $N_e=\abs{\mathcal{E}}$, denote the set of mesh entities describing the discretization of the computational domain of interest.
These entities may be: finite elements in the case of a FE semi-discretization of the PDE underlying (\ref{eqn:hdm0}) and (\ref{eqn:hdm}); primal cells in the case of a cell-centered finite volume (FV) 
semi-discretization; dual cells in the case of a vertex-based FV semi-discretization; or simply vertices corresponding to collocation points in the case of a finite difference (FD) semi-discretization
of the aforementioned PDE. In all cases, the computation of the reduced-order residual (\ref{eqn:prom}) can be written as
\begin{equation}
	\label{eqn:res1}
	\begin{split}
		\bm{r}_{k}\left(\bm{y}_k(t; \bm{\mu}),\dot{\bm{y}}_k(t; \bm{\mu}),t; \bm{\mu}\right) &= \bm{W}_k^T\bm{r}\left(\bm{u}_{o,k}+\bm{V}_k\bm{y}_k(t; \bm{\mu}),\bm{V}_k\dot{\bm{y}}_k(t; \bm{\mu}),t; 
		\bm{\mu}\right) \\
		&= \sum_{e \, \in \, \mathcal{E}} \bm{W}_k^T\bm{L}_e^T\bm{r}_e\left(\bm{L}_{e^+}\left(\bm{u}_{o,k}+\bm{V}_k\bm{y}_k(t; \bm{\mu})\right), \bm{L}_{e^+}\bm{V}_k\dot{\bm{y}}_k(t; \bm{\mu}),t; \bm{\mu}\right)
	\end{split}
\end{equation}
where:
\begin{itemize}
	\item $\bm{L}_e\in\{0,1\}^{d_e\times N}$ is the Boolean matrix that localizes a high-dimensional global vector (dimension $N$) to the $d_e$ DOFs attached to the mesh entity $e$.
	\item $\bm{L}_{e^+}\in\{0,1\}^{d_{e^+}\times N}$ is the Boolean matrix that localizes a high-dimensional global vector (dimension $N$) to the $d_{e^+} \ge d_e$ DOFs attached to the mesh entity $e$
		and a set of neighboring mesh entities determined by the stencil of the chosen semi-discretization method: for example, $d_{e^+} = d_e$ for a FE semi-discretization,
		$d_{e^+} > d_e$ for a cell-centered or vertex-based FV semi-discretization, and $d_{e^+} > d_e$ for an FD semi-discretization.
	\item $\bm{r}_e\left(\bm{L}_{e^+}\left(\bm{u}_{o,k}+\bm{V}_k\bm{y}_k(t; \bm{\mu})\right),\bm{L}_{e^+}\bm{V}_k\dot{\bm{y}}_k(t; \bm{\mu}),t; \bm{\mu}\right) \in \mathbb{R}^{d_e}$ is the 
		contribution of the mesh entity $e$ to the global, HDM-based residual $\bm{r} \in {\mathbb R}^N$.
\end{itemize}
The fully-discrete counterpart of the reduced-order, semi-discrete residual (\ref{eqn:res1}) is
\begin{equation}
	\label{eqn:res1fd}
	\bm{r}_{k}^m\left(\bm{y}_k^m(\bm{\mu}), t^m; \bm{\mu}\right) = \sum_{e \, \in \, \mathcal{E}} \bm{W}_k^T\bm{L}_e^T\bm{r}_e^m\left(\bm{L}_{e^+}\left(\bm{u}_{o,k}+\bm{V}_k\bm{y}_k^m(\bm{\mu})
	\right), t^m; \bm{\mu}\right) \,.
\end{equation}

Now, let $\widetilde{\mathcal{E}} \subset \mathcal{E}$ denote an optimally sampled subset of the set of mesh entities $\mathcal{E}$, with $\widetilde{N}_e = \abs{\widetilde{\mathcal{E}}} \ll N_e$
(see Section \ref{sec:sampling}). A natural approach for approximating the PROM-based residual (\ref{eqn:res1fd}) can be written as
\begin{equation}
	\label{eqn:res2}
	\bm{r}_{k}^m\left(\bm{y}_k^m(\bm{\mu}), t^m; \bm{\mu}\right) \approx \tilde{\bm{r}}_{k}^m\left(\bm{y}_k^m(\bm{\mu}), t^m; \bm{\mu}\right) 
	= \sum_{e \, \in \, \widetilde{\mathcal{E}}} \xi_e \bm{W}_k^T\bm{L}_e^T\bm{r}_e^m\left(\bm{L}_{e^+}\left(\bm{u}_{o,k}+\bm{V}_k\bm{y}_k^m(\bm{\mu})\right), t^m; \bm{\mu}\right)
\end{equation}
where $\tilde{\bm{r}}_{k}^m$ is referred to as the hyperreduced fully-discrete residual vector and the set of weights $\left\{\xi_e \mid e \in \widetilde{\mathcal{E}}\right\}$ associated with the 
set of sampled mesh entities $\widetilde{\mathcal{E}}$ leads to interpreting the approximation scheme (\ref{eqn:res2}) as a generalized quadrature rule. Since $\widetilde{N}_e \ll N_e$, it follows 
from (\ref{eqn:res2}) that the hyperreduced fully-discrete vector $\tilde{\bm{r}}_{k}^m$ can be efficiently evaluated at any time-instance and/or queried parameter point in a number of operations that 
is independent of the dimension $N$ of the HDM. 

\begin{remark}
	It is noted that unlike hyperreduction methods of the approximate-then-project type, the hyperreduction approximation (\ref{eqn:res2}) does not attempt to accurately represent the high-dimensional
	residual $\bm{r}^m$, but only of its projection onto the left subspace spanned by the columns of the left ROB $\bm{W}_k$.
\end{remark}

\subsection{Hyperreduction of a reduced-order Jacobian matrix}
\label{sec:ecsw2}

In the event where Newton's method or a variant is chosen for solving a steady-state counterpart of the nonlinear, PROM-based equation (\ref{eqn:prom}), or for solving the nonlinear system of equations 
arising at each time-step of the implicit time-discretization of this PROM-based equation, the resulting reduced-order Jacobian matrix may need to be reconstructed at each Newton iteration and at least 
at each time-step in the latter case. In such circumstances, computational efficiency calls for hyperreducing this PROM matrix when it arises.

Here, attention is focused on the time-dependent case and an implicit time-discretization, as the explicit scenario is a sub-case of this case and the steady-state scenario is a particular instance of 
this case. At each time-step $t^m$, the reduced-order counterpart of the Jacobian matrix (\ref{eqn:hdmjac}) can be written as
\begin{equation}
	\label{eqn:jac1}
	\bm{J}_{k}^m\left(\bm{y}_k^m(\bm{\mu}); \bm{\mu}\right) = \sum_{e \, \in \, \mathcal{E}} \bm{W}_k^T\bm{L}_e^T\bm{J}_e^m\left(\bm{L}_{e^+}(\bm{u}_{o,k}+\bm{V}_k\bm{y}_k^m(\bm{\mu})); \bm{\mu}\right)\bm{L}_{e^+}\bm{V}_k \,.
\end{equation}
As in Section \ref{sec:ecsw1}, $\bm{L}_e$ and $\bm{L}_{e^+}$ denote here the $d_e \times N$ and $d_e^+ \times N$ Boolean matrices that localize a high-dimensional vector (dimension $N$) to the $d_e$ and
$d_e^+$ DOFs associated with the same mesh entity $e$, respectively, and $\bm{J}_e^m\left(\bm{L}_{e^+}(\bm{u}_{o,k}+\bm{V}_k\bm{y}_k^m(\bm{\mu})); \bm{\mu}\right) \in \mathbb{R}^{d_e \times d_e^+}$
is the Jacobian matrix of ${\bm r}_e$ with respect to $\bm{y}_k$. 

Since the reduced-order matrix $\bm{J}_{k}^m(\bm{y}_k^m(\bm{\mu}); \bm{\mu})$ (\ref{eqn:jac1}) is the Jacobian matrix of the fully-discrete PROM-based residual ${\bm r}_k^m$ (\ref{eqn:res1fd}) with 
respect to ${\bm y}_k$, it is proposed here to perform its hyperreduction by computing the Jacobian of the hyperreduced semi-discrete residual $\tilde{\bm{r}}_{k}^m$ (\ref{eqn:res2}) with respect to 
${\bm y}_k$. This leads to the {\it consistent} hyperreduced Jacobian matrix
\begin{equation}
	\label{eqn:jac2}
	\widetilde{\bm{J}}_{k}^m\left(\bm{y}_k^m(\bm{\mu}); \bm{\mu}\right) = \displaystyle{\frac{\partial \tilde{\bm{r}}_{k}^m}{\partial \bm{y}_k} \left(\bm{y}_k^m(\bm{\mu}), t^m; \bm{\mu}\right)} = 
	\sum_{e \, \in \, \widetilde{\mathcal{E}}} \xi_e \bm{W}_k^T\bm{L}_e^T\bm{J}_e^m(\bm{L}_{e^+}(\bm{u}_{o,k}+\bm{V}_k\bm{y}_k^m(\bm{\mu})); \bm{\mu})\bm{L}_{e^+}\bm{V}_k
\end{equation}
and therefore guarantees in principle a good performance of the application of Newton's method or a variant to the solution of the fully-discrete, PROM-based nonlinear problem at hand.

In many applications and PMOR methods, including PG-based methods such as the residual-minimizing least-squares PG (LSPG) projection method \cite{carlberg2011, carlberg2013}, each of the 
reduced-order Jacobian matrix ${\bm{J}}^m_k\left(\bm{y}^m_k(\bm{\mu}); \bm{\mu}\right)$ and mesh-entity-level reduced-order Jacobian matrix 
${\bm W}_k^T\bm{L}_e^T{\bm J}_e\left({\bm L}_{e^+}\left({\bm u}_{o,k}+\bm{V}_k\bm{y}_k^m(\bm{\mu})\right); \bm{\mu}\right){\bm L}_{e^+}\bm{V}_k$ appearing in (\ref{eqn:jac1}) is SPD. In order to 
guarantee that the above hyperreduction approximation preserves this property, both hyperreduction approximations (\ref{eqn:res2}) and (\ref{eqn:jac2}) are equipped here with the constraint $\xi_e > 0$, 
$\forall e\in\widetilde{\mathcal{E}}$. When this constraint is not needed, it can be removed, or simply enforced to enable the use of a unified computational approach for sampling a set 
of mesh entities $\widetilde{\mathcal{E}}$ and computing the associated set of weights $\Xi_{\widetilde{\mathcal{E}}} = \left\{\xi_e \mid e \in \widetilde{\mathcal{E}}\right\}$ (see 
Section \ref{sec:sampling}).

\subsection{Grouped versus individual hyperreduction approximations}

From (\ref{eqn:hdm}), it follows that the fully-discrete residual (\ref{eqn:res1fd}) can also be written as
\begin{equation} 
	\label{eqn:new1}
	\bm{r}_{k}^m\left(\bm{y}_k^m(\bm{\mu}), t^m; \bm{\mu}\right) = \underbrace{\left(\bm{W}_k^T \bm{M}(\bm{\mu})\bm{V}_k\right)}_{\bm{M}_k}\dot{\bm{y}}^m(\bm{\mu}) + \underbrace{\bm{W}_k^T\bm{f}^m\left(\bm{u}_{o,k}+\bm{V}_k\bm{y}_k^m(\bm{\mu}); \bm{\mu}\right)}_{\bm{f}_k^m} - \underbrace{\bm{W}_k^T\bm{g}^m(t^m; \bm{\mu})}_{\bm{g}_k^m} = 0
\end{equation}
where
\begin{align} 
	\bm{M}_{k}(\bm{\mu}) &= \sum_{e \, \in \, \mathcal{E}} \bm{W}_k^T\bm{L}_e^T\bm{M}_e(\bm{\mu})\bm{L}_e\bm{V}_k \label{eqn:Hyp2a} \\
	\bm{f}_{k}^m\left(\bm{u}_{o,k}+\bm{V}_k\bm{y}_k^m(\bm{\mu}); \bm{\mu}\right) &= \sum_{e \, \in \, \mathcal{E}} \bm{W}_k^T\bm{L}_e^T\bm{f}_e^m\left(\bm{L}_{e^+}\left(\bm{u}_{o,k}+\bm{V}_k\bm{y}_k^m(\bm{\mu}) \right); \bm{\mu}\right)\label{eqn:Hyp2b} \\
	\bm{g}_{k}^m(t^m; \bm{\mu}) &= \sum_{e \, \in \, \mathcal{E}} \bm{W}_k^T\bm{L}_e^T\bm{g}_e^m(t^m; \bm{\mu}) \label{eqn:Hyp2c} \,.
\end{align}
At this point, the reader may ask why in Section \ref{sec:ecsw1} hyperreduction was applied in (\ref{eqn:res2}) to the above reduced-order terms as a group, specifically, collectively as the 
fully-discrete residual $\bm{r}_{k}^m\left(\bm{y}_k^m(\bm{\mu}), t^m; \bm{\mu}\right)$, rather than to each of them individually. Before answering this question, a noteworthy observation that 
differentiates between the hyperreduction of a Galerkin PROM and that of a PG PROM is discussed below. 

Consider first the case of a Galerkin PROM ($\bm{W}_k = \bm{V}_k$) and an HDM where only the nonlinear flux vector is parametric -- that is, where $\bm{f}$ depends on the vector of parameters $\bm{\mu}$
and on the solution $\bm{u}$, but the mass matrix $\bm{M}$ is constant and the source term $\bm{g}$ depends only on time $t$. In this case, both $\bm{M}_{k} = \bm{V}_k^T \bm{M}\bm{V}_k$ and 
$\bm{g}_{k}^m(t^m) = \bm{V}_k^T\bm{g}^m$ are pre-computable and therefore only the reduced-order vector $\bm{f}_{k}^m\left(\bm{u}_{o,k}+\bm{V}_k\bm{y}_k^m(\bm{\mu}); \bm{\mu}\right)$ needs be considered 
for hyperreduction. For many highly nonlinear applications (e.g., high-speed compressible fluid flow and solid mechanics with finite-strain viscoelasticity), the reduced-order term 
$\bm{f}_{k}^m\left(\bm{u}_{o,k}+\bm{V}_k\bm{y}_k^m(\bm{\mu}); \bm{\mu}\right)$ can be expected to exhibit non-polynomial nonlinearities at least with respect to $\bm{u}$ 
\big(specifically, $\left(\bm{u}_{o,k}+\bm{V}_k\bm{y}_k^m(\bm{\mu})\right)$\big) and hence to benefit from hyperreduction. Hence, in this scenario where only the reduced-order term 
$\bm{f}_{k}^m\left(\bm{u}_{o,k}+\bm{V}_k\bm{y}_k^m(\bm{\mu}); \bm{\mu}\right)$ (\ref{eqn:Hyp2b}) needs to be hyperreduced at each time-step $t^m$, the set of mesh entities 
$\widetilde{\mathcal{E}} \subset \mathcal{E}$ and the associated set of weights $\Xi_{\widetilde{\mathcal{E}}} = \left\{\xi_e \mid e \in \widetilde{\mathcal{E}}\right\}$ should be 
determined so that the application to (\ref{eqn:Hyp2b}) of the approximation defined in (\ref{eqn:res2}) is sufficiently accurate. 

Consider next the case of a PG PROM ($\bm{W}_k \ne \bm{V}_k$) and a parametric or non-parametric HDM. In this case, the ability to pre-compute any of the terms (\ref{eqn:Hyp2a}), (\ref{eqn:Hyp2b}), and 
(\ref{eqn:Hyp2c}) also depends on whether the left ROB $\bm{W}_k$ is fixed in the parameter and time spaces and is independent of the solution $\bm{u}_k(t; \mu)$ \big(or its subspace approximation 
$\bm{u}_{o,k}+\bm{V}_k\bm{y}_k^m(\bm{\mu})$\big). For many PG-based PMOR methods -- for example, for the residual-minimizing LSPG projection method \cite{carlberg2011, carlberg2013} -- the left ROB 
evolves together with the reduced-order vector of generalized coordinates $\bm{y}_k^m$. For such a method, none of the reduced-order terms (\ref{eqn:Hyp2a}--\ref{eqn:Hyp2c}) can be 
pre-computed and therefore all of them must be hyperreduced in the presence of non-polynomial nonlinearities with respect to any of their variables. Hence, in this case -- which can be expected to be 
typical in the context of PG PROMs -- and in other cases where at least two of the reduced-order terms defining the PG PROM are to be hyperreduced, the question becomes whether these reduced-order terms 
should be hyperreduced as a group, or individually.

Hyperreducing multiple reduced-order terms as a group as in the hyperreduction of the parametric, nonlinear, fully-discrete residual (\ref{eqn:res2}) is most cost effective. It requires training
a single set of mesh entities $\widetilde{\mathcal{E}} \subset \mathcal{E}$ and associated weights $\Xi_{\widetilde{\mathcal{E}}} = \left\{\xi_e \mid e \in \widetilde{\mathcal{E}}\right\}$ as described in 
the next section. On the other hand, hyperreducing multiple reduced-order terms individually can be performed in two different ways: 
\begin{itemize}
	\item {${\mathcal I}1$}: By sampling multiple sets of mesh entities and computing for each an associated set of weights.
	\item {${\mathcal I}2$}: By computing a single set of mesh entities $\widetilde{\mathcal{E}}$ and its associated set of weights $\Xi_{\widetilde{\mathcal{E}}}$ so that the application of the 
		generalized quadrature rule to the approximation of each reduced-order term individually is sufficiently accurate. 
\end{itemize}
The approach ${\mathcal I}1$ is clearly less computationally efficient than hyperreducing multiple reduced-order terms as a group. The approach ${\mathcal I}2$ is less cumbersome than 
its counterpart ${\mathcal I}1$ as it does not involve manipulating multiple sets of mesh entities and associated sets of weights. Nevertheless, the approach ${\mathcal I}2$ is less computationally
efficient than hyperreducing multiple reduced-order terms as a group as it entails training $\widetilde{\mathcal{E}}$ and $\Xi_{\widetilde{\mathcal{E}}}$ to be accurate for many more approximations. 
However, it can be expected to be more accurate for the same computational cost.

There are a few applications for which achieving a desired level of accuracy requires performing necessary hyperreductions at the individual level. An example is the fast solution of parametric
generalized eigenvalue problems by HPROMs \cite{farhat2019}, where two sets of unknowns -- namely, the eigenvalues and the corresponding eigenvectors -- are governed by a single algebraic 
equation that can be written in residual form. In this case, accuracy dictates hyperreducing the generalized mass and stiffness matrices individually; and for this purpose,
computational efficiency as well as practicality call for adopting the approach labeled above as ${\mathcal I}2$. For the class of problems considered in this paper however, which is represented by the 
HDM (\ref{eqn:hdm0}) and involves for each $k \in  \{1, \ldots, N_c\}$ a single set of unknowns $\bm{y}_k^m$, performing hyperreduction at the group level as in (\ref{eqn:res2}) -- which is most
computationally efficient -- can deliver the desired level of accuracy, as show in Section \ref{sec:app} for several challenging problems. Hence for this class of problems, necessary hyperreduction
is performed here at the group level as in (\ref{eqn:res2}).

\section{Implementation}
\label{sec:sampling}

\subsection{Mesh sampling and weighting}
\label{sec:sampling1}

Given a set of mesh entities $\mathcal{E} = \left\{e_1,\ldots,e_{N_e}\right\}$ describing the discretization of the computational domain of interest, the mesh sampling and weighting problem is defined 
here as that of finding the smallest subset of mesh entities $\widetilde{\mathcal{E}} \subset \mathcal{E}$ and the associated set of weights 
$\Xi_{\widetilde{\mathcal{E}}} = \left\{\xi_e \mid e \in \widetilde{\mathcal{E}}\right\}$ for which the hyperreduction approximations (\ref{eqn:res2}) and (\ref{eqn:jac2}) are sufficiently accurate 
(in the sense specified below). As shown in \cite{farhat2014, farhat2015} for the case of a nonlinear Galerkin PROM based on a single global ROB, this can can be achieved using what is known nowadays 
as a supervised training or machine learning approach. Furthermore, since (\ref{eqn:jac2}) is nothing but the Jacobian of (\ref{eqn:res2}) with respect to $\bm{y}_k$, it suffices to train the pair of 
sets $\left(\widetilde{\mathcal{E}}, \Xi_{\widetilde{\mathcal{E}}}\right)$ only for the hyperreduction approximation (\ref{eqn:res2}).

In the context of local subspace approximations, two different mesh sampling and weighting approaches can be considered:
\begin{itemize}
	\item An approach where a single set of mesh entities $\widetilde{\mathcal{E}}$ and an associated set of weights $\Xi_{\widetilde{\mathcal{E}}}$ is constructed by training both sets for all local 
		approximations $k = 1, \ldots, N_c$, simultaneously.
	\item An alternative approach where $N_c$ pairs of sets $\left(\widetilde{\mathcal{E}}_k, \Xi_{\widetilde{\mathcal{E}}_k}\right)$, $k = 1, \ldots, N_c$, are constructed and each is trained for a 
		single instance of the local approximation defined in (\ref{eqn:subspace}).
\end{itemize}
For a fixed level of accuracy, the first approach can be expected to lead to a set of mesh entities $\widetilde{\mathcal{E}}$ that is larger than any of the counterpart sets $\widetilde {\mathcal E}_k$ 
delivered by the second approach and therefore to deliver a lesser online performance from the wall-clock time viewpoint. Nevertheless, this approach is simpler than the second one which requires an 
elaborate computer implementation and a considerable amount of bookkeeping. For this reason, only the first approach is considered here. It is noted however that the implementation of the second 
approach can be built around that of the first one in a rather straightforward manner.

For both convenience and computational efficiency, the training of the pair of sets $\left(\widetilde{\mathcal{E}}, \Xi_{\widetilde{\mathcal{E}}}\right)$ is not performed directly on generalized coordinates 
snapshots $\bm{y}_{k}^{(s)}$, $k\in\{1,\ldots,N_c\}$, because this would be computationally inefficient for the following reasons:
\begin{itemize}
	\item Such a training would require two different series of offline simulations, namely:
		\begin{itemize} 
			\item A series of HDM-based simulations designed to generate a first set of solution snapshots $\mathcal{S}_k$ (see Section \ref{sec:pg1.1}) for constructing the right and left
				local ROBs $\bm{V}_k$ and $\bm{W}_k$, respectively, $k\in\{1,\ldots,N_c\}$, as well as the associated nonlinear PG PROMs.  
			\item Another series of nonlinear simulations based on these non hyperreduced PG PROMs to generate $N_c$ additional sets of solution snapshots of the form 
				$\mathcal{Y}_k = \left \{\bm{y}_k^m(\bm{\mu}^q)\right \}$, $k\in\{1,\ldots,N_c\}$, for hyperreducing the constructed PG PROMs.  
		\end{itemize}
	\item The second series of PROM-based computations would be very compute intensive due to the absence of hyperreduction. 
\end{itemize}
Instead, the pair of sets $\left(\widetilde{\mathcal{E}}, \Xi_{\widetilde{\mathcal{E}}}\right)$ is efficiently trained here using a subset ${\mathcal S}_H$ of the same set of solution snapshots $\mathcal{S}$ 
used for constructing the local subspace approximations (see Section \ref{sec:pg1}) -- that is, $\mathcal{S}_H \subseteq \mathcal{S}$ -- as follows. For any snapshot $\bm{u}^{(s)}\in\mathcal{S}_H$, let
$k_s \in \{1,\ldots,N_c\}$ denote the index of the unique cluster $\mathcal{S}_k$ containing this snapshot $\left(\mathcal{S}_k \ni \bm{u}^{(s)}\right)$. This solution snapshot can be converted 
{\it on the fly} into its approximation $\bm{u}_{o,k_s}+\bm{V}_{k_s}\bm{y}_{k_s}$ using the orthogonal projector 
\begin{equation*}
	\Pi^{\bot}_{\bm{V}_{k_s}} = \bm{V}_{k_s}\bm{V}_{k_s}^T
\end{equation*}
-- that is, the orthogonal projector onto the subspace spanned by the columns of the local right ROB $\bm{V}_{k_s}$. Hence, each solution snapshot $\bm{u}^{(s)}\in\mathcal{S}_H$ is
transformed here on the fly into the following vector of generalized coordinates $\bm{y}_{k_s}$
\begin{equation*}
	\bm{y}_{k_s} = \bm{V}_{k_s}^T\left(\bm{u}^{(s)}-\bm{u}_{o,k_s}\right) \,.
\end{equation*}
Next, the corresponding discrete residual $\bm{r}^m\left(\bm{u}_{o,k_s}+\Pi^{\bot}_{\bm{V}_{k_s}}(\bm{u}^{(s)}-\bm{u}_{o,k_s}),t^m; \bm{\mu}^q\right)$ is computed and used for training 
$\left(\widetilde{\mathcal{E}}, \Xi_{\widetilde{\mathcal{E}}}\right)$.

Now, let
\begin{equation*}
	\begin{array}{r c l l}
		\bm{c}_{s e} &=& \bm{W}_{k_s}^T\bm{L}_e^T{\bm{r}_e^m\left(\bm{L}_{e^+}\left(\bm{u}_{o,k_s}+\Pi^{\bot}_{\bm{V}_{k_s}}(\bm{u}^{(s)}-\bm{u}_{o,k_s})\right),
		t^m; \bm{\mu}^q\right)} \in \mathbb{R}^{n_{k_s}}, & \qquad s = 1, \ldots, N_H \\
		\bm{d}_{s} &=& \sum\limits_{e \, \in \, \mathcal{E}} \bm{c}_{k_{s_e}} = \bm{r}_{k_s}^m(\bm{y}_{k_s},t^m; \bm{\mu}^q) \in \mathbb{R}^{n_{k_s}}, & \qquad s = 1, \ldots, N_H
	\end{array}
\end{equation*}
where $N_H = \abs{\mathcal{S}_H}$. Using the above notation, the exact assembly of the training data on the original mesh (associated with the HDM) can be written as
\begin{equation*}
	\bm{C}\bm{1}=\bm{d}
\end{equation*}
where 
\begin{equation}
	\label{eqn:Cd}
	\bm{C} = \begin{bmatrix} \bm{c}_{1 1} & \ldots & \bm{c}_{1 N_e} \\ \vdots & \ddots & \vdots \\ \bm{c}_{{N_H} 1} & \ldots & \bm{c}_{{N_H} N_e} \end{bmatrix} 
	\in \mathbb{R}^{\left(\sum\limits_{s=1}^{N_H} n_{k_s}\right) \times N_e} \qquad \qquad \bm{d} = \begin{bmatrix} \bm{d}_1 \\ \vdots \\ \bm{d}_{N_H} \end{bmatrix} \in 
	\mathbb{R}^{\sum\limits_{s=1}^{N_H}} n_{k_s}\vspace{\belowdisplayskip}
\end{equation}
and $\bm{1}$ is the $N_e$-dimensional vector of ones. Then, the hyperreduction approximation (\ref{eqn:res2}) can be written in matrix form as
\begin{equation*}
	\bm{C}\bm{\xi} - \bm{d} \approx 0
\end{equation*}
where $\bm{\xi}\in\mathbb{R}^{N_e}$ is the vector of element weights with a zero in each entry corresponding to a mesh entity $e \in \mathcal{E} \setminus \widetilde{\mathcal{E}}$. This suggests
formulating the problem of mesh sampling and weighting as finding a minimal subset of mesh entities $\widetilde{\mathcal{E}}$ and the corresponding set of weights $\Xi_{\widetilde{\mathcal{E}}}$ for which 
the hyperreduction approximation (\ref{eqn:res2}) retains a specified level of accuracy when applied to the training data set -- that is,
\begin{equation}
	\label{eqn:l0}
	\begin{split}
		\left(\widetilde{\mathcal{E}}, \Xi_{\mathcal E}\right) \; = \; &\argmin \norm{\bm{\zeta}}_0 \\ 
		&\begin{split}
			\, \textnormal{subject to} \; &\norm{\bm{C}\bm{\zeta}-\bm{d}}_2 \leq \varepsilon \norm{\bm{d}}_2 \\ 
			&\bm{\zeta} \geq 0 \,.
		\end{split}
	\end{split}
\end{equation}
In (\ref{eqn:l0}) above , $\norm{\diamond}_0$ denotes the $\ell^0$-``norm'' which counts the number of non-zero entries in $\diamond$ and $\varepsilon \in [0,1]$ is a small, user-specified, relative error 
tolerance for controlling the accuracy of the hyperreduction approximation.

Unfortunately, the optimization problem underlying (\ref{eqn:l0}) is well known to be NP-hard and thus is computationally intractable, even for mesh sizes that are small by computational mechanics 
standards. For this reason, it was proposed in \cite{chapman2017} to solve instead a convex approximation of this problem that promotes sparsity in the solution. For this purpose, several choices were 
identified and discussed in \cite{chapman2017}, including the following non-negative least-squares (NNLS) formulation
\begin{equation}
	\label{eqn:nnls}
	\begin{split}
		\textnormal{minimize} \; &\norm{\bm{C}\bm{\zeta}-\bm{d}}_2^2 \\
		\textnormal{subject to} \; &\bm{\zeta} \geq 0
	\end{split}
\end{equation}
equipped with the threshold-based early termination criterion
\begin{equation}
	\label{eqn:eps}
	\norm{\bm{C}\bm{\zeta}-\bm{d}}_2 \leq \varepsilon \norm{\bm{d}}_2 \,.
\end{equation}
This criterion replaces the Karush-Kuhn-Tucker optimality conditions and is necessary for promoting sparsity in the solution. Here, the convex optimization problem (\ref{eqn:nnls}) is chosen
as the alternative for the NP-hard problem (\ref{eqn:l0}) and solved using a distributed implementation of the active set algorithm proposed in \cite{lawson1995}. This implementation, which was
developed in \cite{chapman2017}, features a parallel, updatable, QR factorization kernel for efficiently solving at each iteration the overdetermined system of equations.

In summary, mesh sampling and weighting is performed here by solving offline the optimization problem (\ref{eqn:nnls}), as in the original ECSW hyperreduction method for Galerkin PROMs 
\cite{farhat2014, farhat2015}. The solution of this problem delivers simultaneously the subset of mesh entities $\widetilde{\mathcal{E}}$ and the associated set of weights 
$\Xi_{\widetilde{\mathcal{E}}}$. The resulting ECSW hyperreduction method for PG PROMs differentiates itself from alternatives such as EIM \cite{barrault2004}, DEIM \cite{chaturantabut2010}, 
GNAT \cite{carlberg2011, carlberg2013}, and other hyperreduction methods of the approximate-then-project type in at least two major ways:
\begin{itemize}
	\item From an algorithmic viewpoint, it determines the pair of sets $\left(\widetilde{\mathcal{E}}, \Xi_{\widetilde{\mathcal{E}}}\right)$ by efficiently solving a rigorous, convex approximation
		of the true sampling problem (\ref{eqn:l0}). The aforementioned alternatives rely on a heuristic, suboptimal greedy procedure to sample the original mesh (associated with the HDM).
	\item From a practical viewpoint, it has only one tuning parameter -- that is, the tolerance $\varepsilon$ whose value can be chosen to trade the sparsity of $\widetilde{\mathcal{E}}$ for the 
		accuracy of the resulting hyperreduction approximation. In general, the aforementioned alternatives depend on multiple tunable parameters. In their simplest form, many of them require 
		specifying {\it a priori} the size of the sampled set of mesh entities -- that is, their own equivalent of $\widetilde{\mathcal{E}}$ -- which is impractical as it typically requires a less 
		intuitive trial and error approach for determining a mesh sample that delivers the desired level of hyperreduction accuracy.
\end{itemize}

As stated in the introduction of this paper, a general concern for hyperreduction is the computational cost of its offline phase. For the hyperreduction method presented so far, this cost is dominated
by the computational cost of solving problem (\ref{eqn:nnls}), which itself is determined by the size of the matrix $\bm{C}$ characterizing this problem. The size of $\bm{C}$ is determined by 
$\sum\limits_{s=1}^{N_H} n_{k_s}$, which depends on the number of training solution snapshots and the average dimension of the local subspace approximation, and by the number of elements $N_e$ of the 
original mesh (associated with the HDM). For structural dynamics applications, FE models with $N_e = O(10^5)$, a number of training solution snapshots $N_H = O(10^2)$, and nonlinear Galerkin PROMs, it 
was shown in \cite{chapman2017} that the parallel implementation of the NNLS solver developed in that reference is scalable on up to 64 cores of a Linux cluster and solves problem (\ref{eqn:nnls}) on 
such a computing system in roughly 2 \si{\minute} wall-clock time. Another significant contribution of this paper is to show that even for CFD problems with $N_e = O(10^8)$, which is typical for 
industrial-scale applications, and a number of training solution snapshots $N_H = O(10)$, the matrix $\bm{C}$ can be built offline on a computing system with 3,584 cores in less than 15 \si{\minute} and 
the aforementioned NNLS solver can compute on this system in roughly 16 \si{\minute} wall-clock time a pair of sets $\left(\widetilde{\mathcal{E}}, \Xi_{\widetilde{\mathcal{E}}}\right)$ for which the 
hyperreduction approximation performed on nonlinear PG PROMs is very accurate.

\begin{remark}
	When hyperreducing a fully-discrete residual as in (\ref{eqn:res2}), the training solution snapshots chosen to construct offline the training matrix $\bm{C}$ and training vector 
	$\bm{d}$ (\ref{eqn:Cd}) should be collected such that each {\it training residual} $\bm{r}^m\left(\bm{u}_{o,k_s}+\Pi^{\bot}_{\bm{V}_{k_s}}(\bm{u}^{(s)}-\bm{u}_{o,k_s}),t^m; \bm{\mu}^q\right)$ is 
	sufficiently large in magnitude. Hence, a recommended strategy is to collect a solution snapshot for $\mathcal{S}_H$ at the beginning of a time-step $t^m$, where the Newton 
	iteration is typically initialized using the converged solution at the previous time-step, $\bm{u}^{m-1}$ -- and therefore where $\norm{\bm{r}^m}_2$ is the largest during that time-step.
	Collecting instead a training solution snapshot at the end of a time-step, after the Newton iterations during that time-step have converged, would lead to near-zero training residuals.
	This in turn would lead to a poor training of the pair $\left(\widetilde{\mathcal{E}}, \Xi_{\widetilde{\mathcal{E}}}\right)$ due to very small amplitude and therefore nearly redundant data.
\end{remark}

\subsection{Construction of a reduced mesh}
\label{sec:sampling2}

From (\ref{eqn:res2}) and (\ref{eqn:jac2}), it follows that as far as mesh entities are concerned, evaluating the hyperreduced residual vectors and Jacobian matrices associated with PG PROMs using
the ECSW method described above requires access only to $\widetilde{\mathcal{E}}$ and to the neighbors of each mesh entity $e\in\widetilde{\mathcal{E}}$ defining at $e$ the spatial discretization stencil 
underlying the HDM. Hence, it is convenient to support such evaluations with a ``reduced mesh'' defined by the augmented set $\widetilde{\mathcal{E}}^+ \subset \mathcal{E}$ consisting of the mesh
entities underlying the non-zero entries of the Boolean localization matrix $\bm{L}_{e^+}$ (hence, $\widetilde{\mathcal{E}} \subset \widetilde{\mathcal{E}}^+$). This is because just like the HDM-based 
residuals and Jacobians can be conveniently computed on the mesh associated with the HDM, their hyperreduced counterparts (\ref{eqn:res2}) and (\ref{eqn:jac2}) can be conveniently computed on the 
reduced mesh. Due to the locality properties of most spatial discretization schemes, $\abs{\widetilde{\mathcal{E}}} \ll N_e \Rightarrow \abs{\widetilde{\mathcal{E}}^+} \ll N_e$, 
which means that the online evaluations of the hyperreduction approximations (\ref{eqn:res2}) and (\ref{eqn:jac2}) can be efficiently performed on the reduced mesh.

\section{Applications}
\label{sec:app}

Here, the proposed adaptation to PG PROMs of the ECSW hyperreduction method and its performance are illustrated and assessed, respectively, for three CFD applications:
\begin{itemize}
	\item An academic problem of laminar flow over a circular cylinder that has the added benefit of being easily reproducible by the reader. For this problem, the performance 
		of the proposed hyperreduction method is compared to that of two alternative hyperreduction methods for PG PROMs: the GNAT method \cite{carlberg2011, carlberg2013}; and 
		the least-squares collocation method \cite{legresley2006} equipped with the same reduced mesh produced by the gappy-POD-based GNAT method.
	\item An unsteady wake flow problem associated with the Ahmed body geometry -- which is considered by the automotive industry to be a benchmark problem for CFD.
	\item A very large-scale, turbulent CFD problem of industrial relevance -- specifically, the prediction of the turbulent flow past an F-16C/D Block 40 aircraft configuration with external stores
		at a high angle of attack.
\end{itemize}
In particular, the second and third applications outlined above demonstrate the combination of local subspace approximations for mitigating the well-known Kolmogorov $n$-width barrier for the PMOR of first-order
nonlinear hyperbolic problems such as those associated with convection-dominated CFD problems \cite{ohlberger2016} and the proposed ECSW hyperreduction method for achieving computational efficiency.

In each application introduced above, the HDM-based problem (\ref{eqn:hdm0}) is constructed by semi-discretizing the nondimensional form of the three-dimensional, compressible, 
Navier-Stokes equations by a mixed FV/FE method. In this method, the spatial approximation of the convective fluxes is performed using an asymptotically third-order, upwind, vertex-based FV scheme 
based on the MUSCL-type approach and Roe's approximate Riemann solver \cite{anderson1986, farhat1993}, while the diffusive term and all source terms are treated by a piecewise-linear, Galerkin FE scheme.
In the second and third applications, turbulence modeling is performed using the detached-eddy simulation (DES) approach based on the Spalart-Allmaras one-equation model \cite{strelets2001}.
The resulting nonlinear flux vector $\bm{f}\left(\bm{u}(t; \bm{\mu}); \bm{\mu}\right)$ \big(see (\ref{eqn:hdm})\big) is not characterized by a polynomial dependence on the high-dimensional solution
vector $\bm{u}(t; \bm{\mu})$: for this reason, hyperreduction is required to reduce the computational complexity of each constructed PROM. 

It is noted that for each application highlighted above, the PROM and HPROM are constructed for a single parameter instance $\bm{\mu}^{\star}$ of the HDM, for two different reasons: 1) the highly 
nonlinear, time-dependent nature of the flow problem already necessitates hyperreduction due to the computational bottlenecks discussed in Section \ref{sec:pg2}; and 2) the sheer size of the HDM --
of the order of $10^7$ to $10^8$ DOFs -- is such that considering instead a parametric instance of the flow problem would lead to an unnecessary consumption of precious CPU resources without adding 
significant value to the performance assessment intended here.

In each application: 
\begin{itemize}
	\item Time-discretization of the HDM is performed using an implicit scheme and the nonlinear system of equations arising at each computational time-step is solved by a Newton-Krylov 
		method equipped with an additive Schwarz preconditioned GMRES \cite{cai1998} algorithm as the linear equation solver.  
	\item The PG PROM is constructed using the LSPG method \cite{carlberg2011, carlberg2013}. Hence, $\bm{W}_k$ is defined at each time-instance $t^m$ by 
		\begin{equation*} 
			\bm{W}_k=\bm{W}_k^m\left(\bm{y}^m_k(\bm{\mu}); \bm{\mu}\right) = \bm{J}^m\left(\bm{u}_{o,k}+\bm{V}_k\bm{y}^m_k(\bm{\mu}); \bm{\mu}\right)\bm{V}_k\in\mathbb{R}^{N\times n_k} 
		\end{equation*} 
		where the HDM-based Jacobian matrix $\bm{J}^m\left(\bm{u}_{o,k}+\bm{V}_k\bm{y}_k^m(\bm{\mu}); \bm{\mu}\right)$ is given in (\ref{eqn:hdmjac}). 
	\item The parallel variant of the Lawson and Hanson NNLS algorithm based on an updatable QR factorization developed in \cite{chapman2017} is used to solve the convex optimization
		problem (\ref{eqn:nnls}) equipped with (\ref{eqn:eps}) and $\varepsilon = 1\times10^{-2}$, and produce the reduced mesh $\widetilde{\mathcal{E}}$ and associated set of weights 
		$\Xi_{\widetilde{\mathcal{E}}}$.
	\item The accuracy of the constructed LSPG HPROM -- that is, the accuracy of the approximate solution of the HDM-based problem (\ref{eqn:hdm0}) obtained by reconstructing the solution of the 
		PROM-based problem (\ref{eqn:prom}) equipped with the hyperreduction approximations (\ref{eqn:res2}) and (\ref{eqn:jac2}), or an alternative hyperreduction approximation method -- is 
		assessed for a selection of quantities of interest (QoIs) using in each case the relative error with respect to the counterpart value obtained from the solution of the HDM-based 
		problem 
		\begin{equation} 
			\label{eq:RE}
			\mathbb{RE}_Q = \frac{\sqrt{\sum\limits_{t\,\in\,\mathcal{P}}\left(\widetilde{Q}(t)-Q(t)\right)^2}}{\sqrt{\sum\limits_{t \, \in \, 
			\mathcal{P}} Q(t)^2}} \times 100\,(\%)
		\end{equation}
		where $\widetilde{Q}(t)$ is the approximation of the QoI $Q(t)$ obtained by reconstructing the solution of the PROM-based problem (\ref{eqn:prom}) equipped with the hyperreductions 
		(\ref{eqn:res2}) and (\ref{eqn:jac2}) and employing the same postprocessing procedures as for the HDM-based solution, $\mathcal{P}$ is the set of time-stamps used for evaluating $\mathbb{RE}_Q$ -- that is
		\begin{equation*}
			\mathcal{P} = \left \{t \in \{0, \Delta s_{\mathbb{RE}}, 2\Delta s_{\mathbb{RE}}\,...\} : t \leq T_f \right\} \,,
		\end{equation*}
		and $\Delta s_{\mathbb{RE}}$ denotes the sampling time-interval for the computation of $\mathbb{RE}_Q$.
	\item The performance of the LSPG HPROM is assessed by computing the speedup factors it delivers with respect to both the CPU time and wall-clock time elapsed in the solution of the HDM-based 
		problem.  For this purpose, each speedup factor is defined here as the ratio of the CPU (wall-clock) time elapsed in the solution of the LSPG HPROM-based problem and the CPU (wall-clock) 
		time elapsed in the solution of the underlying HDM problem.
\end{itemize}
Finally, it is specified that all numerical simulations reported herein are performed in double precision arithmetic.

\subsection{Laminar flow over a circular cylinder}
\label{sec:cylinder}

\subsubsection{High-dimensional model}

The first application considered here focuses on the computation of a two-dimensional, laminar flow over a right circular cylinder at the Reynolds number $Re = 100$ and the free-stream Mach number 
$M_\infty=0.2$. At this Reynolds number, the flow exhibits periodic vortex shedding after a transient startup phase: it demonstrates well the von K\'arm\'an vortex street. This problem is a commonly 
studied academic problem and therefore has the added benefit of being easily reproducible by the reader. 

\begin{figure}[h!]
	\centering
	%\vglue 0.1 truein
	\begin{subfigure}[c]{0.25\textwidth}%
	\centering
	\includegraphics[width=\textwidth]{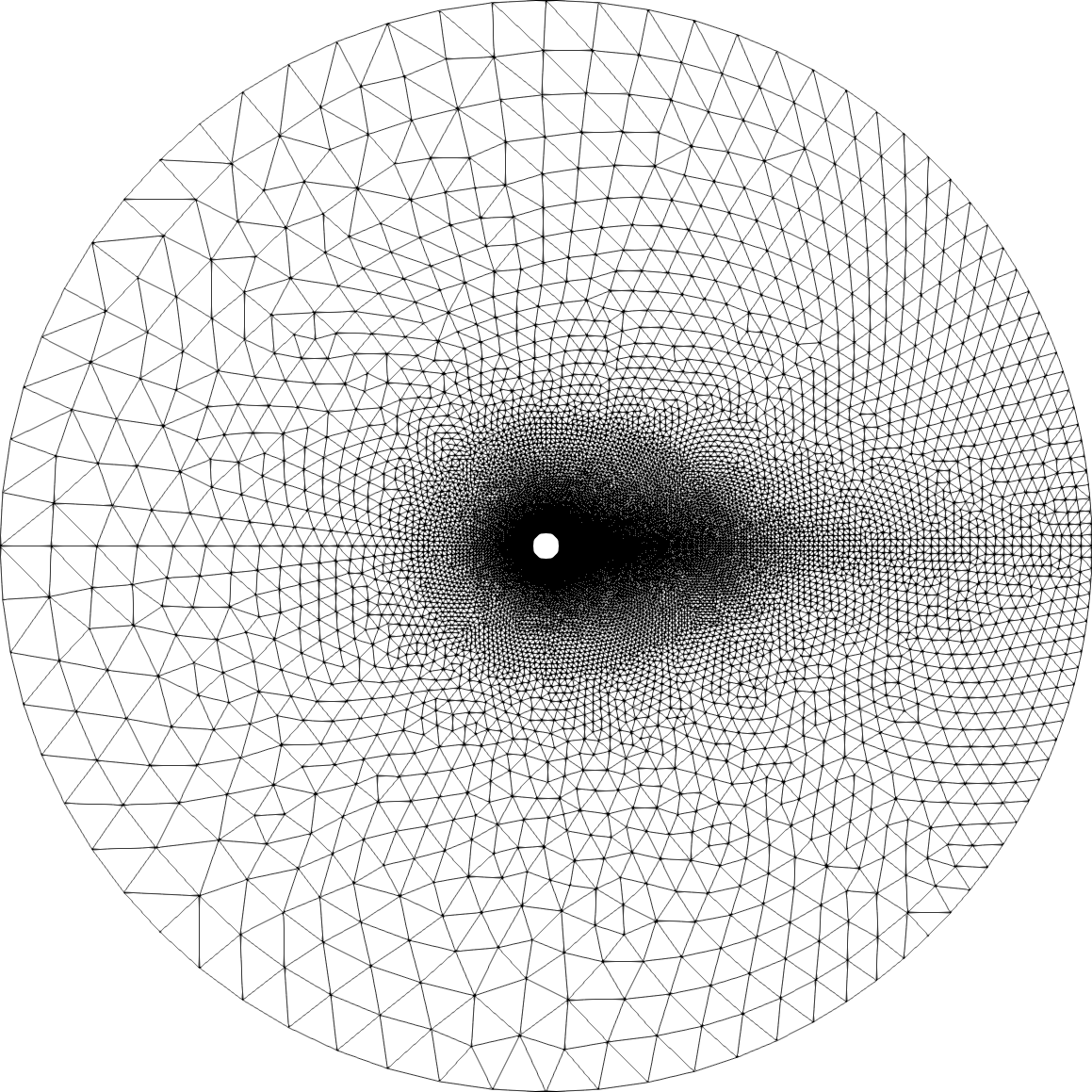}
	\caption{}
	\label{fig:cylindermesha}
	\end{subfigure}%
	\hspace{1em}
	\begin{subfigure}[c]{0.25\textwidth}%
	\centering
	\includegraphics[width=\textwidth]{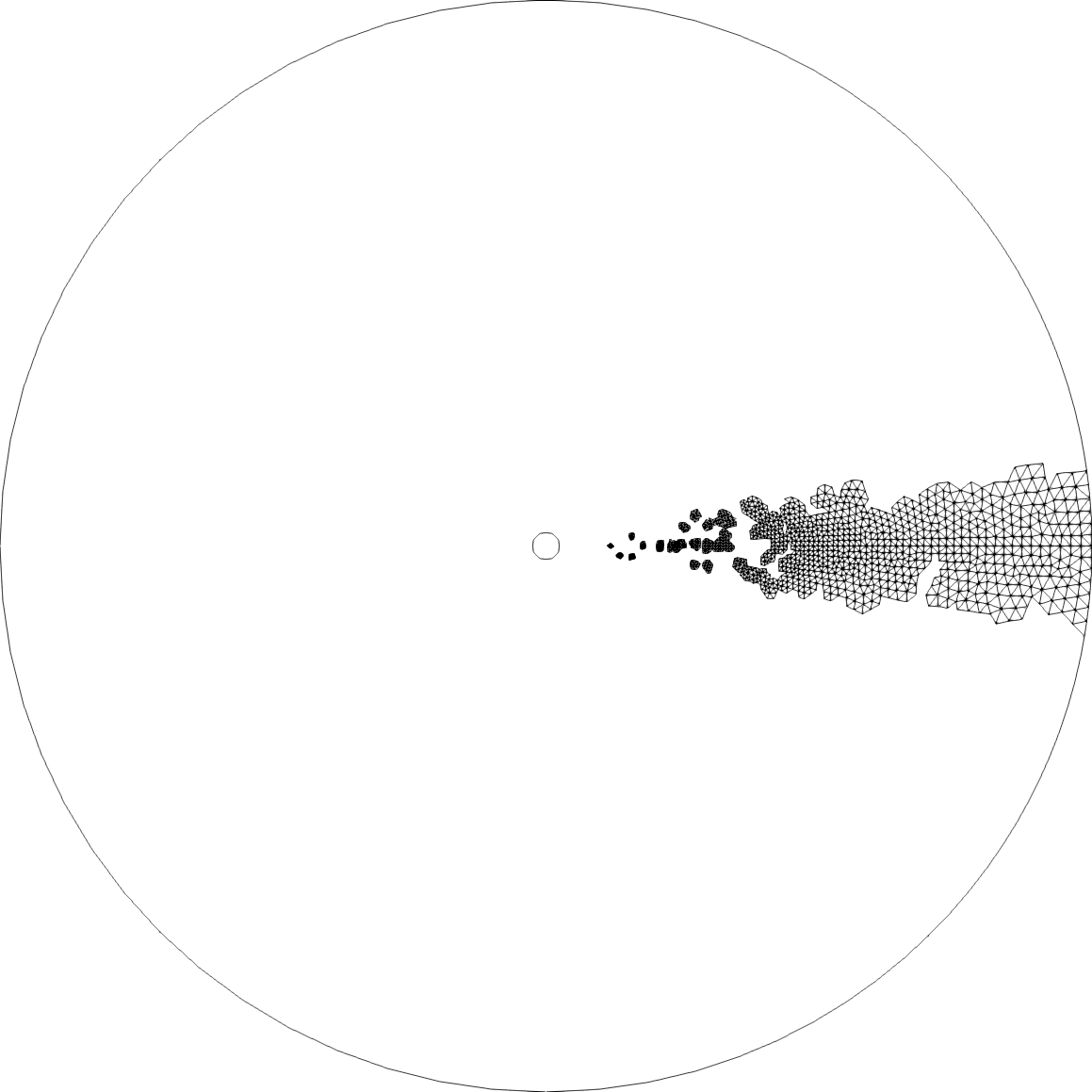}
	\caption{}
	\label{fig:cylindermeshb}
	\end{subfigure}%
	\hspace{1em}
	\begin{subfigure}[c]{0.25\textwidth}%
	\centering
	\includegraphics[width=\textwidth]{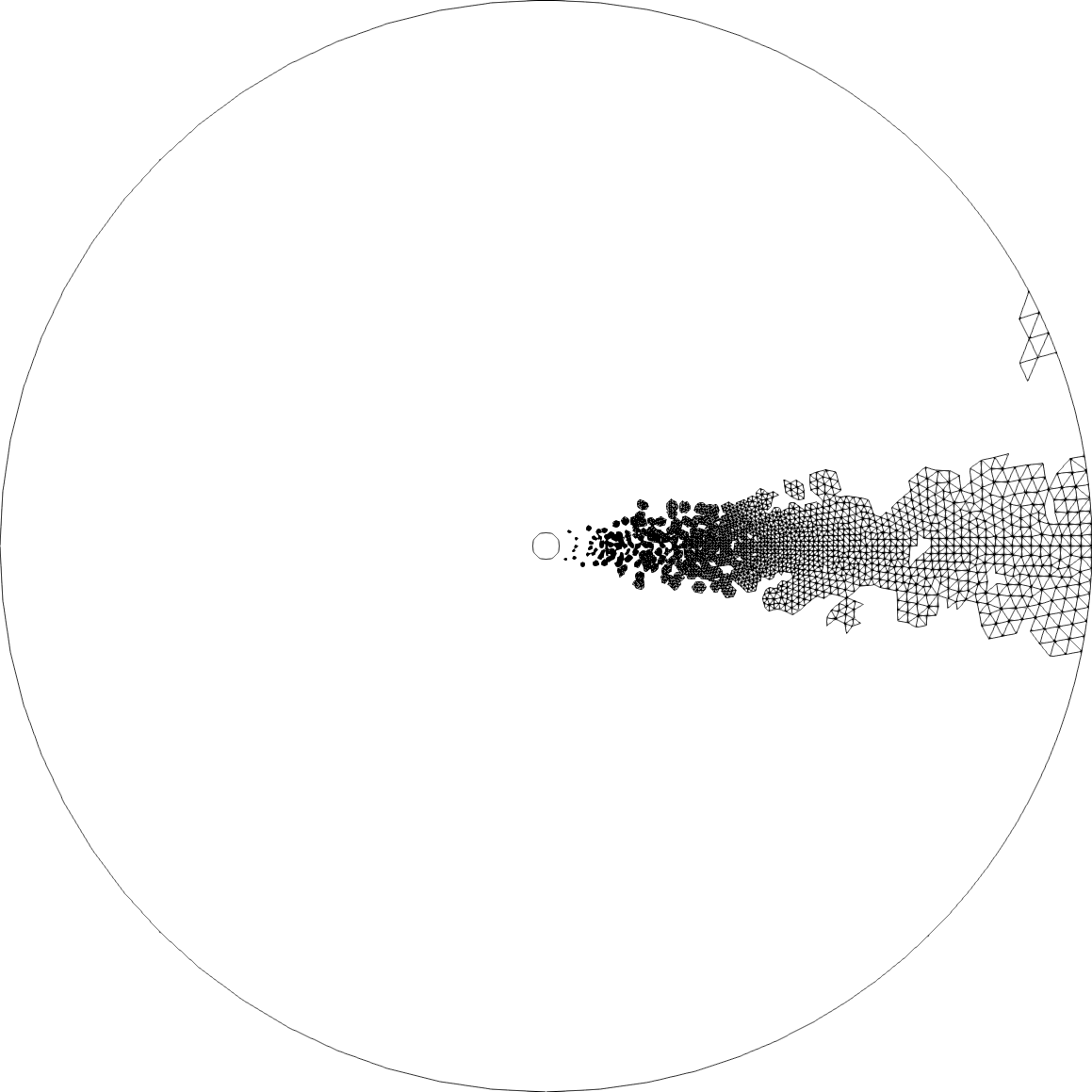}
	\caption{}
	\label{fig:cylindermeshc}
	\end{subfigure}%
	%\vglue 0.1 truein
	\caption{Two-dimensional flow problem over a right circular cylinder: discretized computational domain (a); reduced mesh constructed using the proposed ECSW method (b); and reduced mesh 
	constructed using the gappy-POD-based GNAT method (c).}
	\label{fig:cylindermesh}
\end{figure}

For this application, the computational domain is chosen as the disk of diameter $40D$, where $D$ denotes the cylinder diameter. This domain is discretized using a one-element-thick, unstructured,
three-dimensional mesh with $98,140$ vertices and $284,700$ tetrahedral elements. Figure \ref{fig:cylindermesha} shows the computational mesh.  Symmetry boundary conditions are applied on its spanwise 
faces to ensure that the resulting flow is two-dimensional. A no-slip adiabatic wall boundary condition is applied on the cylinder surface. The dimension of the resulting semi-discrete HDM is $N=490,700$.
This HDM is time-discretized using a second-order DIRK scheme and a fixed nondimensional time-step $\Delta t=1\times10^{-1}$, which, for the aforementioned CFD mesh, corresponds to a CFL number of 
approximately $7,600$. The initial condition for the HDM-based simulation is computed by impulsively starting the flow from a uniform state and time-integrating the semi-discrete HDM until the onset of 
vortex shedding. Then, the semi-discrete HDM is time-integrated from this initial condition until the end of the nondimensional time-interval $[0,200]$. At roughtly $t=100$, the flow becomes periodic.
The computed flow solution using the discrete HDM is in good agreement with the results of experimental and other numerical studies of this problem \cite{wieselsberger1922, henderson1995}.

\subsubsection{Hyperreduced LSPG projection-based reduced-order models}

In the first nondimensional time subinterval $[0, 150]$, the HDM-based solution is sampled at the nondimensional rate defined by $\Delta s=2\times10^{-1}$ (recall that after
$t = 150$, the solution remains periodic). Two different {\it global} ($N_c = 1$) affine subspace approximations are constructed by compressing the $N_s=751$ solution snapshots collected
in $[0, 150]$ using SVD and building two right ROBs of dimension $n=9$ and $n=35$. These two right ROBs correspond to the singular value energy truncation tolerances of $99\%$ and $99.99\%$, 
respectively. In both cases, the affine offset $\bm{u}_{o}$ is set to the initial condition of the problem. Then, the ECSW method described in this paper is applied to construct for each aforementioned 
subspace approximation an associated reduced mesh using $N_H=376$ training residual snapshots and the corresponding LSPG HPROM. The training residual snapshots are computed from the sampling of the 
HDM-based solution in $[0, 150]$ at the nondimensional rate $\Delta s_H=4\times10^{-1}$. The size $\widetilde{N}_e$ of each constructed reduced mesh is given in Table \ref{tab:cylinder}. The reader can 
observe that this size is a tiny fraction of the size $N_e$ of the CFD mesh outlined above.

\begin{table}[h!]
	\small
	\centering
	\caption{Two-dimensional flow problem over a right circular cylinder: number of mesh cells sampled by ECSW for each constructed global subspace approximation ($\epsilon = 1\times 10^{-2}$).}
	%\vglue 0.1 truein
	\begin{tabular}{ccc} 
		\toprule
		$n$ & $\widetilde{N}_e$ & $\widetilde{N}_e/N_e$ ($\%$) \\ \midrule
		$9$ & $75$ & $0.076$ \\
		$35$ & $321$ & $0.327$ \\ \bottomrule
	\end{tabular}
	%\vglue 0.1 truein
	\label{tab:cylinder}
\end{table}

Next, for each size $\widetilde{N}_e$ reported in Table \ref{tab:cylinder}, an additional reduced mesh is constructed using the gappy-POD-based GNAT method \cite{carlberg2011, carlberg2013}. This
additional reduced mesh is to be used in both the GNAT hyperreduction method and the least-squares collocation method (which in principle can be equipped with any reduced mesh).
Figure \ref{fig:cylindermeshb} shows the reduced mesh delivered by ECSW for $n=35$ and Figure \ref{fig:cylindermeshc} shows the alternative reduced mesh of the same size constructed using the 
alternative gappy-POD-based GNAT method. Using these reduced meshes, each constructed LSPG HPROM is discretized using the same second-order DIRK scheme applied to discretize the 
underlying HDM and the same nondimensional time-step $\Delta t=1\times10^{-1}$, and time-integrated on a single core of the same Linux cluster on which the HDM-based simulation is performed.

\subsubsection{Performance of various hyperreduction methods}

\begin{figure}[h!]
	\centering
	%\vglue 0.1 truein
	\begin{subfigure}[c]{0.45\textwidth}%
	\hspace{0.5em} \includegraphics[width=0.85\textwidth]{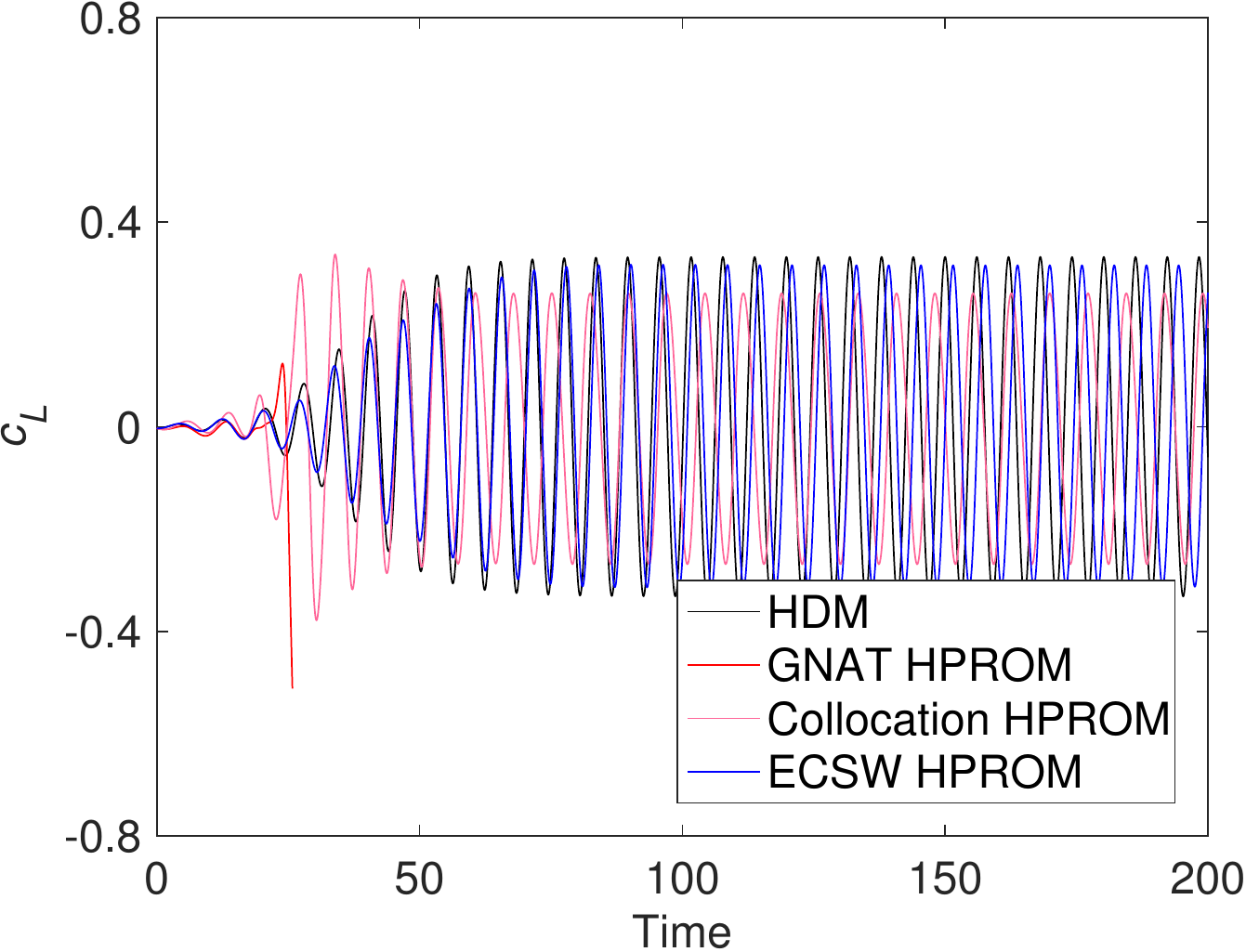}
	\caption{}
	\end{subfigure}%
	\hspace{1em}
	\begin{subfigure}[c]{0.45\textwidth}%
	\hspace{0.5em} \includegraphics[width=0.85\textwidth]{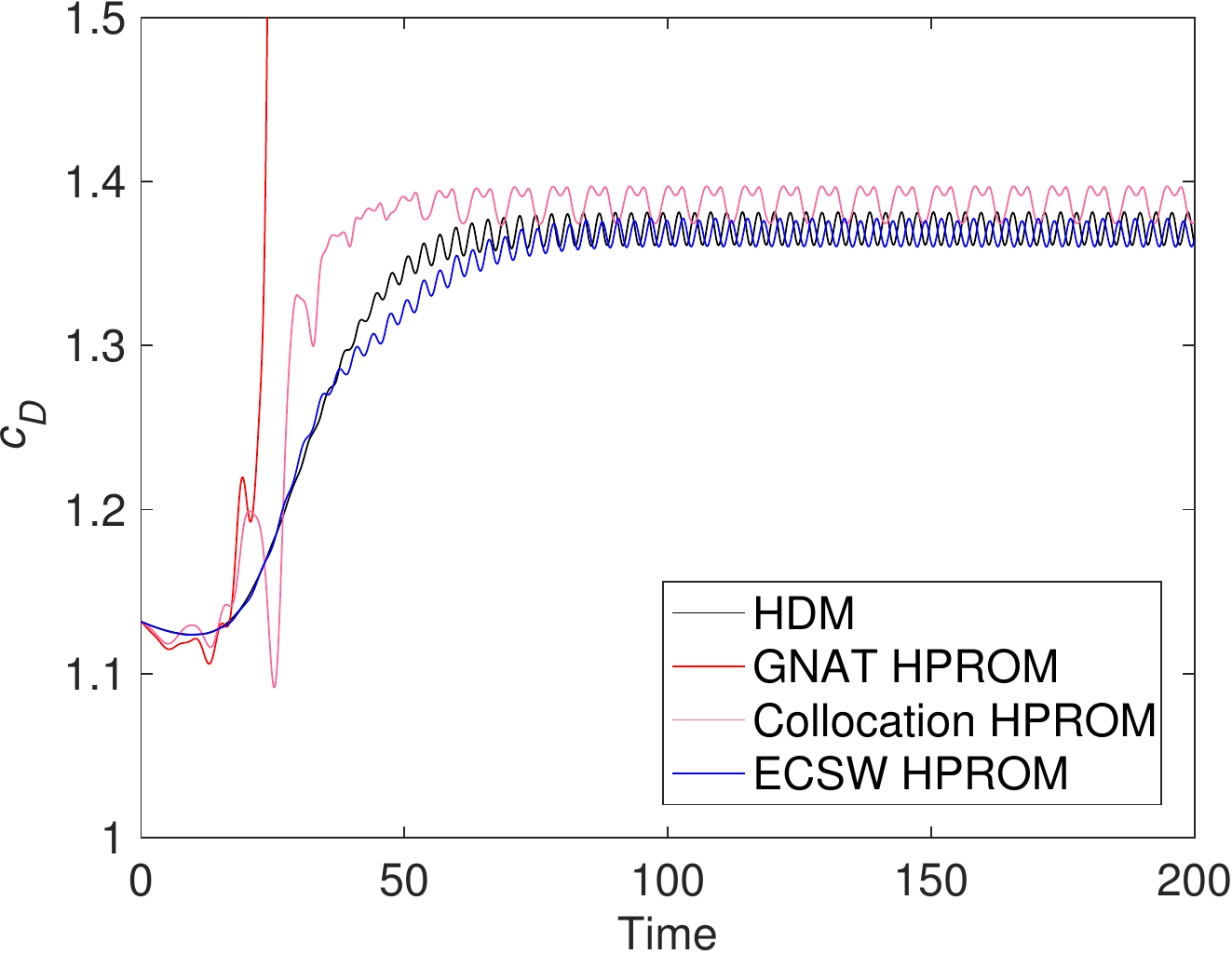}
	\caption{}
	\end{subfigure}%
	
	\medskip
	\begin{subfigure}[c]{0.45\textwidth}%
	\hspace{0.5em} \includegraphics[width=0.85\textwidth]{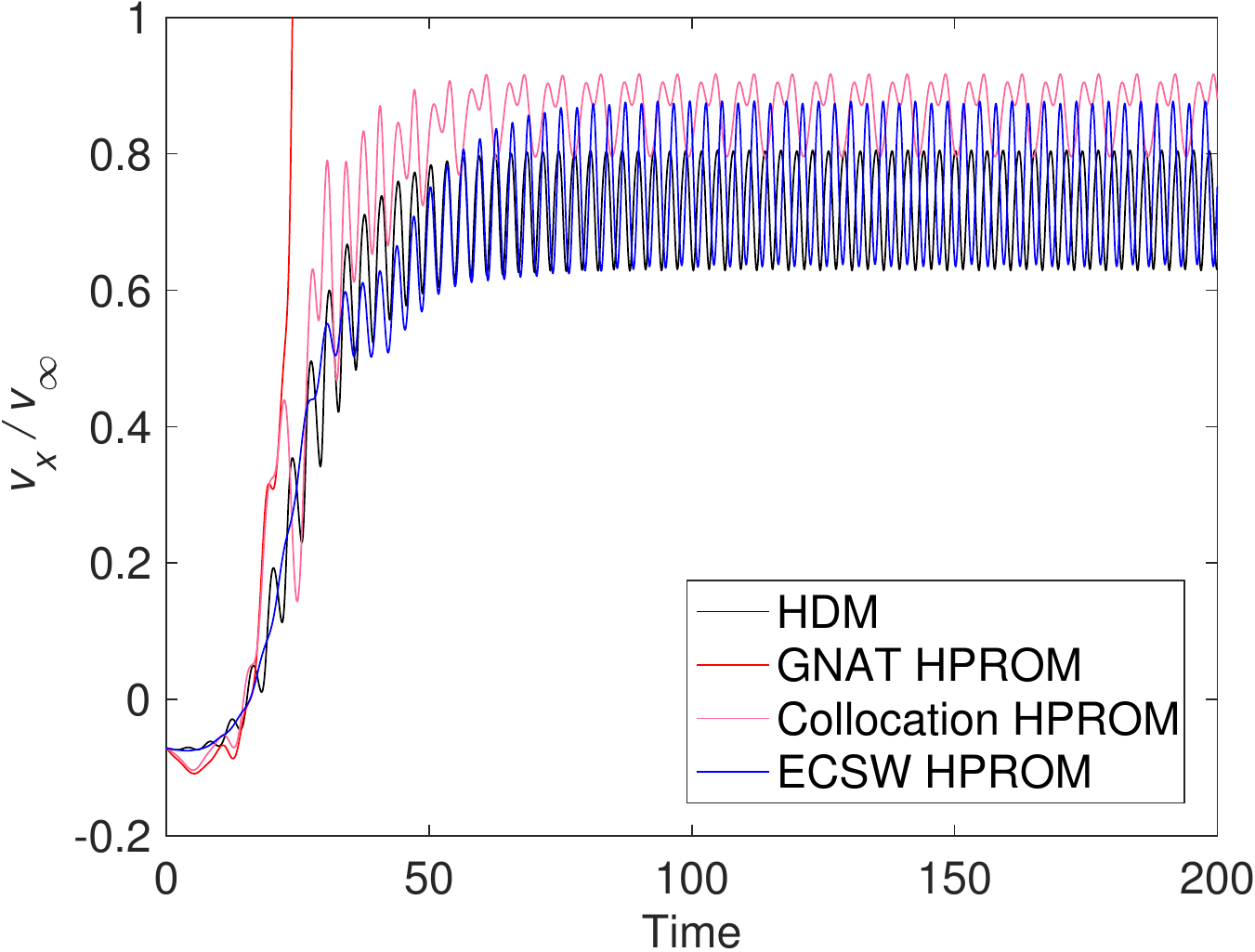}
	\caption{}
	\end{subfigure}%
	\hspace{1em}
	\begin{subfigure}[c]{0.45\textwidth}%
	\hspace{0.5em} \includegraphics[width=0.85\textwidth]{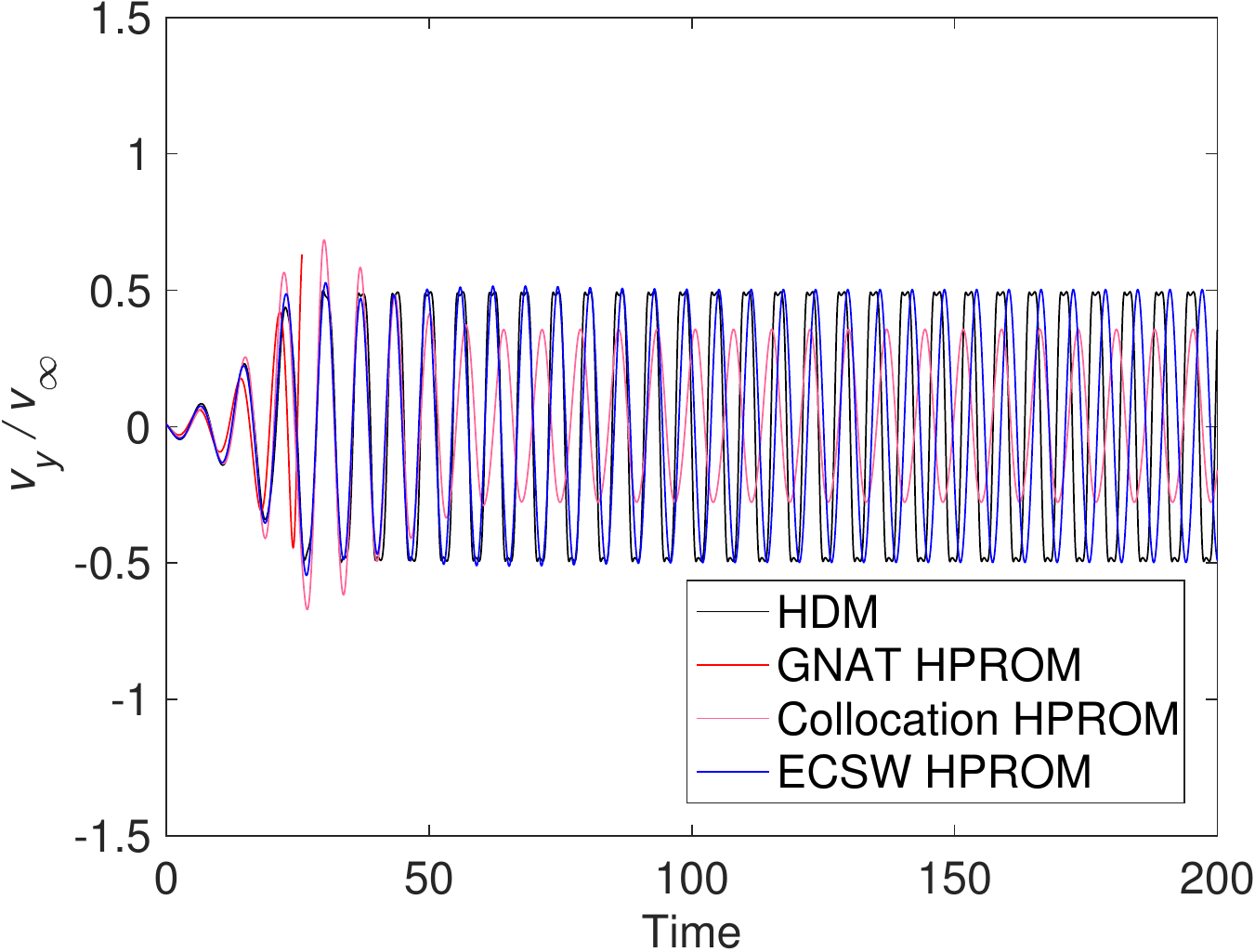}
	\caption{}
	\end{subfigure}%
	%\vglue 0.1 truein	
	\caption{Two-dimensional flow problem over a right circular cylinder -- time-histories computed using the HDM and various LSPG HPROMs of dimension $n=9$ of: the lift coefficient (a); drag 
	coefficient (b); streamwise velocity component at a probe (c); and the normal velocity component at the same probe (d).}
	\label{fig:cylinder99}
\end{figure}

For the case $n=9$, the time-histories of the lift coefficient, drag coefficient, and velocity components at a probe located $5D$ downstream from the cylinder's trailing edge computed using
the following computational models are reported in Figure \ref{fig:cylinder99}:
\begin{itemize}
		\item The HDM.
		\item The LSPG HPROM constructed using the ECSW hyperreduction method.
		\item The LSPG HPROM constructed using the GNAT hyperreduction method.
		\item The LSPG HPROM constructed using the least-squares collocation hyperreduction method equipped with the reduced mesh generated by the gappy-POD-based GNAT method.
\end{itemize}

\begin{figure}[h!]
	\centering
	%\vglue 0.1 truein
	\begin{subfigure}[c]{0.45\textwidth}%
	\hspace{0.5em} \includegraphics[width=0.85\textwidth]{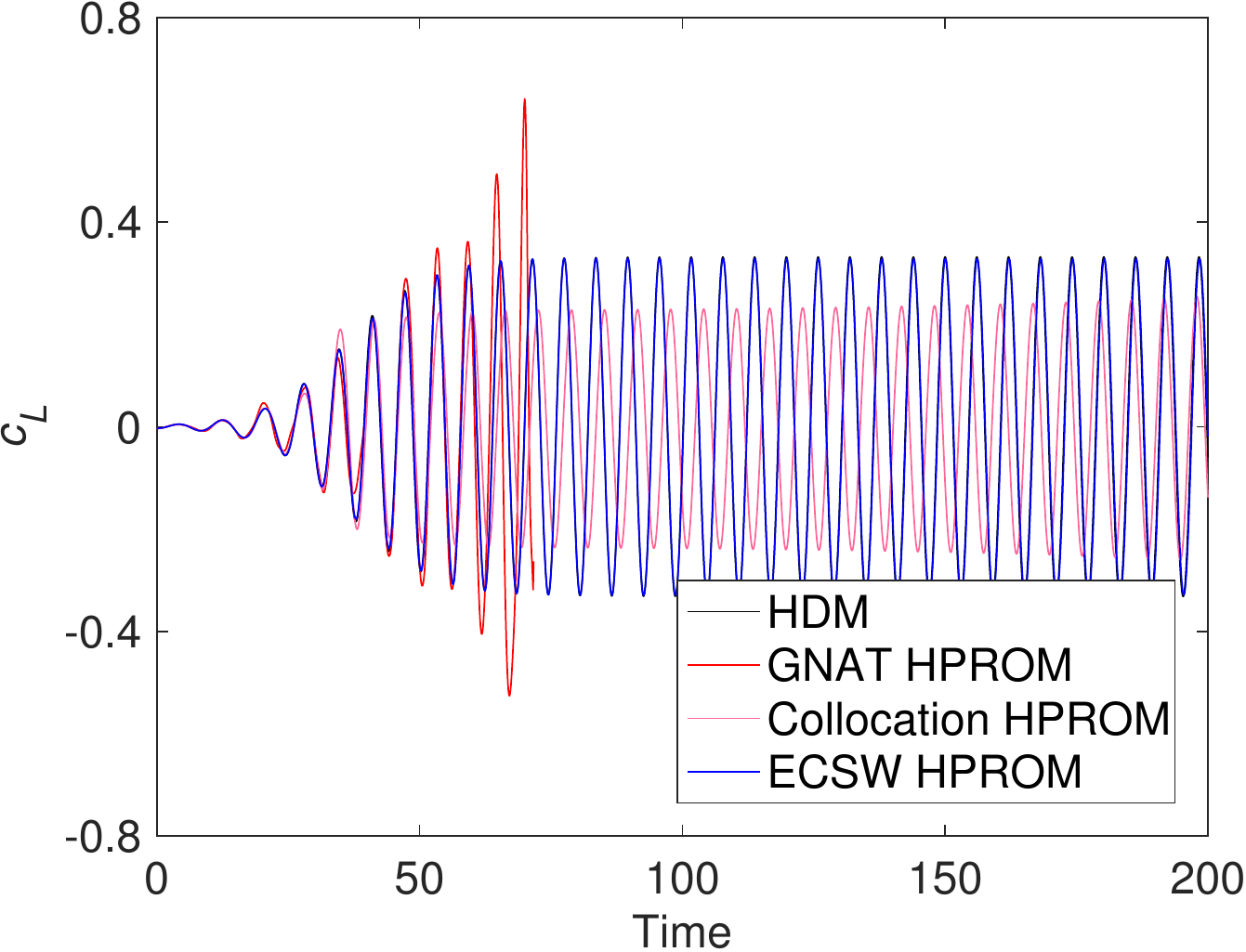}
	\caption{}
	\end{subfigure}%
	\hspace{1em}
	\begin{subfigure}[c]{0.45\textwidth}%
	\hspace{0.5em} \includegraphics[width=0.85\textwidth]{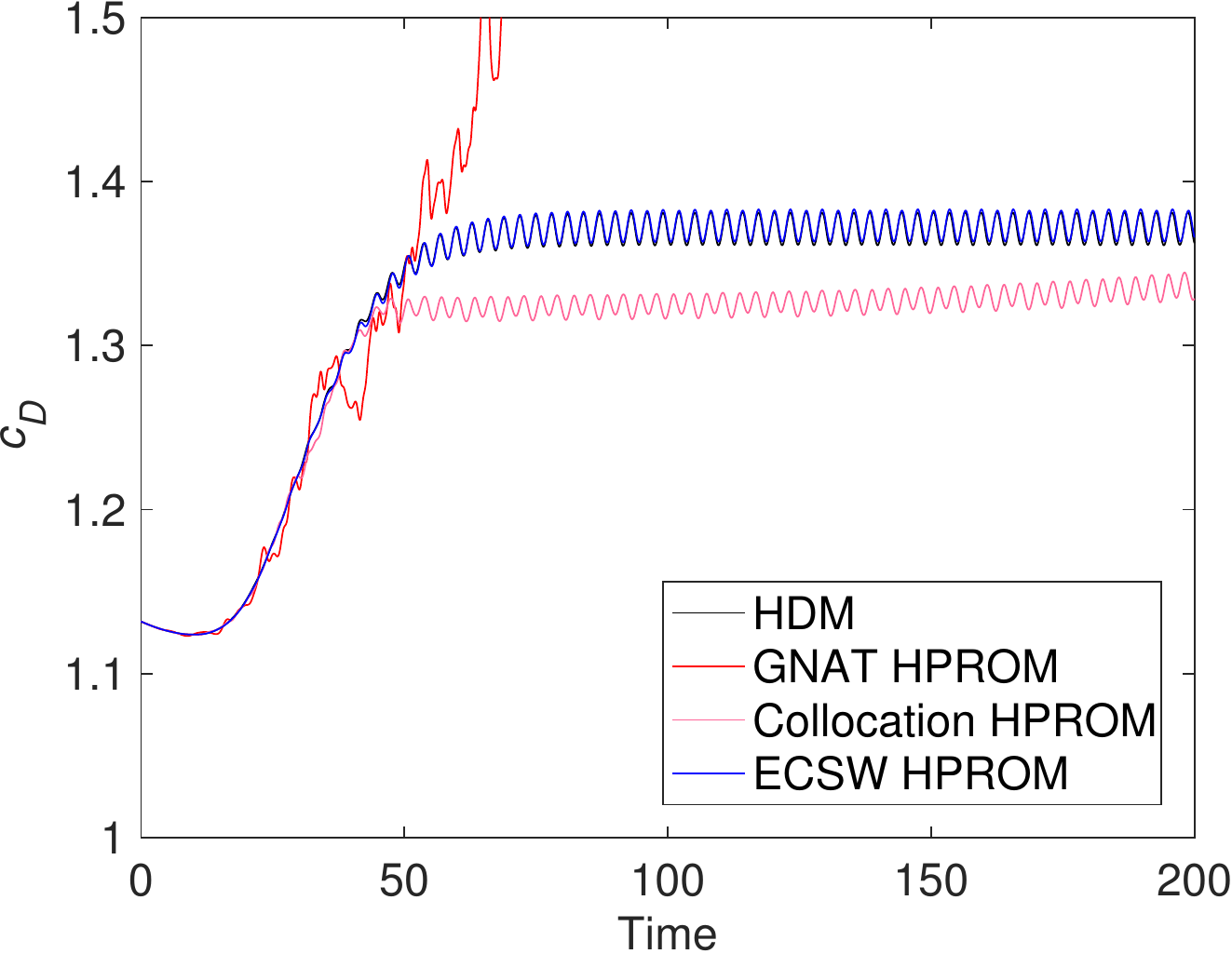}
	\caption{}
	\end{subfigure}%
	
	\medskip
	\begin{subfigure}[c]{0.45\textwidth}%
	\hspace{0.5em} \includegraphics[width=0.85\textwidth]{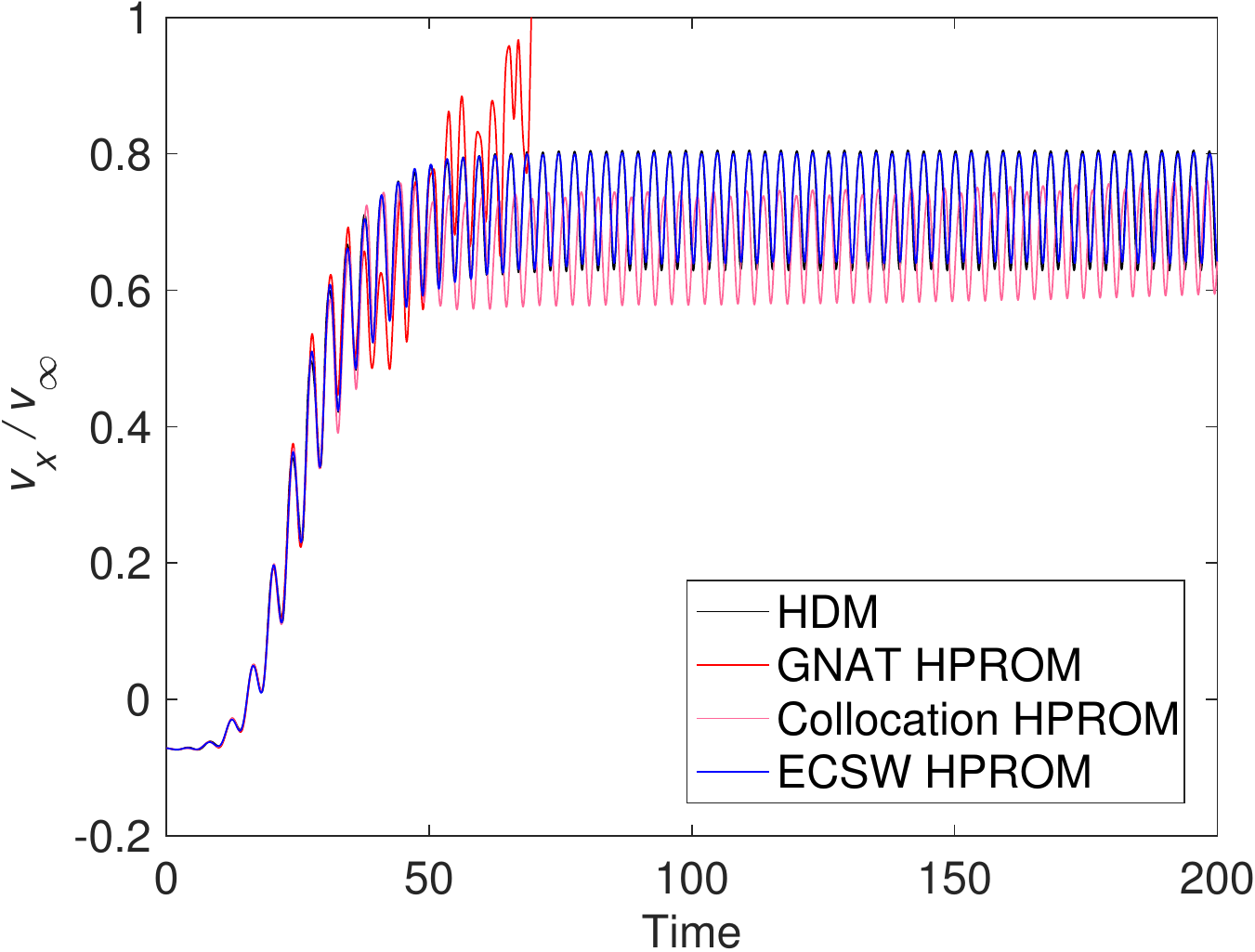}
	\caption{}
	\end{subfigure}%
	\hspace{1em}
	\begin{subfigure}[c]{0.45\textwidth}%
	\hspace{0.5em} \includegraphics[width=0.85\textwidth]{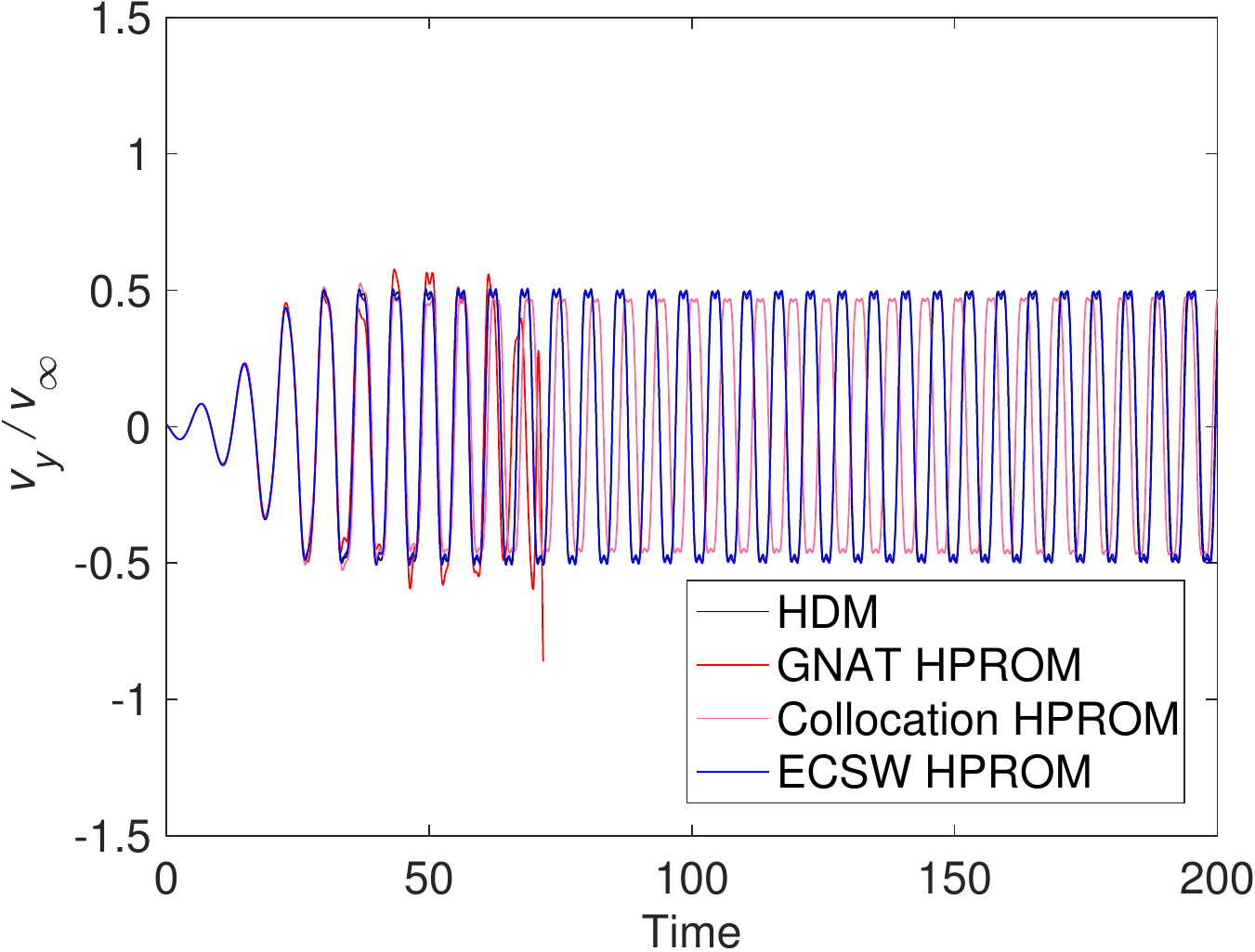}
	\caption{}
	\end{subfigure}%
	%\vglue 0.1 truein
	\caption{Two-dimensional flow problem over a right circular cylinder -- time-histories computed using the HDM and various LSPG HPROMs of dimension $n=35$ of: the lift coefficient (a); drag 
	coefficient (b); streamwise velocity component at a probe (c); and the normal velocity component at the same probe (d).}
	\label{fig:cylinder9999}
\end{figure}

Figure \ref{fig:cylinder9999} reports the counterpart results for the case $n=35$. The reader can observe that in both cases, the LSPG HPROM built using the GNAT method
is numerically unstable and fails to complete the simulation in the specified time-interval. In both cases, the counterpart LSPG HPROM built using the least-squares collocation method
is numerically stable, but fails to deliver an acceptable accuracy. On the other hand, the ECSW-based LSPG HPROM is numerically stable in both cases and delivers exceptional accuracy.

Table \ref{tab:cylindererr} quantifies for this problem the relative errors incurred for the chosen QoIs using $\Delta s_{\mathbb{RE}}=\Delta t=1\times10^{-1}$. The relative errors for the GNAT-based
LSPG HPROM are not included due to the numerical instability exhibited by this HPROM in the time-interval of interest. These errors demonstrate further the superior accuracy of the ECSW-based LSPG HPROM.

\begin{table}[h!]
	\small
	\centering
	\caption{Two-dimensional flow problem over a right circular cylinder: computational accuracy of the various constructed LSPG HPROMs.}
	%\vglue 0.1 truein
	\begin{tabular}{clcccc}
		\toprule
		$n$ & Model & $\mathbb{RE}_{c_D}$ ($\%$) & $\mathbb{RE}_{c_L}$ ($\%$) & $\mathbb{RE}_{v_x}$ ($\%$) & $\mathbb{RE}_{v_y}$ ($\%$) \\ \midrule
		$9$ & GNAT HPROM & $--$ & $--$ & $--$ & $--$ \\
		& Collocation HPROM & $113$ & $125$ & $94.4$ & $103$ \\
		& ECSW HPROM & $0.891$ & $2.45$ & $16.2$ & $22.3$ \\ \midrule
		$35$ & GNAT HPROM & $--$ & $--$ & $--$ & $--$ \\
		& Collocation HPROM & $6.94$ & $125$ & $5.84$ & $130.$ \\
		& ECSW HPROM & $0.130$ & $2.86$ & $1.65$ & $13.7$ \\ \bottomrule
	\end{tabular}
	%\vglue 0.1 truein
	\label{tab:cylindererr}
\end{table}

\begin{figure}[h!]
	\centering
	%\vglue 0.1 truein
	\begin{subfigure}[c]{0.475\textwidth}%
	\centering
	\includegraphics[width=\textwidth]{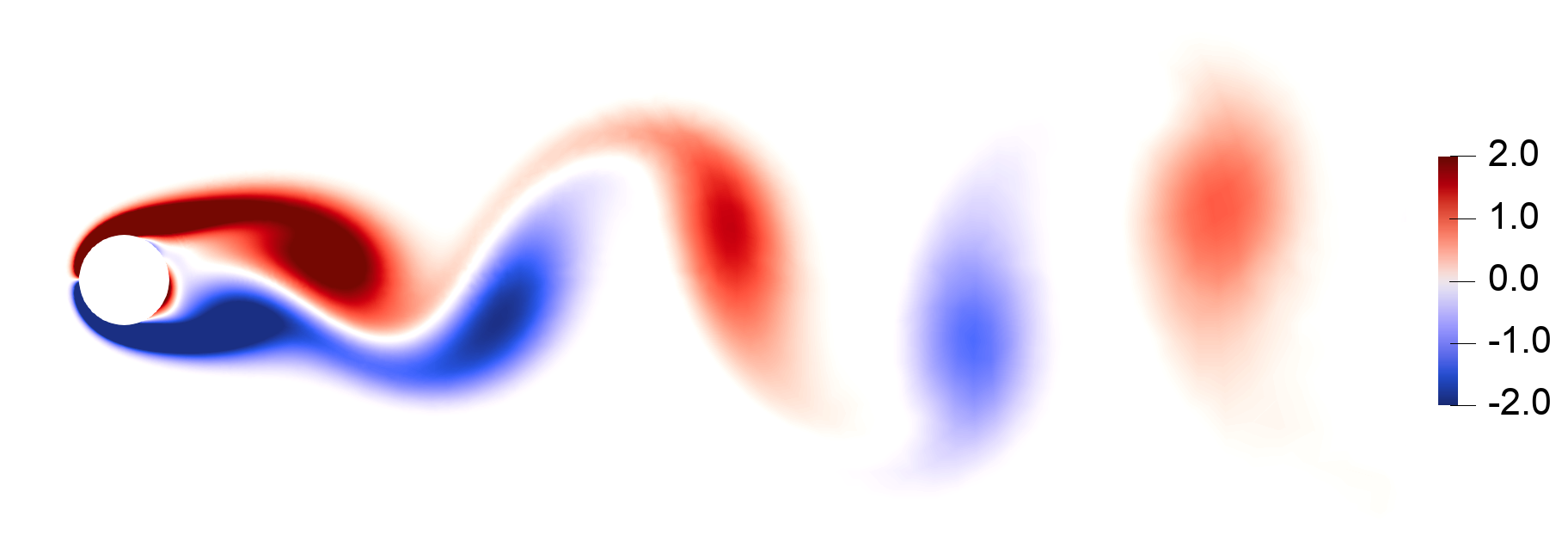}
	\caption{}
	\end{subfigure}%
	\hspace{1em}
	\begin{subfigure}[c]{0.475\textwidth}%
	\centering
	\includegraphics[width=\textwidth]{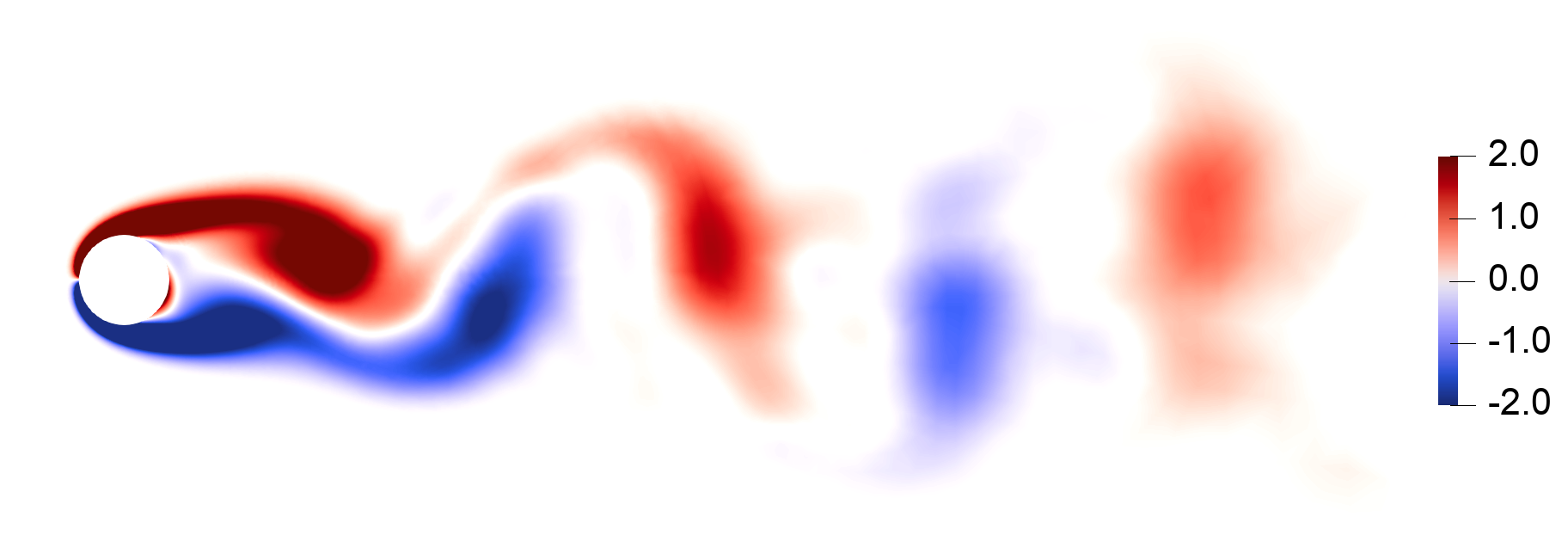}
	\caption{}
	\end{subfigure}%
	
	\medskip
	\begin{subfigure}[c]{0.475\textwidth}%
	\centering
	\includegraphics[width=\textwidth]{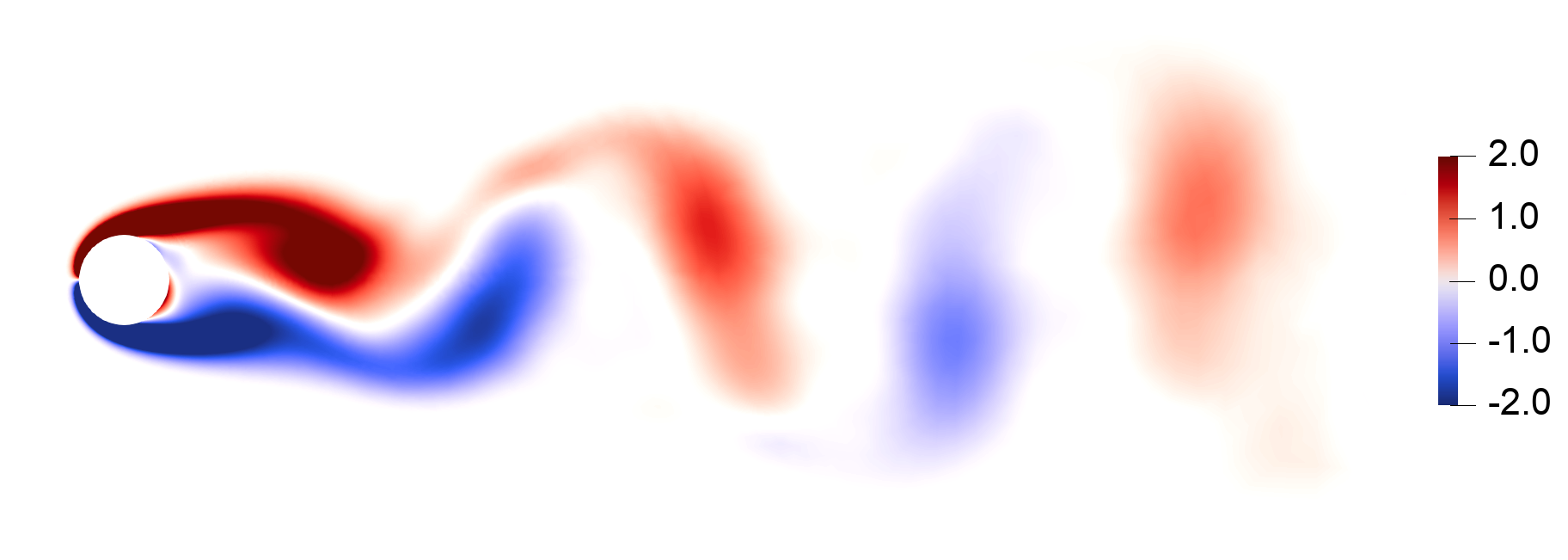}
	\caption{}
	\end{subfigure}%
	\hspace{1em}
	\begin{subfigure}[c]{0.475\textwidth}%
	\centering
	\includegraphics[width=\textwidth]{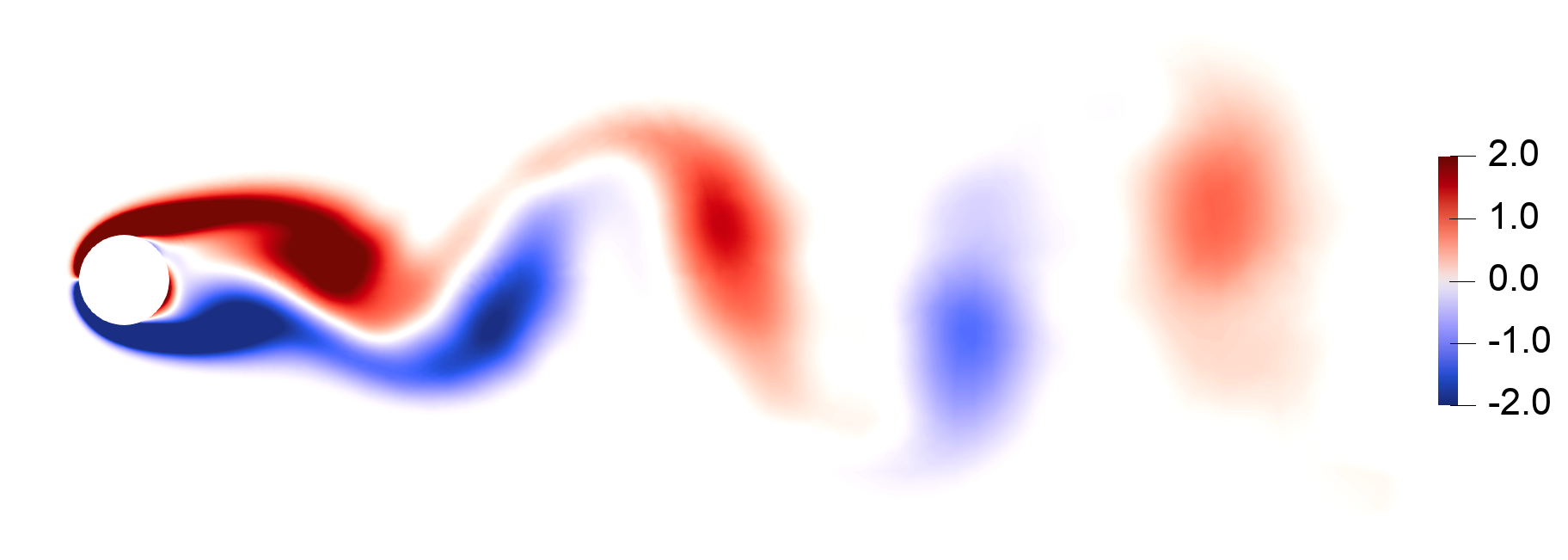}
	\caption{}
	\end{subfigure}%
	%\vglue 0.1 truein
	\caption{Two-dimensional flow problem over a right circular cylinder -- snapshots of the solution vorticity computed at $t = 50$ using: the HDM (a); LSPG-GNAT HPROM (b); LSPG-least-squares 
	collocation HPROM (c); and the LSPG-ECSW HPROM.}
	\label{fig:cylindervort}
\end{figure}

Figure \ref{fig:cylindervort} compares a snapshot of the solution vorticity computed at $t = 50$ using the HDM with counterparts computed using the various LSPG HPROMs of dimension $n=35$ described above.
Note that $t = 50$ is the time-instance just before that at which the simulation performed using the GNAT hyperreduction method terminates prematurely. Again, the reader can observe
that for the same PROM -- namely, the LSPG PROM -- the ECSW-based LSPG HPROM delivers the best reconstruction of the solution at this time-instance.

Using a single core of an Intel Xeon Gold 6128 processor running at 3.40 GHz on a machine equipped with 192 GB of memory, the HDM-based simulation takes 65.4 \si{\hour} wall-clock time to complete. 
For $n=9$ ($n=35$), the construction on the same processor of the ECSW reduced mesh is performed in 27.1 (78) \si{\minute}: 26.1 (67.8) \si{\minute} for forming the convex optimization problem defined 
by (\ref{eqn:nnls}) and (\ref{eqn:eps}), and 1.0 (10.2) \si{\minute} for solving it using the parallel NNLS algorithm \cite{chapman2017}. As for the computational performance of each constructed 
ECSW-based HPROM, it is summarized in Table \ref{tab:cylinderspeed} below. For $n=9$, the reported wall-clock time speedup factors exceed four orders of magnitude; for $n=35$, they exceed three orders 
of magnitude. 

\begin{table}[h!]
	\small
	\centering
	\caption{Two-dimensional flow problem over a right circular cylinder: computational performance of the ECSW-based LSPG HPROMs.}
	%\vglue 0.1 truein
	\begin{tabular}{ccc} 
		\toprule
		$n$  & Wall-clock time (\si{\second}) & Speedup factor   \\ \midrule
		$9$  & $23.5$ & $1.00\times10^4$ \\
		$35$ & $140.0$ & $1.68\times10^3$ \\ \bottomrule
	\end{tabular}
	%\vglue 0.1 truein
	\label{tab:cylinderspeed}
\end{table}

%\clearpage

\subsection{Turbulent flow in the wake of Ahmed body}
\label{sec:ahmed}

\subsubsection{High-dimensional model}

\begin{figure}[h!]
	\centering
	%\vglue 0.1 truein
	\includegraphics[width=0.4\textwidth]{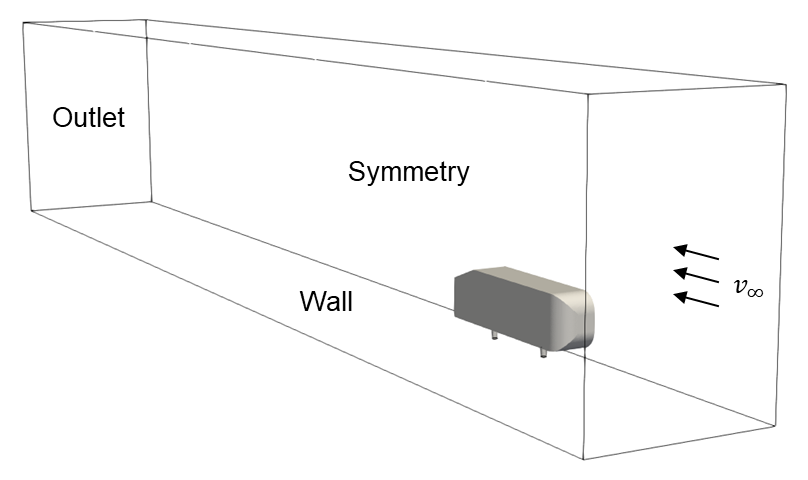}
	%\vglue 0.1 truein
	\caption{Ahmed body wake flow problem: computational domain.}
	\label{fig:ahmedmesh}
\end{figure}

The second application presented here concerns the DES-based numerical simulation of a turbulent flow past the Ahmed body geometry -- a popular benchmark CFD problem in the automotive industry 
\cite{ahmed1984}. The considered flow problem is for the slant angle of $20$\si{\degree}, the free-stream velocity $v_\infty = 60$ \si{\meter/\second}, and the Reynolds number (based on the 
body length) $Re = 4.29\times10^6$. For this problem, which has a symmetry plane, the computational domain associated with one half of the model is shown in Figure \ref{fig:ahmedmesh}. It is discretized 
by an unstructured CFD mesh with $2,890,434$ vertices and $17,017,090$ tetrahedral elements. Adiabatic boundary conditions are applied on all wall boundary surfaces of the computational domain using
Reichardt's law of the wall. As previously mentioned, the DES approach is based here on the Spalart-Allmaras one-equation model \cite{strelets2001}: therefore, the compuational model has
$d_e = 6$ DOFs per vertex, leading to a semi-discrete HDM of dimension $N = 17,342,604$.

The semi-discrete HDM is time-integrated using the implicit, second-order, three-point backward difference scheme and the fixed computational time-step $\Delta t = 8\times10^{-5}$ \si{\second}. 
The unsteady, HDM-based simulation is initialized using a quasi-steady flow around the geometry at the same flow conditions outlined above in the time-interval $[0, 2\times10^{-1}]$ \si{\second}. 
It is performed on 240 cores of a Linux cluster where each compute node is equipped with two 6-core Intel Xeon Gold 6128 processors running at 3.40 GHz and 192 GB of memory. Using this parallel 
computational platform, the HDM-based simulation takes 12.1 \si{\hour} wall-clock time to complete -- which corresponds to an aggregate CPU time of $2.91\times10^3$ \si{\hour}. 
It predicts the time-averaged drag coefficient $\bar{c}_D = 0.263$ -- with time-averaging performed in the time-interval $[5\times10^{-2},2\times10^{-1}]$ \si{\second}, after the 
startup transients have vanished. This value of $\bar{c}_D$ differs by just $3.1\%$ from its experimental value $\bar{c}_D^{\:exp} = 0.255$ \cite{ahmed1984}.

\subsubsection{Local ECSW-based hyperreduced LSPG projection-based reduced-order models}

\begin{figure}[h!]
	\centering
	%\vglue 0.1 truein
	\begin{subfigure}[c]{0.45\textwidth}%
	\hspace{0.95em} \includegraphics[width=0.85\textwidth]{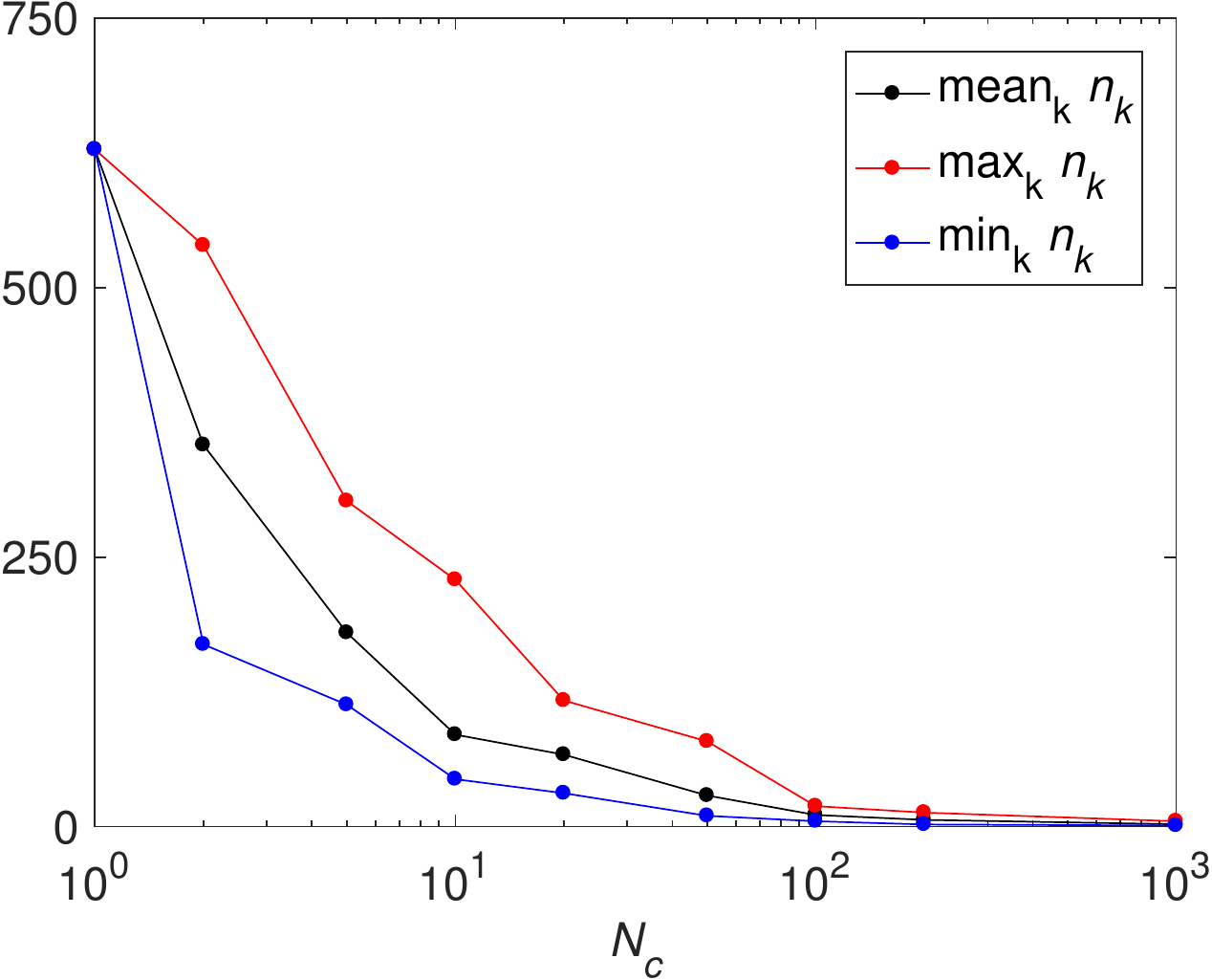}
	\caption{}
	\end{subfigure}%
	\hspace{1em}
	\begin{subfigure}[c]{0.45\textwidth}%
	\hspace{0.75em} \includegraphics[width=0.85\textwidth]{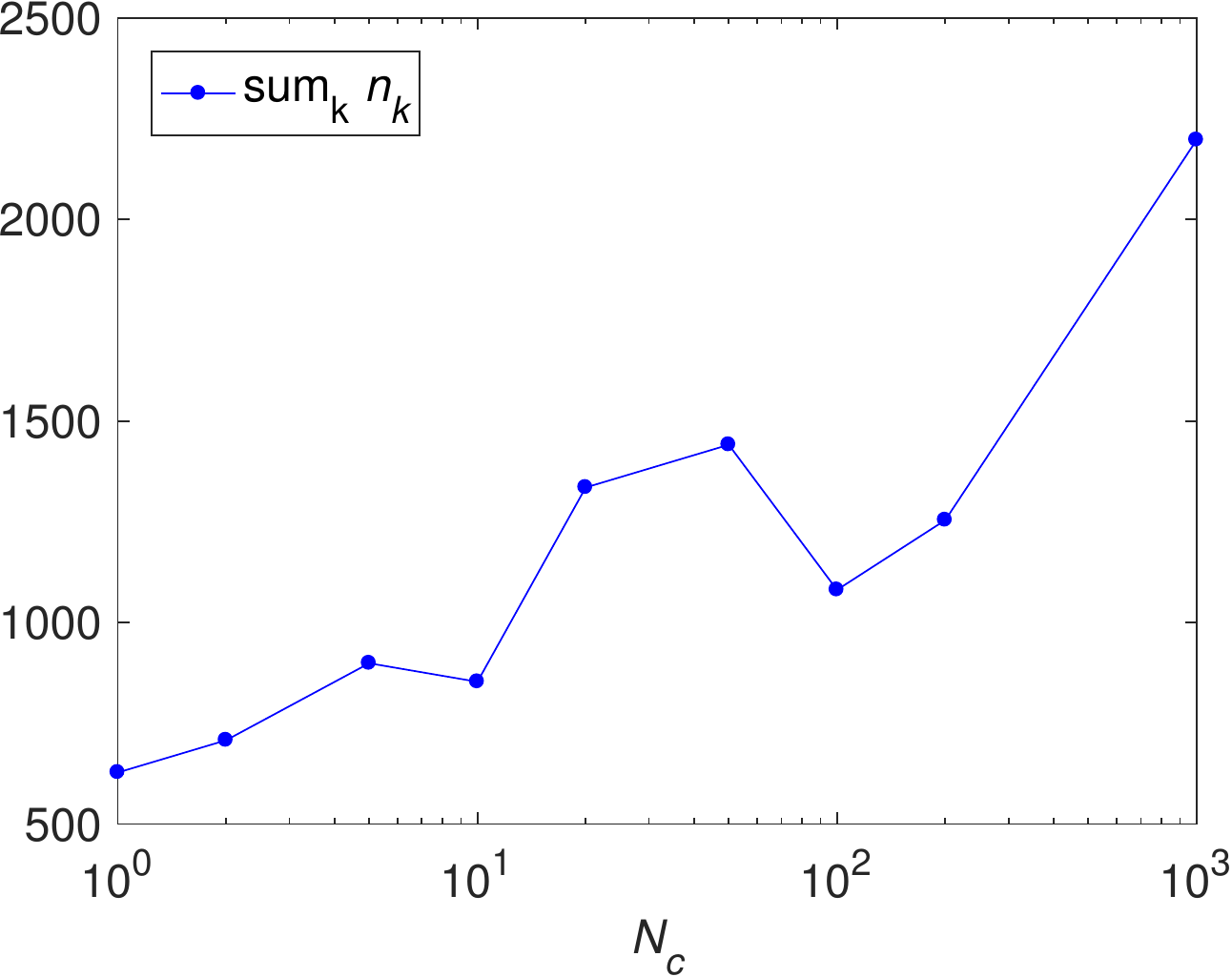}
	\caption{}
	\end{subfigure}%
	%\vglue 0.1 truein
	\caption{Ahmed body wake flow problem -- variations with the number of clusters of: the mean, maximum, and minimum ROB dimensions (a); and the total number of ROB modes across all clusters (b).}
	\label{fig:ahmedpod}
\end{figure}

For this problem, $N_s = 1,251$ HDM-based solution snapshots are collected in the time-interval $[0, 2\times10^{-1}]$ \si{\second} using the sampling rate defined by $\Delta s = 1.6\times10^{-4}$ 
\si{\second}. The $k$-means clustering algorithm is applied, along with the POD method of snapshots, to construct piecewise-affine local approximation subspaces. To this end, various values of $N_c$ in 
the range of $2$ to $100$ are considered (see below) and the overlapping factor in Algorithm \ref{alg:overlap} is set to $\phi = 0.1$. For the fixed singular value energy threshold of $99.99\%$, Figure 
\ref{fig:ahmedpod} reports the variation with $N_c$ of the mean, maximum, and minimum dimensions of the right ROBs defining the local approximation subspaces as well as the growth of the total number of 
retained POD modes $\sum\limits_{k=1}^{N_c} n_k$ with the number of clusters $N_c$. In general, choosing an appropriate value for $N_c$ amounts to trading the part of the online cost associated with the 
solution at each time-step of an implicit time-integration scheme of the dense systems of algebraic equations incurred by Newton's method, which scales as $O(n_k^3)$, with the remaining part of the 
online cost, which scales as $O(N_c^2)$. More importantly, as $N_c$ is an input to many standard clustering algorithms including the $k$-means algorithm considered here, it is desirable for the 
accuracy of the resulting local PROM to be largely insensitive to this parameter, in order to avoid as much as possible its tuning. Therefore to verify this behavior, three different values of $N_c$ are 
considered here: $N_c=10$, $N_c=50$, and $N_c=100$.

To accelerate the processing of each constructed local LSPG PROM, the proposed ECSW framework is applied for transforming such a PROM into an LSPG HPROM that employs the same time-discretization 
scheme as the HDM. Mesh sampling is performed in this case using $N_H = 17$ training residuals sampled in the same time-interval $[0, 2\times10^{-1}]$ \si{\second} with 
$\Delta s_H = 1.2\times10^{-2}$ \si{\second}. Table \ref{tab:ahmedecsw} reports for each considered value of $N_c$, the mean, maximum, and minimum dimension of each local ROB and the size 
$\widetilde{N}_e$ of the reduced mesh defining the constructed ECSW-based LSPG HPROM. Note that when $N_c > N_H$, the offline ECSW training procedure exploits data from only a subset of the $N_c$ 
clusters and therefore bets on the fact that the reduced mesh constructed using this data will perform well for the local subspace approximations whose data was not accounted for in the training process.

\begin{table}[tb]
	\small
	\centering
	\caption{Ahmed body wake flow problem -- variations with the number of clusters of: the mean, maximum, and minimum ROB dimensions; and the size of the ECSW reduced mesh.}
	%\vglue 0.1 truein
	\begin{tabular}{cccccc} 
		\toprule
		$N_c$ & $\mean_k n_k$ & $\max_k n_k$ & $\min_k n_k$ & $\widetilde{N}_e$ & $\widetilde{N}_e/N_e$ ($\%$) \\ \midrule
		$10$ & $85.3$ & $229$ & $44$ & $1,620$ & $0.056$ \\
		$50$ & $28.8$ & $79$ & $10$ & $347$ & $0.012$ \\
		$100$ & $10.8$ & $19$ & $5$ & $137$ & $0.0047$ \\ \bottomrule
	\end{tabular}
	%\vglue 0.1 truein
	\label{tab:ahmedecsw}
\end{table}

\subsubsection{Performance of the ECSW hyperreduction method}

\begin{figure}[h!]
	\centering
	%\vglue 0.1 truein
	\begin{subfigure}[c]{0.475\textwidth}%
	\centering
	\includegraphics[width=\textwidth]{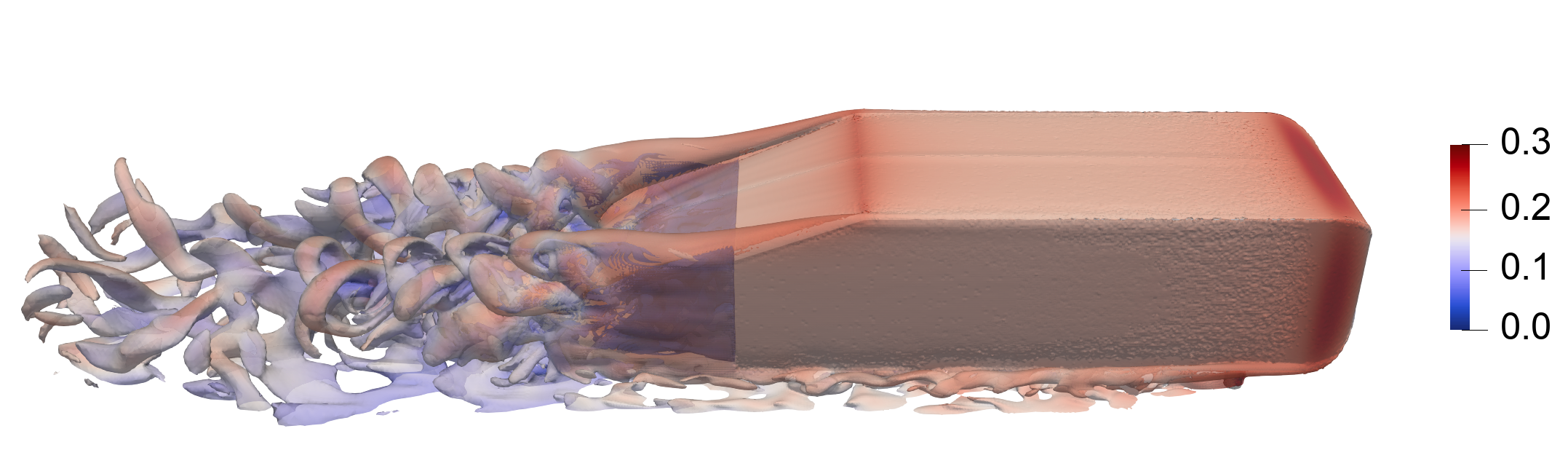}
	\caption{}
	\end{subfigure}%
	\hspace{1em}
	\begin{subfigure}[c]{0.475\textwidth}%
	\centering
	\includegraphics[width=\textwidth]{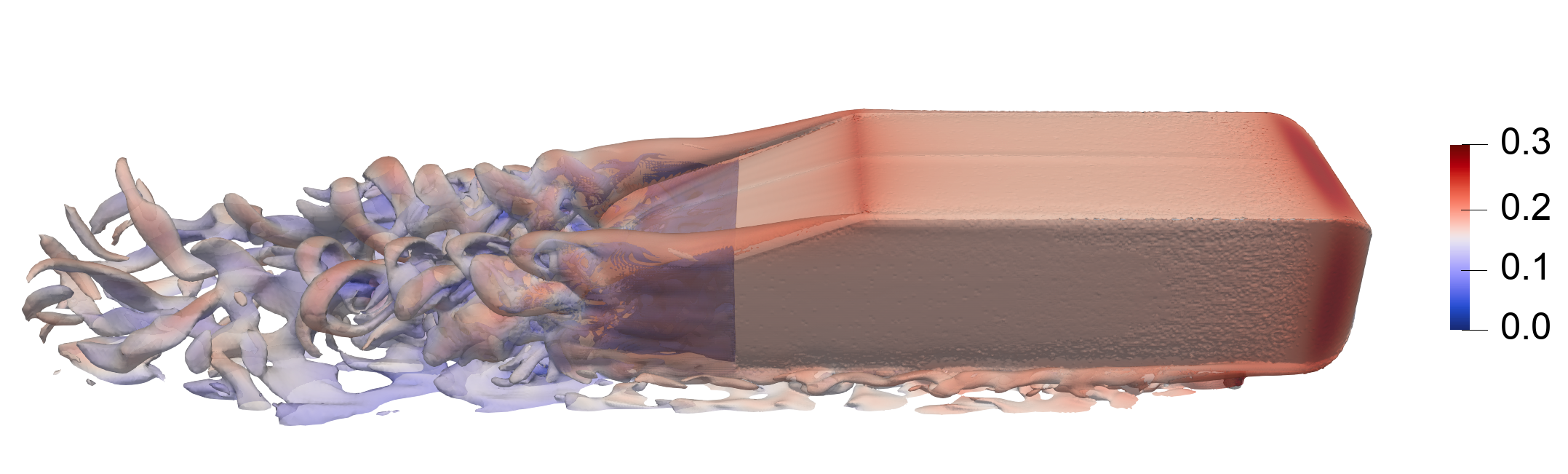}
	\caption{}
	\end{subfigure}%
	
	\medskip
	\begin{subfigure}[c]{0.475\textwidth}%
	\centering
	\includegraphics[width=\textwidth]{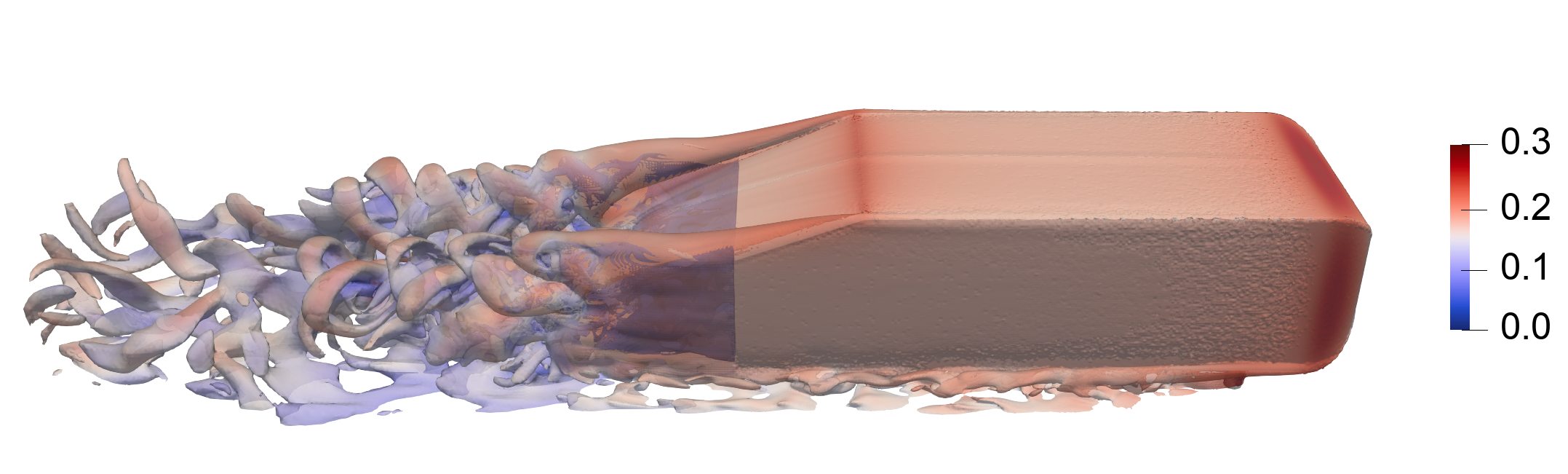}
	\caption{}
	\end{subfigure}%
	\hspace{1em}
	\begin{subfigure}[c]{0.475\textwidth}%
	\centering
	\includegraphics[width=\textwidth]{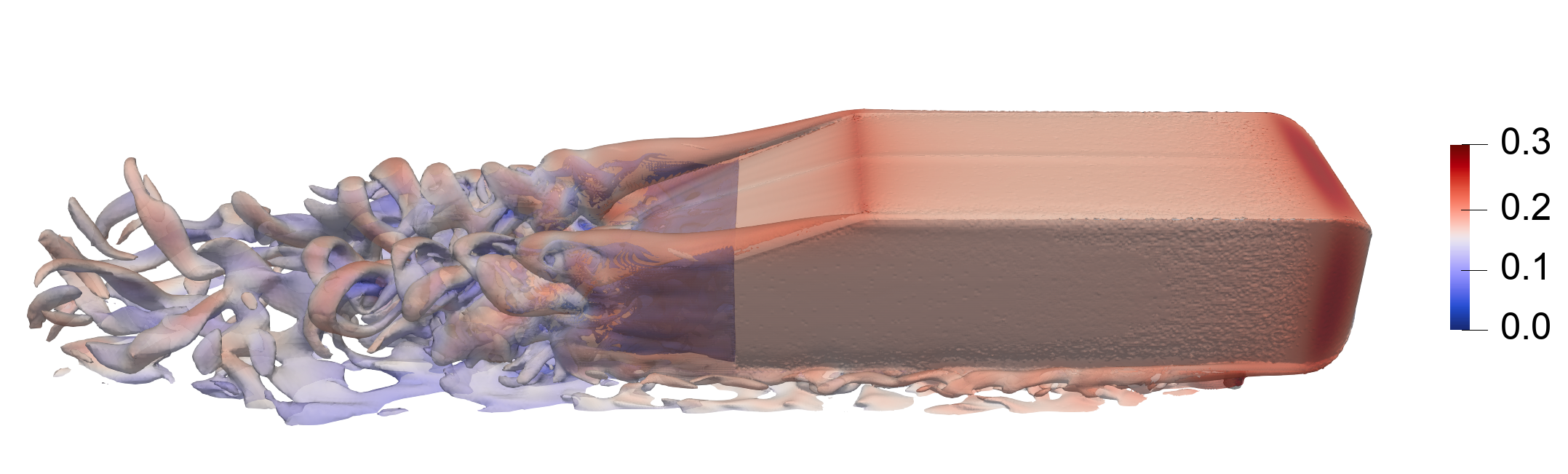}
	\caption{}
	\end{subfigure}%
	%\vglue 0.1 truein
	\caption{Ahmed body wake flow problem -- isosurfaces of the vorticity magnitude colored by the Mach number computed at $t = T_f = 2\times10^{-1}$ \si{\second} using: the HDM (a); and the ECSW-based local LSPG 
	HPROMs with $N_c=10$ (b), $N_c=50$ (c), and $N_c=100$ (d).}
	\label{fig:ahmedvort}
\end{figure}

Figure \ref{fig:ahmedvort} displays snapshots of isosurfaces of the flow solution vorticity magnitude computed at the end of the simulation time-interval $t = T_f = 2\times10^{-1}$ \si{\second} 
using the HDM and local ECSW-based LSPG HPROMs, and colored by the local Mach number. The reader can observe that the local-HPROM-based solutions do not exhibit any major discrepancy with their 
HDM-based counterparts: for each considered value of $N_c$, they capture well the highly turbulent structures present in the wake of the body.

\begin{figure}[h!]
	\centering
	%\vglue 0.1 truein
	\begin{subfigure}[c]{0.45\textwidth}%
	\hspace{0.5em} \includegraphics[width=0.85\textwidth]{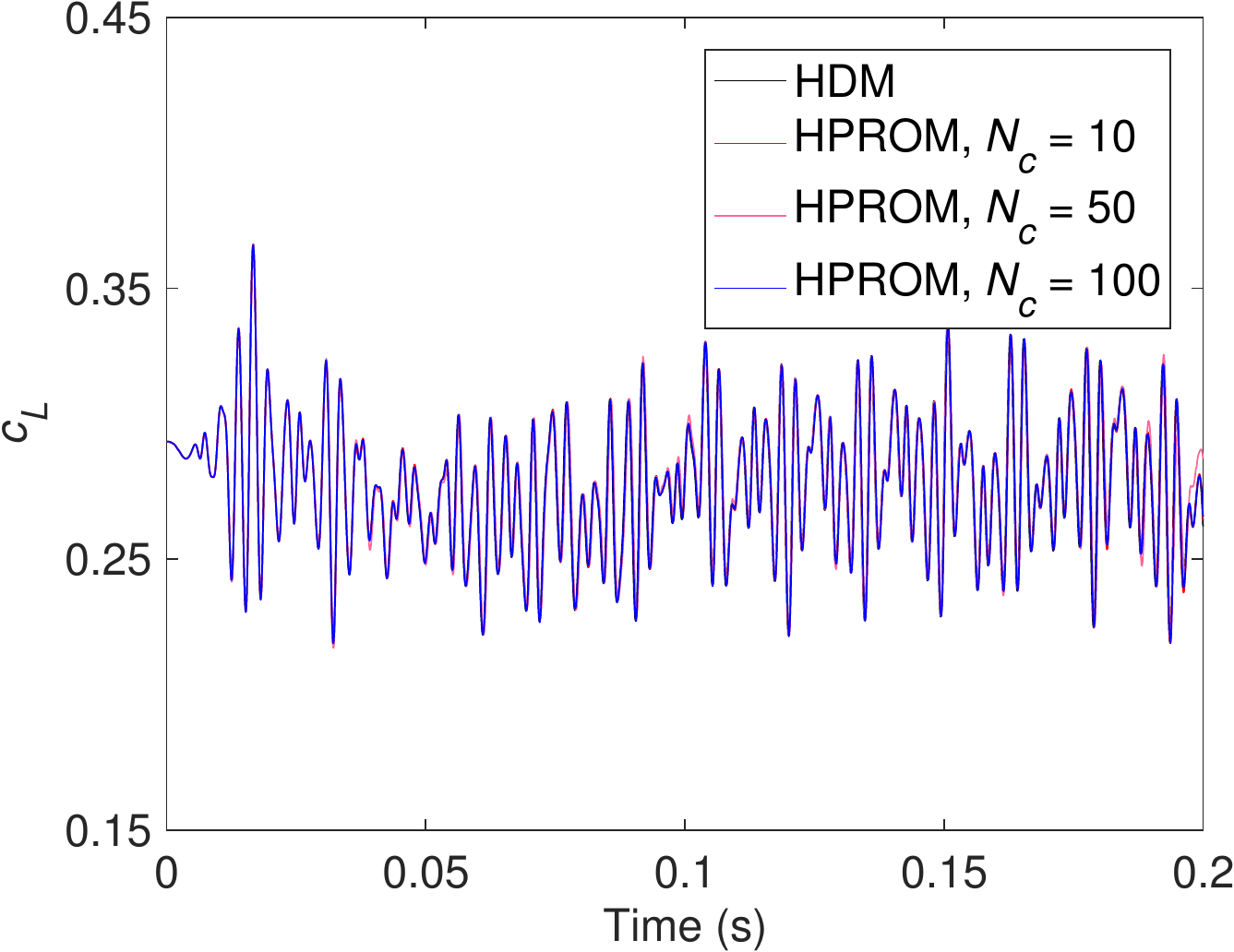}
	\caption{}
	\end{subfigure}%
	\hspace{1em}
	\begin{subfigure}[c]{0.45\textwidth}%
	\hspace{0.5em} \includegraphics[width=0.85\textwidth]{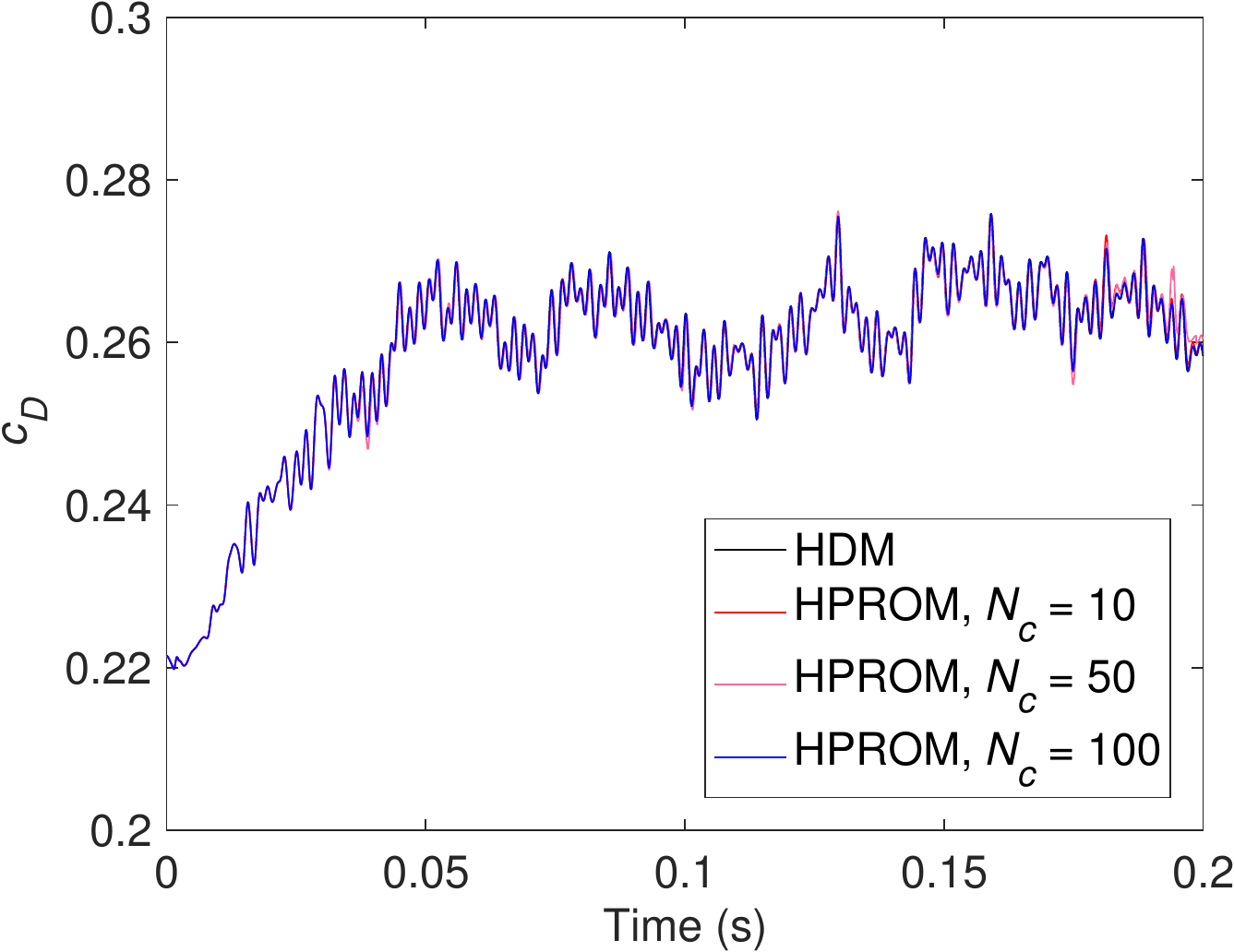}
	\caption{}
	\end{subfigure}%
	
	\medskip
	\begin{subfigure}[c]{0.45\textwidth}%
	\hspace{0.5em} \includegraphics[width=0.85\textwidth]{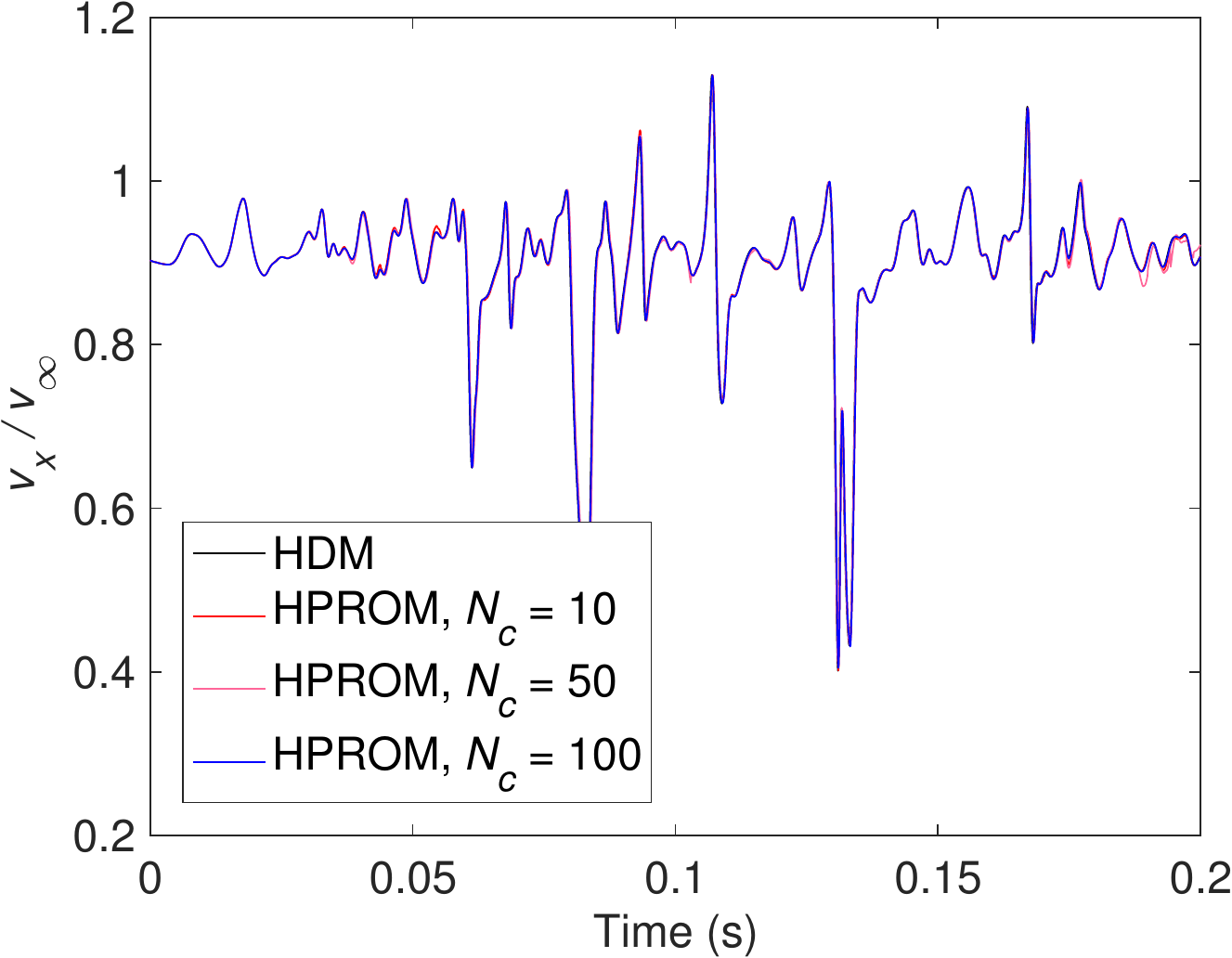}
	\caption{}
	\end{subfigure}%
	\hspace{1em}
	\begin{subfigure}[c]{0.45\textwidth}%
	\hspace{0.5em} \includegraphics[width=0.85\textwidth]{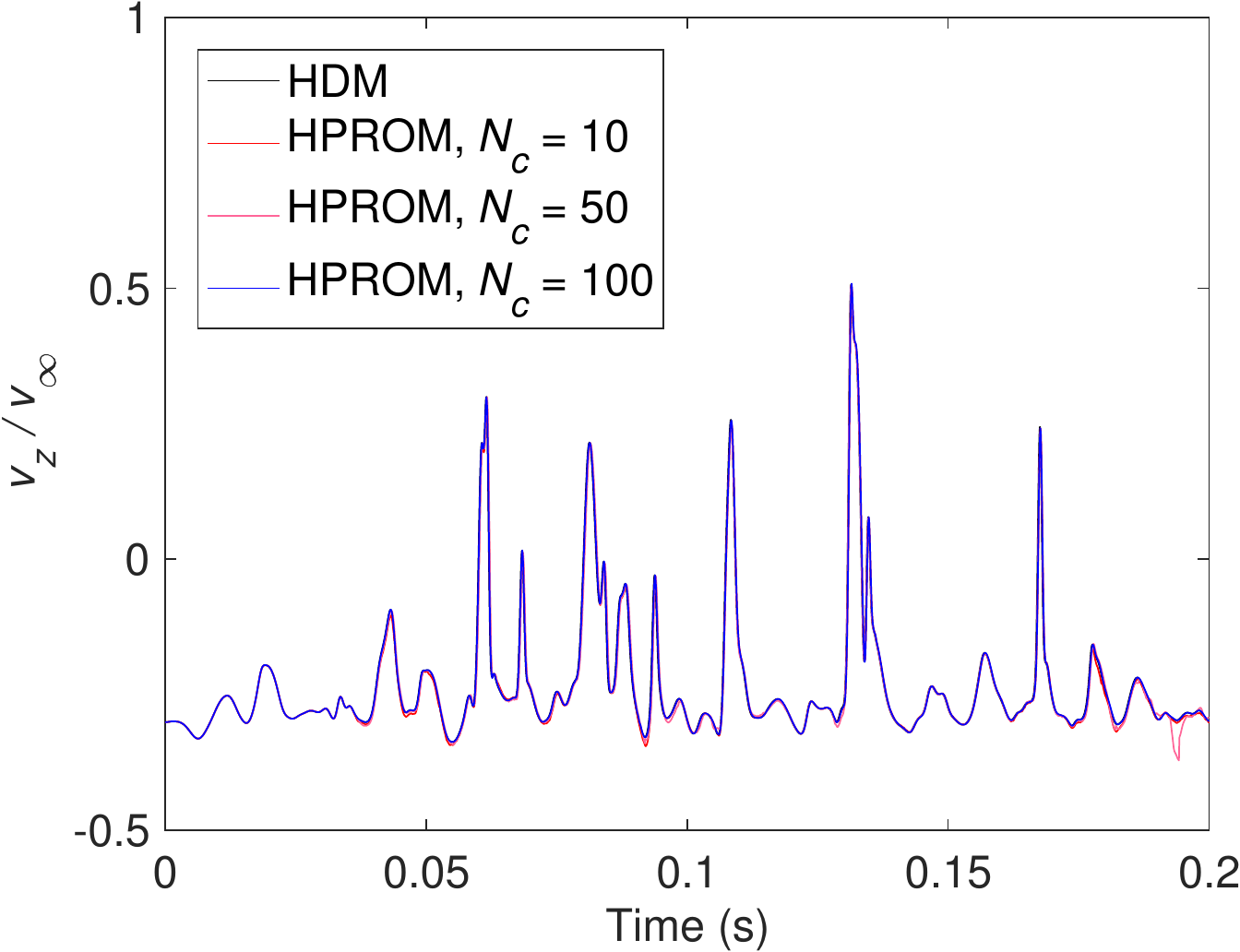}
	\caption{}
	\end{subfigure}%
	%\vglue 0.1 truein
	\caption{Ahmed body wake flow problem -- time-histories computed using the HDM and various constructed local ECSW-based LSPG HPROMs of: the lift coefficient (a); drag coefficient (b); streamwise 
	velocity component at a probe (c); and the normal velocity component at the same probe (d).}
	\label{fig:ahmeddrag}
\end{figure}

Figure \ref{fig:ahmeddrag} compares the time-histories computed using the HDM and local ECSW-based LSPG HPROMs of the lift coefficient, drag coefficient, and streamwise and vertical velocity components 
at a probe located one half body length downstream at half of the body height and along the body centerline. Again, the reader can observe that all constructed local ECSW-based HPROMs perfectly reproduce 
the time-history of each considered QoI over the entire simulation time-interval. These results suggest that the accuracy of the local ECSW-based HPROMs is largely insensitive to the choice of $N_c > 1$.
For each QoI predicted using a local ECSW-based HPROM, the relative error computed using $\Delta s_{\mathbb{RE}}=\Delta t=8\times10^{-5}$ \si{\second} is reported in Table \ref{tab:ahmederr}. The 
reported errors demonstrate a high level of achieved accuracy and an insensitivity with respect to the choice of $N_c$. For the specified singular value energy threshold, a global ROB ($N_c = 1$) would 
include over 600 POD modes (see Figure \ref{fig:ahmedpod}); on the other hand, local subspace approximations enable the same level of accuracy using a database of ROBs of much lower dimensions.

\begin{table}[h!]
	\small
	\centering
	\caption{Ahmed body wake flow problem -- computational accuracy of the various constructed LSPG HPROMs.}
	%\vglue 0.1 truein
	\begin{tabular}{ccccc}
		\toprule
		$N_c$ & $\mathbb{RE}_{c_D}$ ($\%$) & $\mathbb{RE}_{c_L}$ ($\%$) & $\mathbb{RE}_{v_x}$ ($\%$) & $\mathbb{RE}_{v_y}$ ($\%$) \\ \midrule
		$10$ & $0.262$ & $1.58$ & $1.06$ & $5.33$ \\
		$50$ & $0.419$ & $2.18$ & $1.17$ & $6.12$ \\
		$100$ & $0.212$ & $1.35$ & $0.839$ & $4.07$ \\ \bottomrule
	\end{tabular}
	%\vglue 0.1 truein
	\label{tab:ahmederr}
\end{table}

For this application, all ECSW reduced meshes are constructed on 240 cores of the aforementioned Linux cluster:
\begin{itemize}
	\item For $N_c = 10$, the mesh sampling is performed in 34.8 \si{\minute}: 11.7 \si{\minute} for forming the convex optimization problem defined by (\ref{eqn:nnls}) and (\ref{eqn:eps}), and 23.1 
		\si{\minute} for solving 
		it using the parallel NNLS algorithm \cite{chapman2017}. 
	\item For $N_c = 50$, the mesh sampling is performed in 2.92 \si{\minute}: 0.95 \si{\minute} for forming the aforementioned convex optimization problem and 1.97 \si{\minute} for solving it using the parallel NNLS algorithm.
	\item For $N_c = 100$, the mesh sampling is performed in 2.76 \si{\minute}: 2.06 \si{\minute} for forming the convex optimization problem and 0.70 \si{\minute} for solving it using the same parallel algorithm.
\end{itemize}
These performance results demonstrate that the proposed hyperreduction procedure is computationally tractable even for a large-scale HDM.

\begin{figure}[h!]
	\centering
	%\vglue 0.1 truein
	\includegraphics[width=0.3825\textwidth]{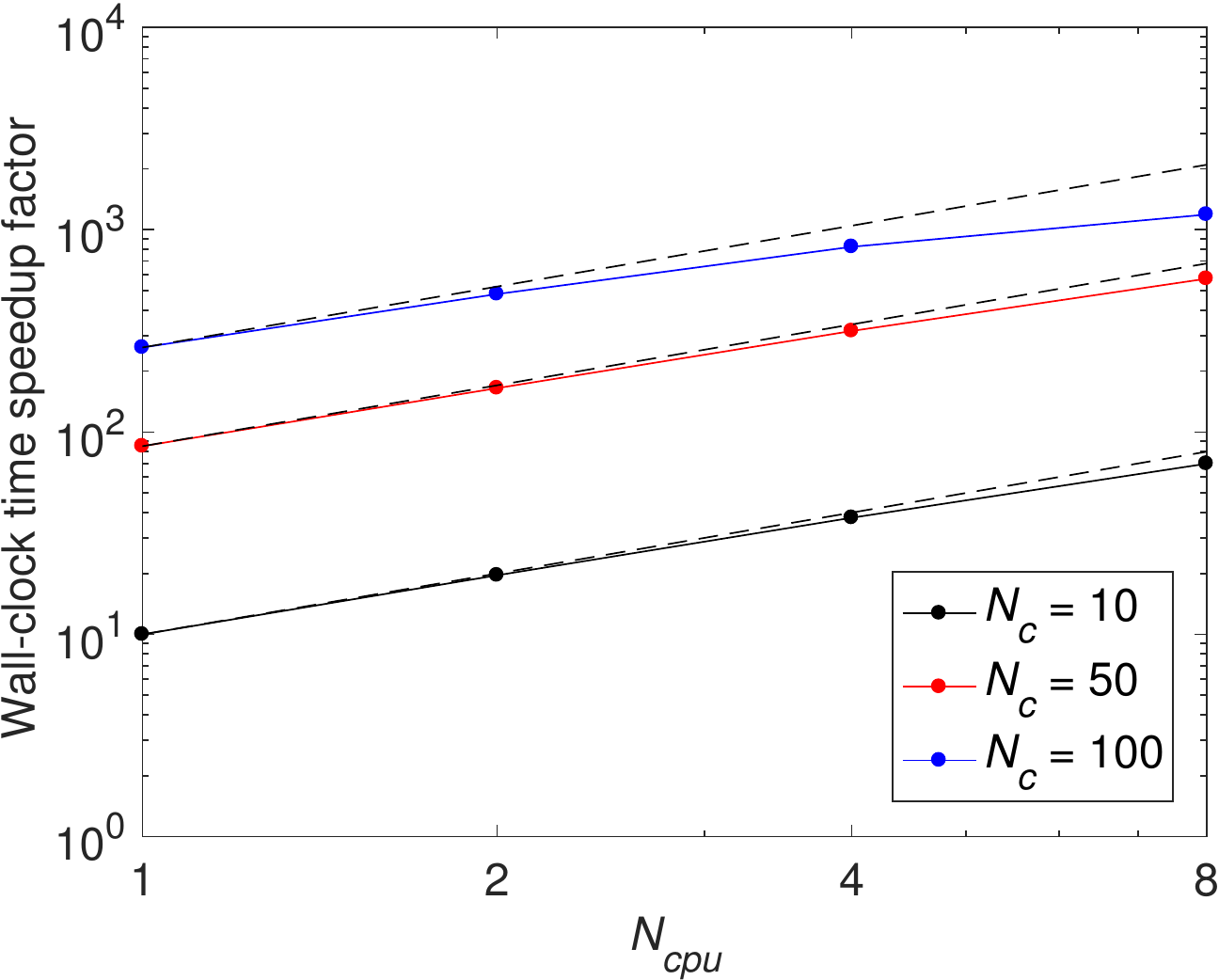}
	%\vglue 0.1 truein
	\caption{Ahmed body wake flow problem -- parallel scalability of the constructed ECSW-based LSPG HPROMs and delivered speedup factors with respect to the HDM 
	(black dashed lines represent the ideal speedup factors).}
	\label{fig:ahmedcpu}
\end{figure}

All HPROM-based simulations are performed on $N_{cpu} = 1$, $2$, $4$, and $8$ cores of a single node of the aforementioned Linux cluster. Figure \ref{fig:ahmedcpu}, which reports
for each considered value of $N_c$ the obtained wall-clock time speedup factor, demonstrates the scalability of the constructed ECSW-based LSPG HPROMs on up to 8 cores -- despite the small sizes of the
considered HPROMs (see Table \ref{tab:ahmedecsw}). Table \ref{tab:ahmedcpu} focuses on the case $N_{cpu} = 8$. It reports for this case the wall-clock execution time and both wall-clock time and 
CPU time speedup factors delivered by all constructed ECSW-based LSPG HPROMs. The results show that all ECSW-based LSPG HPROMs deliver excellent speedup factors -- up to three orders of magnitude for the 
wall-clock time and exceeding four orders of magnitude for the CPU time.  

\begin{table}[h!]
	\small
	\centering
	\caption{Ahmed body wake flow problem -- computational performance of the ECSW-based LSPG HPROMs on $N_{cpu}=8$ cores: wall-clock execution time; and wall-clock time and CPU time speedup factors.}
	%\vglue 0.1 truein
	\begin{tabular}{cccc} 
		\toprule
		$N_c$ & Wall-clock time (\si{\second}) & Wall-clock time speedup factor & CPU time speedup factor \\ \midrule
		$10$ & $626.6$ & $6.97\times10^1$ & $2.09\times10^3$ \\
		$50$ & $76.3$ & $5.72\times10^2$ & $1.72\times10^4$ \\
		$100$ & $36.8$ & $1.19\times10^3$ & $3.56\times10^4$ \\ \bottomrule
	\end{tabular}
	%\vglue 0.1 truein
	\label{tab:ahmedcpu}
\end{table}

%\clearpage

\subsection{Turbulent flow past an F-16C/D Block 40 aircraft configuration at a high angle of attack}
\label{sec:f16}

\subsubsection{High-dimensional model}

\begin{figure}[h!]
	\centering
	%\vglue 0.1 truein
	\begin{subfigure}[b]{0.45\textwidth}%
	\centering
	\includegraphics[width=0.9\textwidth]{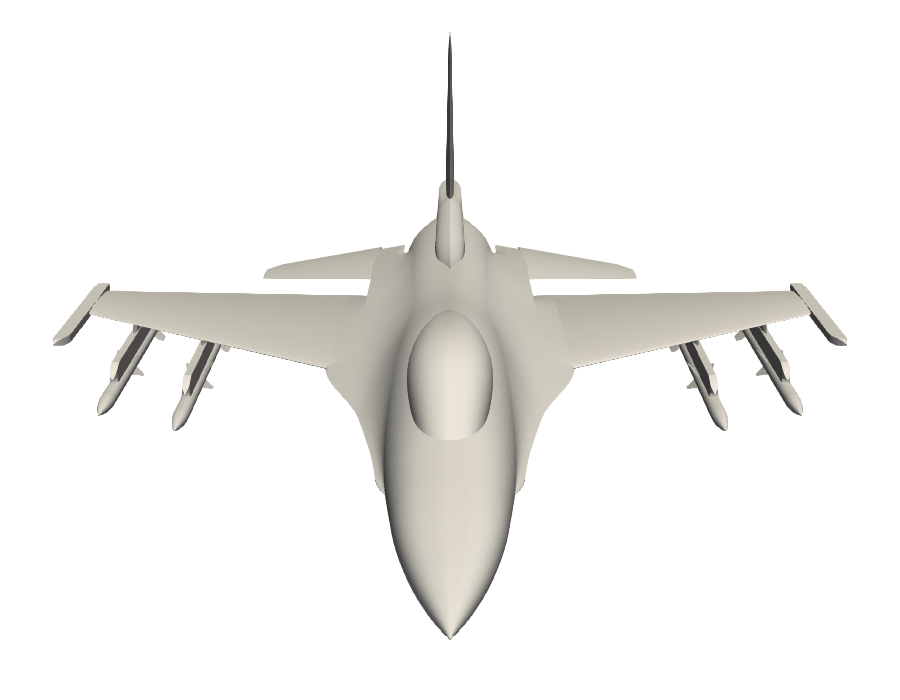}
	\caption{}
	\label{fig:f16mesha}
	\end{subfigure}%
	\hspace{1em}
	\begin{subfigure}[b]{0.45\textwidth}%
	\centering
	\includegraphics[width=\textwidth]{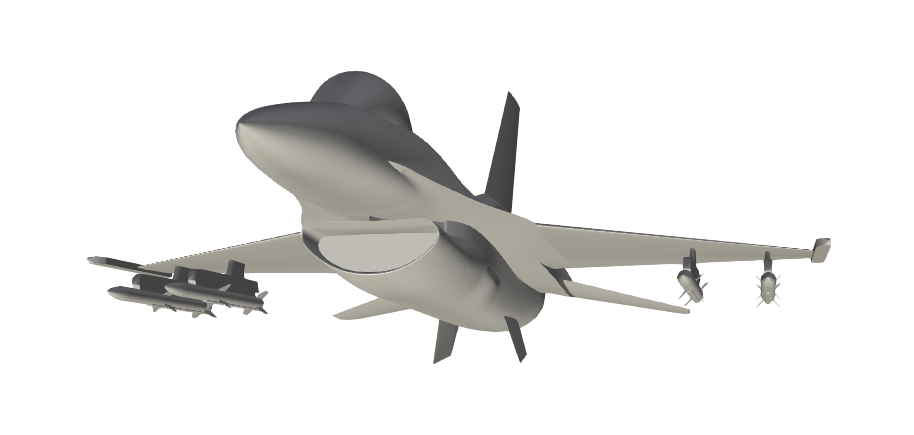}
	\vspace{0.5em}
	\caption{}
	\label{fig:f16meshb}
	\end{subfigure}%
	%\vglue 0.1 truein
	\caption{Turbulent flow past an F-16C/D aircraft -- model geometry: top view (a); and underside view (b).}
	\label{fig:f16mesh}
\end{figure}

The final application considered here focuses on the prediction of the turbulent flow past a configuration of the F-16C/D Block 40 Fighting Falcon with external stores (see Figure \ref{fig:f16mesh})
at the free-stream Mach number $M_\infty = 0.3$, $30$\si{\degree} angle of attack, and $10,000$ \si{\foot} altitude. For this geometry and flight conditions, the resulting Reynolds number based on the 
mean aerodynamic chord (MAC) of $11.32$ \si{\foot} is $Re = 1.82\times10^7$. Specifically, this aircraft configuration includes wing-tip launch rails and four occupied under-wing stations, which adds
to the complexity of the geometry, resulting unsteady flow field, and resulting wake dynamics.

As in the previous example, the flow is modeled here using a DES approach and only half of the geometry is considered for the purpose of computational efficiency. Hence, a symmetry plane along the
fuselage midsection is incorporated in the computational model. Reichardt's law of the wall is again employed for enforcing the adiabatic wall boundary conditions on the 
aircraft surface to further limit the dimensionality of the HDM. The half-aircraft model surface is discretized using $575,951$ vertices and $1,148,092$ triangular elements. The computational domain 
around the aircraft is discretized by a CFD mesh with $26,919,879$ vertices and $158,954,429$ tetrahedral elements, leading to an HDM of dimension $N = 161,519,274$. 

Time-discretization of the HDM is performed in the time-interval $[0, 5\times10^3]$ using a second-order DIRK scheme and the fixed, nondimensional time-step $\Delta t=5\times10^{-2}$
(note that $T_f = 5\times10^3$ corresponds in this case to the dimensional value of $1.29$ \si{\second}). The initial condition is set to a computed quasi-steady flow around the aircraft geometry.

Due to its very large scale, the HDM-based simulation is performed on 3,584 cores of the Frontera supercomputer at The University of Texas at Austin. Each compute node of this supercomputed is 
equipped with two 26-core Intel Xeon Platinum 8280 processors with a 2.70 GHz nominal clock speed and 192 GB of memory. The interconnect is based on Mellanox HDR InfiniBand. Using this computing
platform, the HDM-based simulation completes in 100.3 \si{\hour} wall-clock time -- or equivalently, $3.60\times10^5$ \si{\hour} CPU time.  

In summary, this application is representative of industrial-scale CFD applications. The underlying CFD problem is a showcase problem for demonstrating the feasibility of the offline part of the proposed
hyperreduction method as well as the potential of its online part for drastically accelerating DES simulations using LSPG HPROMs.

\subsubsection{Local ECSW-based hyperreduced LSPG projection-based reduced-order models}

First, an LSPG PROM is constructed for this problem. For this purpose, HDM-based solution snapshots are collected in the nondimensional time-interval $[0, 5\times10^3]$ using the nondimensional 
sampling rate defined by $\Delta s = 1$, resulting in the storage of a total of $N_s = 5,000$ snapshots. Using the same procedure based on POD and the $k$-means clustering algorithm as in the previous 
application (see Section \ref{sec:ahmed}), a piecewise-affine local subspace approximation with $N_c = 50$ is constructed. For the fixed singular value energy threshold of $99.99\%$ and the overlapping 
factor $\phi = 0.1$ in Algorithm \ref{alg:overlap}, the resulting mean, maximum, and minimum right ROB dimensions are reported in Table \ref{tab:f16pod}.  

Next, hyperreduction is performed using the proposed ECSW method with $N_H = 21$ training snapshots collected for $t \in [0, 5\times10^3]$ using $\Delta s_H = 250$. The characteristics of the generated
reduced mesh are summarized in Table \ref{tab:f16pod}. The resulting ECSW-based LSPG HPROM is discretized using the same second-order DIRK time-integration scheme used for discretizing the HDM. For this 
problem, the computational time-step of the HDM is not limited by the solution accuracy but by its numerical stability. Since the ECSW-based LSPG HPROM is less stiff than its underlying HDM, its numerical
stability can sustain in this case a larger computational time-step. Consequently, the HPROM-based simulation is performed using the nondimensional time-step $\Delta t=5\times10^{-1}$ -- that is,
10 times larger than its counterpart employed in the HDM-based simulation.

\begin{table}[h!]
	\small
	\centering
	\caption{Turbulent flow past an F-16C/D aircraft: mean, maximum, and minimum ROB dimensions of the constructed local approximation subspaces; and size of the ECSW reduced mesh.}
	%\vglue 0.1 truein
	\begin{tabular}{cccccc} 
		\toprule
		$N_c$ & $\mean_k n_k$ & $\max_k n_k$ & $\min_k n_k$ & $\widetilde{N}_e$ & $\widetilde{N}_e/N_e$ ($\%$) \\ \midrule
		$50$ & $53.7$ & $115$ & $26$ & $787$ & $0.0029$ \\ \bottomrule
	\end{tabular}
	%\vglue 0.1 truein
	\label{tab:f16pod}
\end{table}

\subsubsection{Performance of the ECSW hyperreduction method}

\begin{figure}[h!]
	\centering
	%\vglue 0.1 truein
	\begin{subfigure}[c]{0.475\textwidth}%
	\centering
	\includegraphics[width=\textwidth]{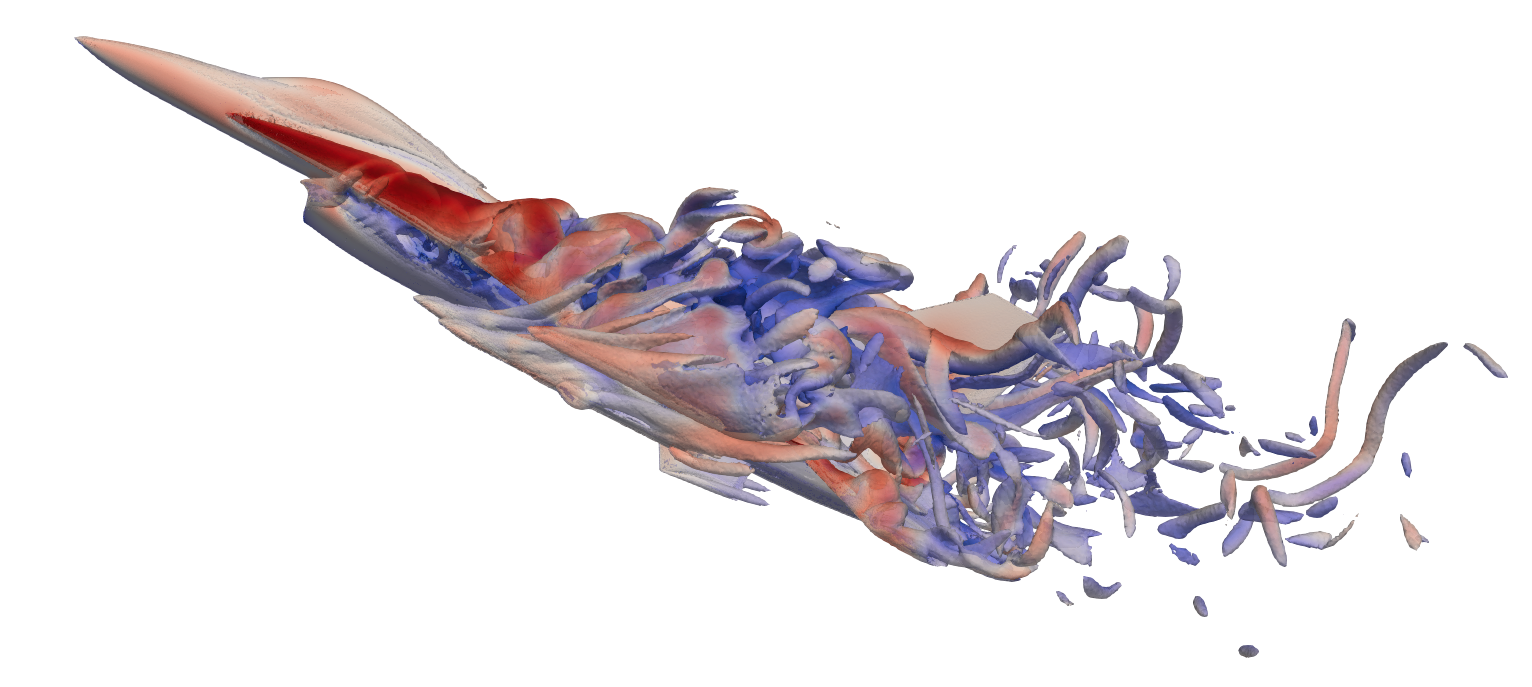}
	%\caption{}
	\end{subfigure}%
	\hspace{1em}
	\begin{subfigure}[c]{0.475\textwidth}%
	\centering
	\includegraphics[width=\textwidth]{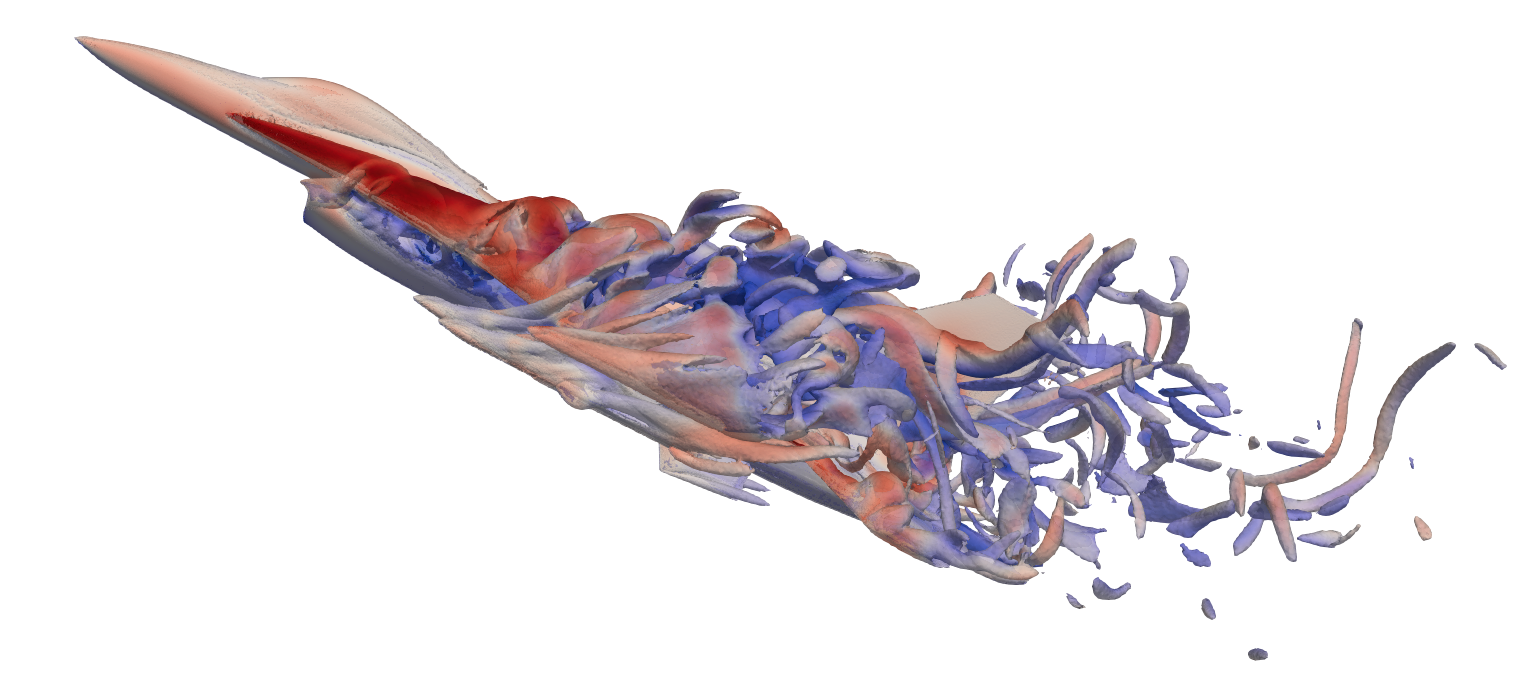}
	%\caption{}
	\end{subfigure}%
	\vspace{-0.1 truein}
	
	\begin{subfigure}[c]{0.475\textwidth}%
	\centering
	\includegraphics[width=\textwidth]{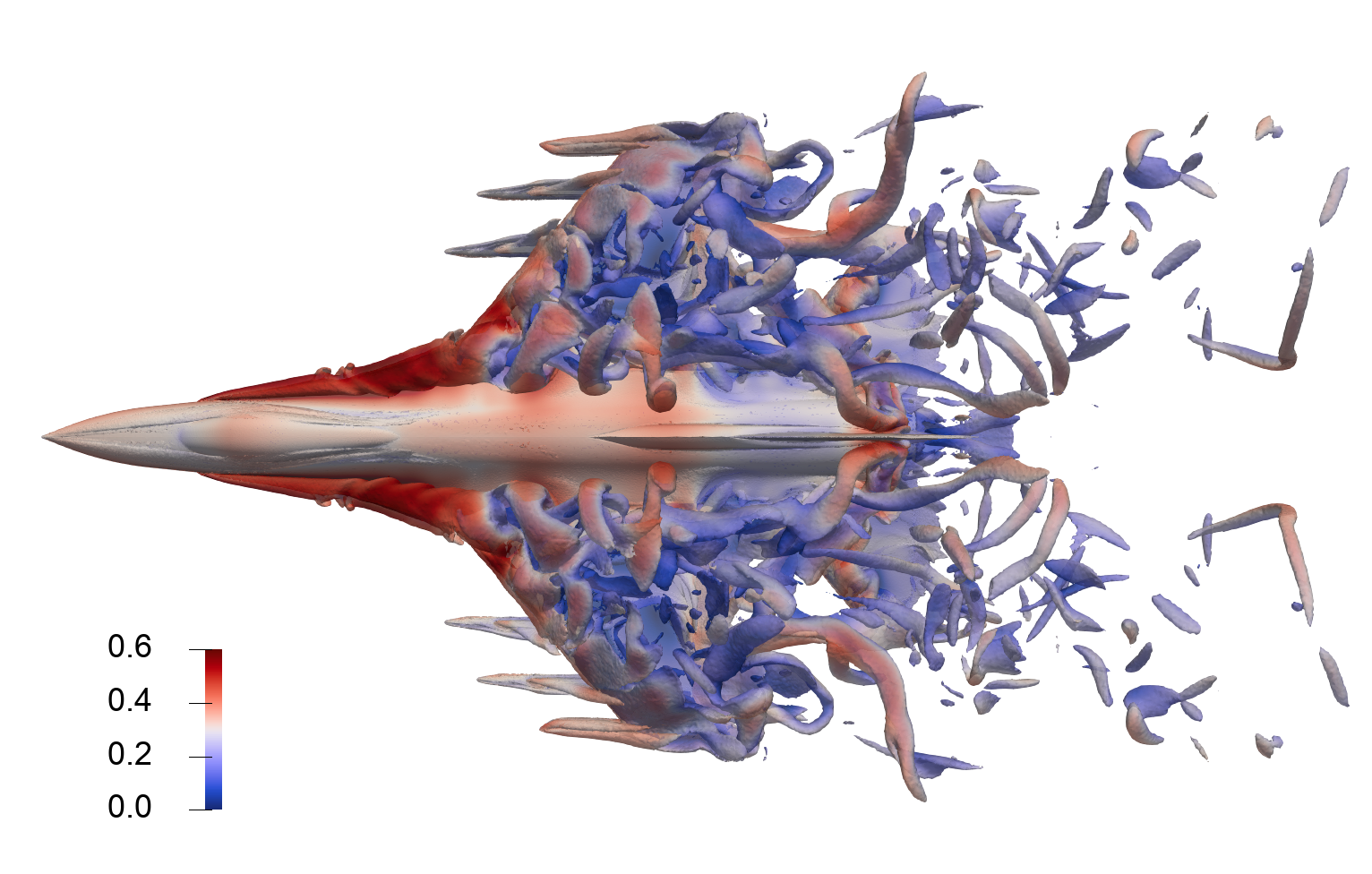}
	\caption{}
	\end{subfigure}%
	\hspace{1em}
	\begin{subfigure}[c]{0.475\textwidth}%
	\centering
	\includegraphics[width=\textwidth]{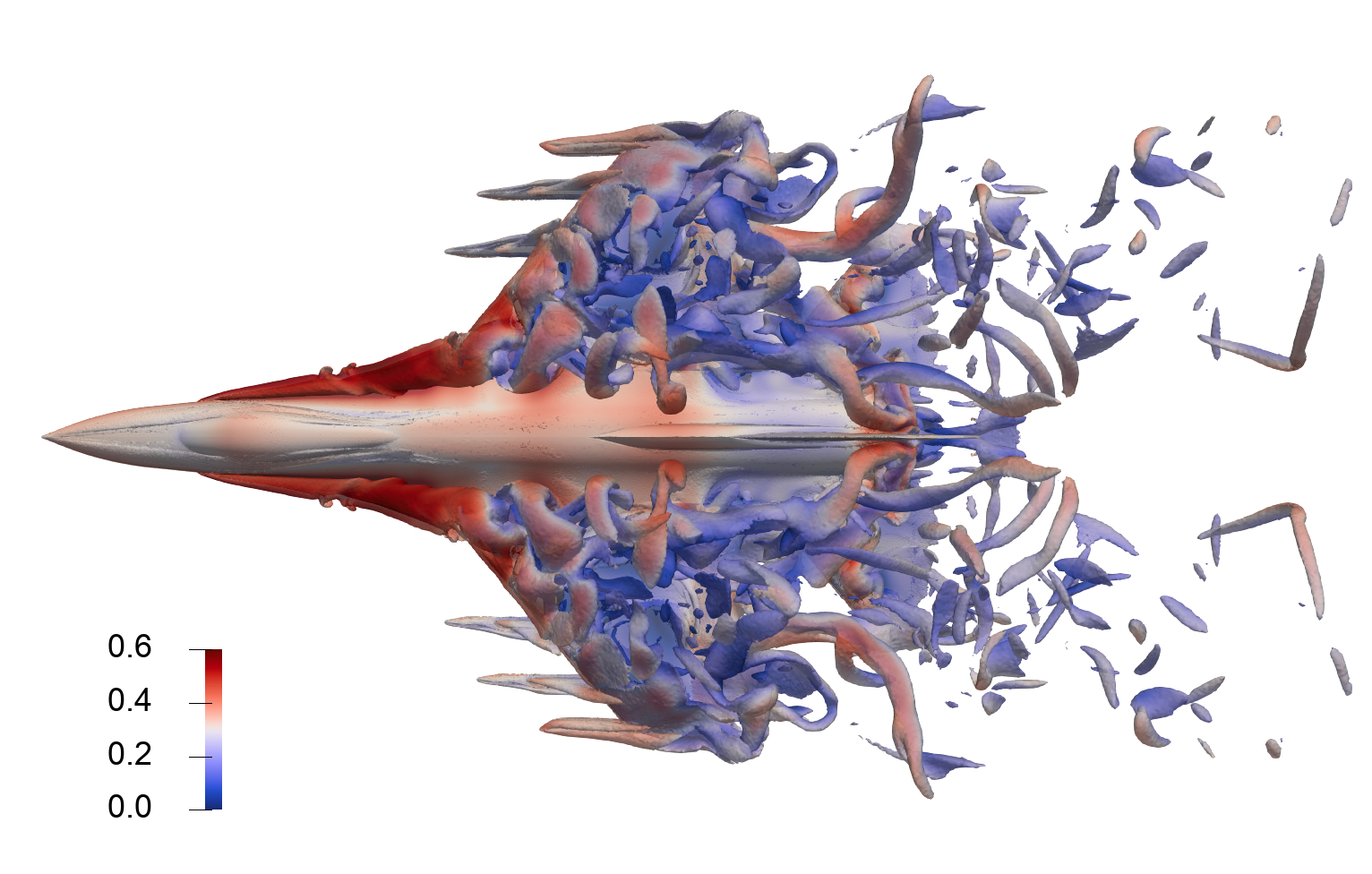}
	\caption{}
	\end{subfigure}%
	%\vglue 0.1 truein
	\caption{Turbulent flow past an F-16C/D aircraft -- isosurfaces of the solution vorticity magnitude colored by Mach number computed at $t = T_f = 1.29$ \si{\second} using: the HDM (a); and the
	ECSW-based LSPG HPROM (b).}
	\label{fig:f16vort}
\end{figure}

\begin{figure}[h!]
	\centering
	%\vglue 0.1 truein
	\begin{subfigure}[c]{0.475\textwidth}%
	\centering
	\includegraphics[width=\textwidth]{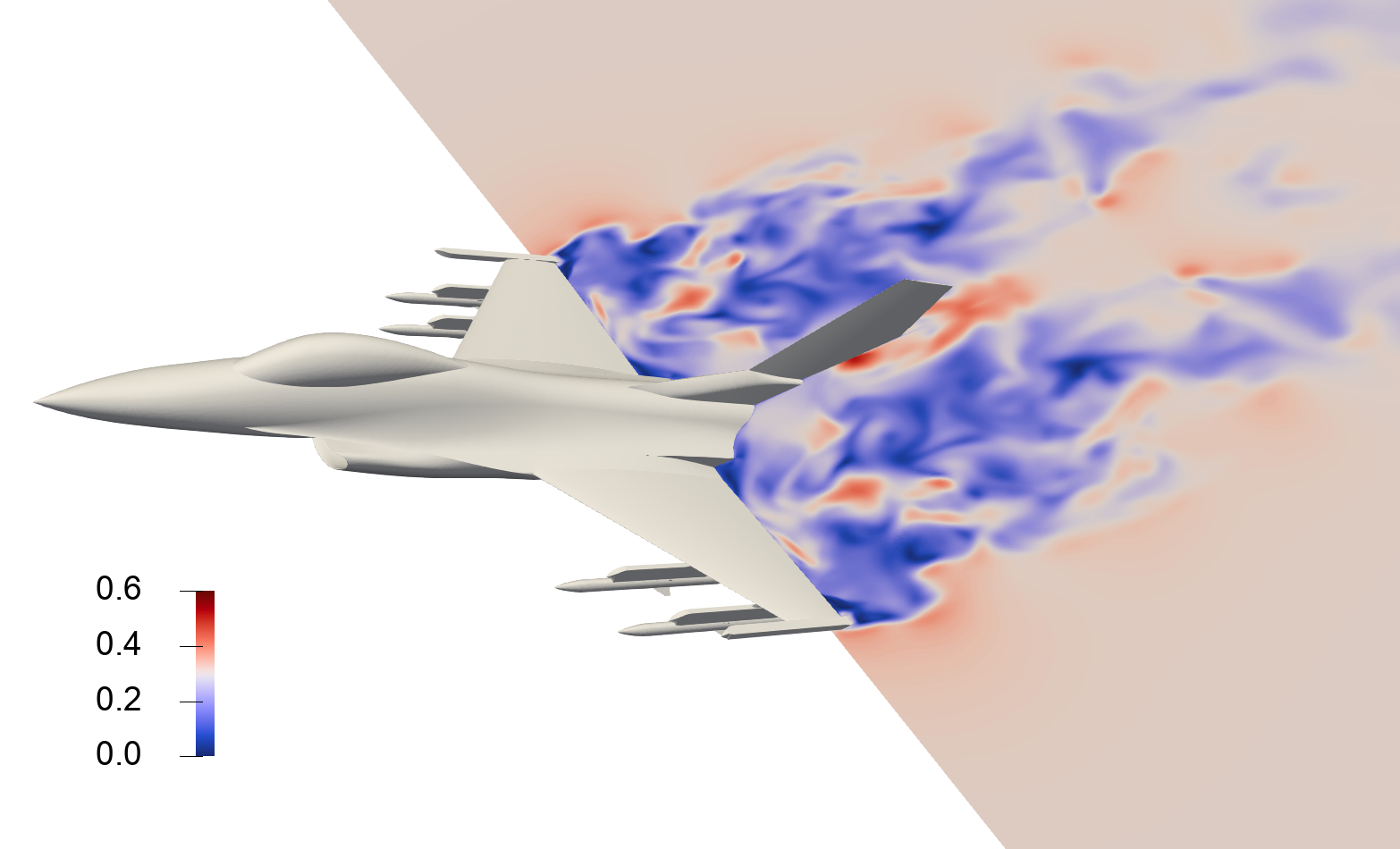}
	\caption{}
	\end{subfigure}%
	\hspace{1em}
	\begin{subfigure}[c]{0.475\textwidth}%
	\centering
	\includegraphics[width=\textwidth]{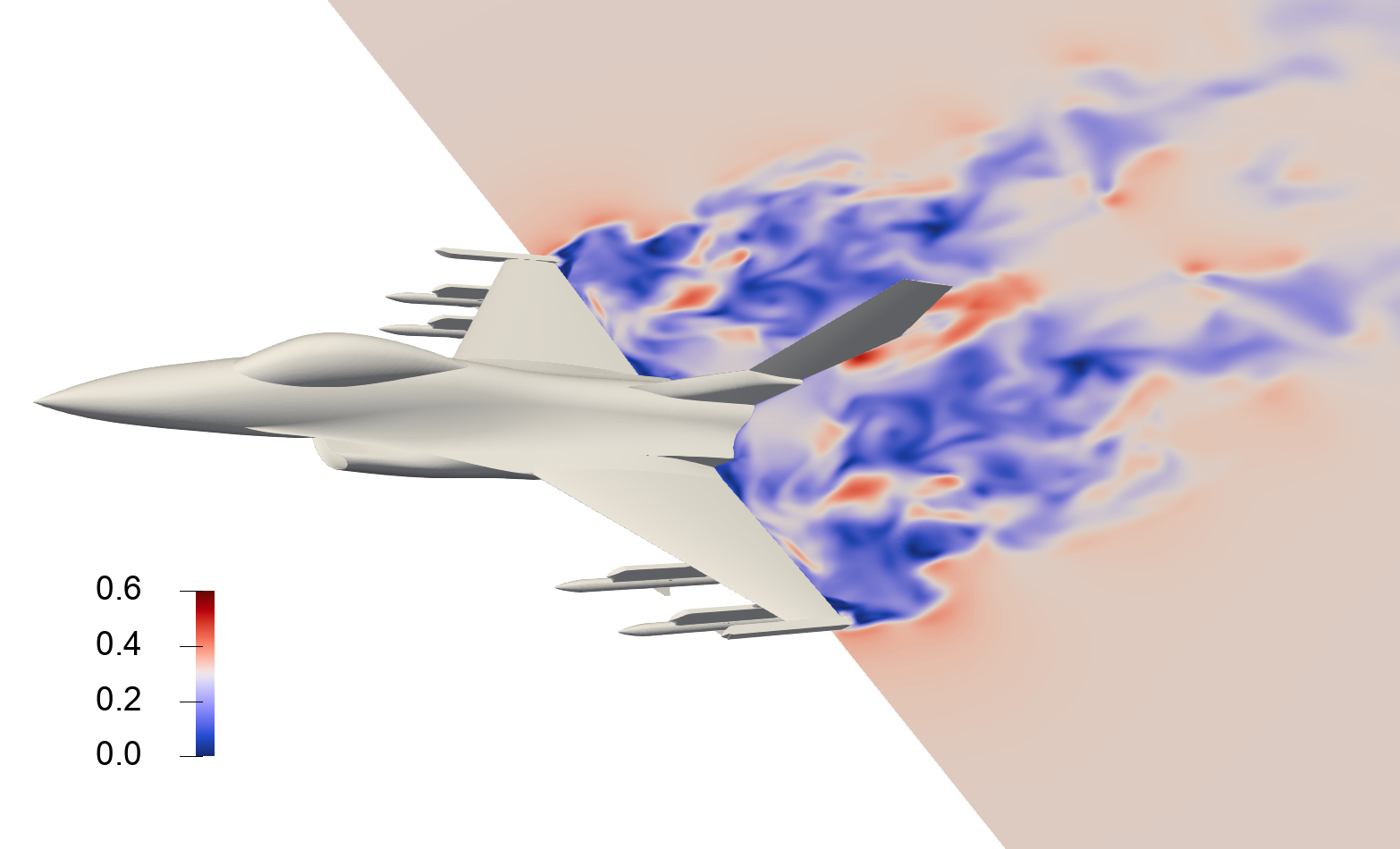}
	\caption{}
	\end{subfigure}%
	%\vglue 0.1 truein
	\caption{Turbulent flow past an F-16C/D aircraft -- snapshots of the solution Mach number in the aircraft wake computed at $t = T_f = 1.29$ \si{\second} using: the HDM (a); and the ECSW-based LSPG 
	HPROM (b).}
	\label{fig:f16mach}
\end{figure}

\begin{figure}[h!]
	\centering
	%\vglue 0.1 truein
	\begin{subfigure}[c]{0.45\textwidth}%
	\hspace{0.25em} \includegraphics[width=0.85\textwidth]{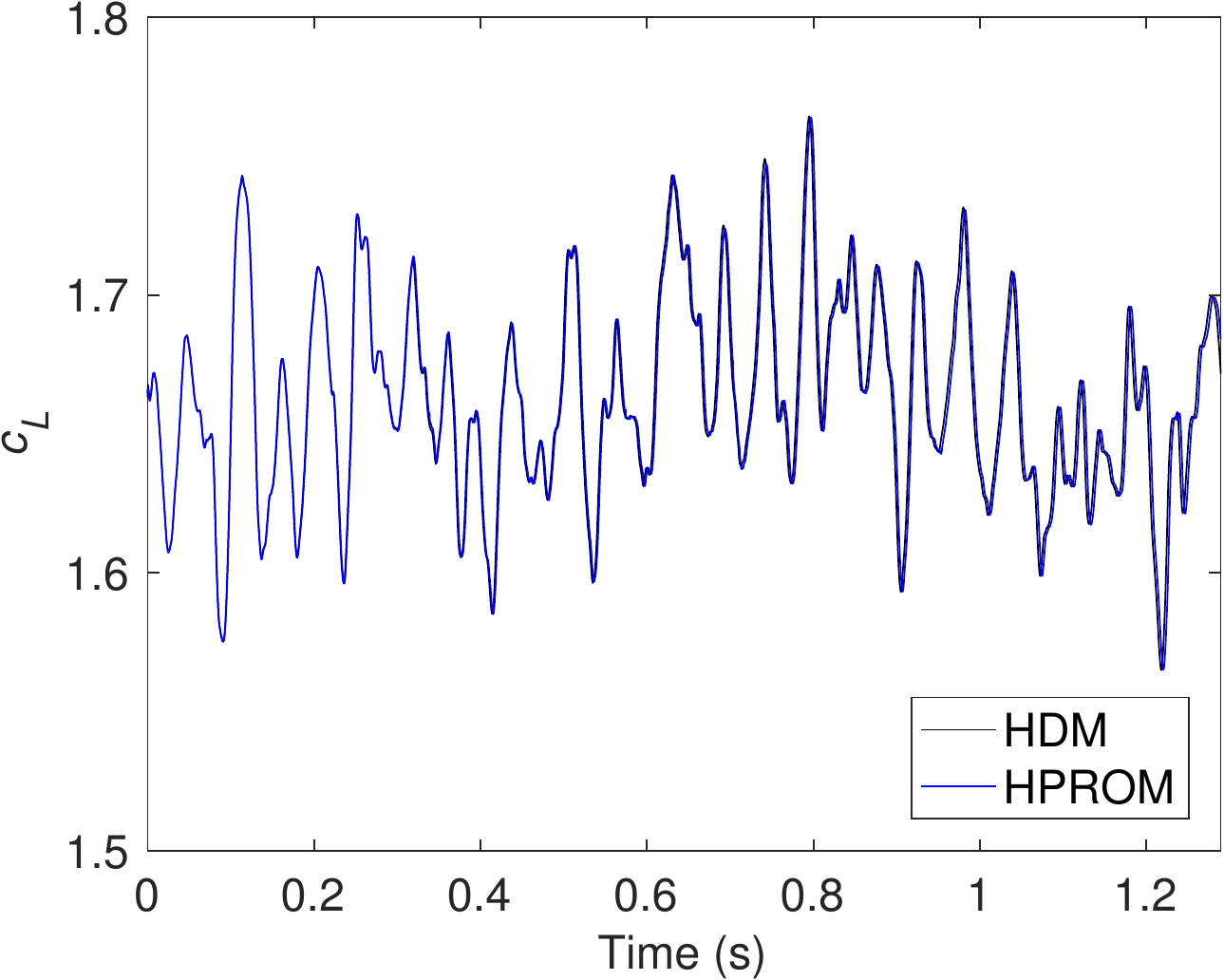}
	\caption{}
	\end{subfigure}%
	\hspace{1em}
	\begin{subfigure}[c]{0.45\textwidth}%
	\hspace{0.25em} \includegraphics[width=0.85\textwidth]{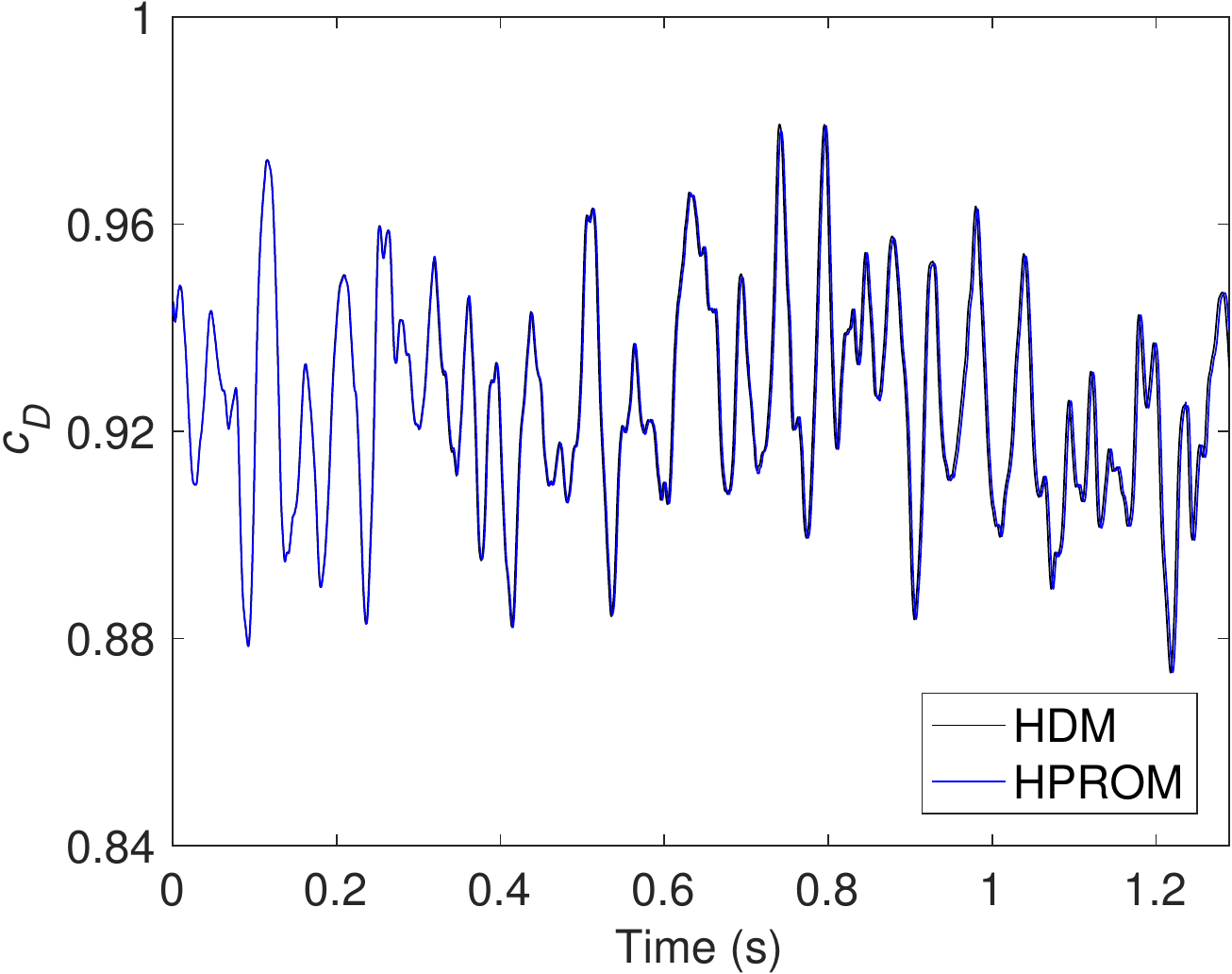}
	\caption{}
	\end{subfigure}%
	%\vglue 0.1 truein
	\caption{Turbulent flow past an F-16C/D aircraft: time-histories of the lift and drag coefficients computed using the HDM and ECSW-based LSPG HPROM.}
	\label{fig:f16drag}
\end{figure}

Figures \ref{fig:f16vort} and \ref{fig:f16mach} compare isosurfaces of vorticity magnitude and visualizations of the flow Mach number in the wake of the aircraft, respectively, computed at  
$t = T_f = 1.29$ \si{\second} using the results of the HDM-based and ECSW-based LSPG HPROM simulations. They suggest that the ECSW-based LSPG HPROM is capable of reproducing the flow features captured by
the HDM. Figure \ref{fig:f16drag}, which contrasts the time-histories of the lift and drag coefficients computed using the HDM with their counterparts computed using the ECSW-based 
LSPG HPROM, and Table \ref{tab:f16err}, which reports the relative errors computed using $\Delta s_{\mathbb{RE}}=5\times10^{-1}$ of the HPROM-based predictions of the lift and drag coefficients 
support this suggestion.

\begin{table}[h!]
	\small
	\centering
	\caption{Turbulent flow past an F-16C/D aircraft: computational accuracy of the ECSW-based LSPG HPROM.}
	%\vglue 0.1 truein
	\begin{tabular}{ccc}
		\toprule
		$N_c$ & $\mathbb{RE}_{c_D}$ ($\%$) & $\mathbb{RE}_{c_L}$ ($\%$) \\ \midrule
		$50$ & $0.489$ & $0.459$ \\ \bottomrule
	\end{tabular}
	%\vglue 0.1 truein
	\label{tab:f16err}
\end{table}

\begin{figure}[h!]
	\centering
	%\vglue 0.1 truein
	\includegraphics[width=0.3825\textwidth]{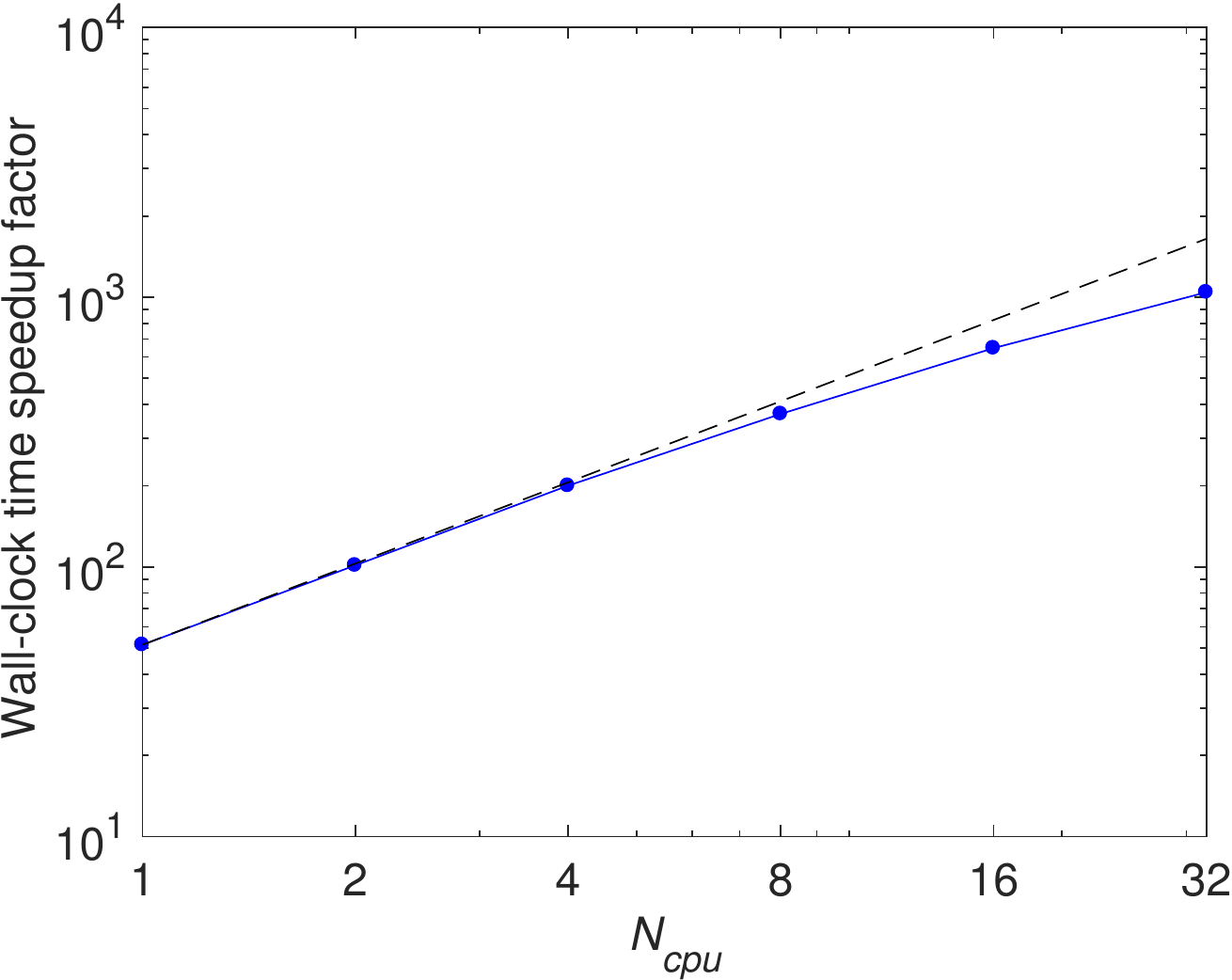}
	%\vglue 0.1 truein
	\caption{Turbulent flow past an F-16C/D aircraft: parallel scalability of the ECSW-based LSPG HPROM and delivered wall-clock time speedup factors.}
	\label{fig:f16cpu}
\end{figure}

For this problem, which features a CFD mesh with more than 25 million vertices (and thus more than 25 million FV dual cells), 
ECSW generates the reduced mesh in 30.6 \si{\minute} on 3,548 cores of the Frontera supercomputer: 14.6 \si{\minute} of these are 
elapsed in forming the convex optimization problem defined by (\ref{eqn:nnls}) and (\ref{eqn:eps}), and the other 16 \si{\minute} 
are consumed by the parallel NNLS algorithm \cite{chapman2017} for solving this 
problem. In order to highlight the fact that even very-low-dimensional computations can be parallelized to some degree, the HPROM-based simulation is repeated on $N_{cpu} = 1$, $2$, $4$, $8$,
$16$, and $32$ cores of a single node of the Frontera supercomputer. The obtained performance results are reported in Figure \ref{fig:f16cpu}, which highlights the parallel scalability of the
ECSW-based LSPG HPROM -- due to the concept of a reduced mesh described in Section \ref{sec:sampling} --  and Table \ref{tab:f16cpu}, which reports for $N_{cpu}= 32$ cores a wall-clock time speedup 
factor of three orders of magnitude and a CPU time speedup factor of five orders of magnitude. 

\begin{table}[h!]
	\small
	\centering
	\caption{Turbulent flow past an F-16C/D aircraft -- computational performance of the ECSW-based LSPG HPROM on $N_{cpu}=32$ cores: wall-clock execution time; and wall-clock time and CPU time 
	speedup factors.}
	%\vglue 0.1 truein
	\begin{tabular}{cccc}
		\toprule
		$N_c$ & Wall-clock time (s) & Wall-clock time speedup factor & CPU time speedup factor \\ \midrule
		$50$ & $ 346$ & $1.04\times10^3$ & $1.17\times10^5$ \\ \bottomrule
	\end{tabular}
	%\vglue 0.1 truein
	\label{tab:f16cpu}
\end{table}

\section{Conclusions}
\label{sec:conclusions}

In this paper, the ECSW hyperreduction method -- one of the most popular hyperreduction methods of the project-then-approximate type developed for accelerating Galerkin PMOR methods -- is extended to PG PMOR methods. Its computational framework is also generalized to cover not only FE spatial discretizations, but also FV 
and FD semi-discretization methods. Its scope is extended to PMOR methods based on local, piecewise-affine approximation subspaces designed for addressing the Kolmogorov $n$-width barrier
issue associated with many highly nonlinear problems such as high-speed, convection-dominated flow problems. For large-scale turbulent flow problems with $O(10^7)$ and $O(10^8)$ DOFs, the
offline part of the resulting ECSW method is shown to be not only computationally tractable, but also computationally fast. For such large-scale applications, the online part of the resulting ECSW
method for PG PMOR methods is also demonstrated to be robust, accurate, and most importantly, to enable wall-clock time and CPU time speedup factors of several orders of magnitude.

%\backmatter

\section*{Acknowledgments}
Sebastian Grimberg, Radek Tezaur, and Charbel Farhat acknowledge partial support by the Air Force Office of Scientific Research under grant FA9550-17-1-0182,
partial support by the Office of Naval Research under Grant N00014-17-1-2749, partial support by a research grant from the King Abdulaziz City for Science and Technology (KACST),
and partial support by The Boeing Company under Contract Sponsor Ref. 45047. All authors acknowledge the Texas Advanced Computing Center (TACC) at The University of Texas at Austin
for providing computing resources that have contributed to producing the numerical results reported in this paper. This document however does not necessarily reflect the position of
any of these institutions and therefore no official endorsement should be inferred.

\bibliography{main}
\bibliographystyle{elsarticle-num-names}

\end{document}